\documentclass[preprint,3p,12pt]{elsarticle}

\usepackage{amssymb}
\usepackage{mathrsfs}
\usepackage{amsfonts}

\usepackage{amsmath}
\usepackage{shortvrb,psfrag}
\usepackage{enumerate}
\usepackage{color}

\usepackage{bbding}
\usepackage{dcolumn}
\usepackage{graphicx}
\usepackage{psfrag}
\usepackage{wrapfig}
\usepackage{subfigure}
\usepackage{makeidx}
\usepackage{bm}
\usepackage{epsf}
\usepackage{epsfig}

\usepackage{footnote}

\newtheorem{proposition}{Proposition}[section]
\newtheorem{corollary}{Corollary}[section]

\newtheorem{remark}{Remark}[section]

\newtheorem{assumption}{Assumption}[section]

\makeatletter
\@addtoreset{figure}{section}

\makeatother

\makeatletter
\@addtoreset{table}{section}

\makeatother

\let\al=\alpha
\let\be=\beta
\let\si=\sigma

\let\f=\frac

\let\ga=\gamma
\let\La=\Lambda
\let\Ga=\Gamma
\let\la=\lambda
\let\De=\Delta

\let\ka=\kappa
\let\p=\partial

\def\no{\noindent}
\def\nonu{\nonumber}

\let\l=\left
\let\r=\right

\newcommand{\na}{\nabla}

\newcommand{\beq}{\begin{equation}}
\newcommand{\eeq}{\end{equation}}

\newcommand{\ben}{\begin{eqnarray}}
\newcommand{\een}{\end{eqnarray}}
\newcommand{\beno}{\begin{eqnarray*}}
\newcommand{\eeno}{\end{eqnarray*}}

\newcommand\dif{\mathrm{d}}     
\newcommand\diff{\,\dif}        

\renewcommand{\v}[1]{\textbf{\textit{#1}}}

\def\ef{\hphantom{MM}\hfill\llap{$\square$}\goodbreak}

\begin{document}


\title{ The generalized Riemann problems for hyperbolic balance laws: A unified formulation towards high order\tnoteref{ac}
}

\tnotetext[ac]{Jianzhen Qian is supported by Postdoctoral Science Foundation of China No. 2012M510366; Jiequan Li is supported by NSFC with Nos. 91130021 and 11031001; Shuanghu Wang is supported by NSFC with No. 91130021.
}

\author[iapcm]{Jianzhen Qian
}
\ead{qianjzmath@gmail.com}
\author [bnu]{Jiequan Li\corref{cor1}}
\ead{jiequan@bnu.edu.cn}
\author [iapcm,iapcm2]{Shuanghu Wang}
\ead{wang\_shuanghu@iapcm.ac.cn}

\address[iapcm]{Institute of Applied Physics and Computational Mathematics, Beijing 100088,
China}
\cortext[cor1]{Corresponding author}
\address[bnu]{ School of Mathematical Science, Beijing Normal University, Beijing, 100875, P. R. China }
\address[iapcm2]{Key Laboratory of Computational Physics,
Institute of Applied Physics and Computational Mathematics, Beijing
100088, China}

\begin{abstract}
The Generalized Riemann Problems (GRP) for nonlinear hyperbolic systems of balance laws
in one space dimension are now well-known and can be formulated as follows: Given initial-data which are smooth on two sides of a discontinuity, determine the time evolution of the solution near the discontinuity.
In particular, the GRP of $(k+1)$th order high-resolution is based on an
analytical evaluation of the time derivative up to $k$th order, which turns out to be dependent only on the spatial derivatives up to $k$th order.
While the classical Riemann problem serves as a primary ``building block" in the construction of many numerical
schemes (most notably the Godunov scheme), the analytic study of GRP will lead to an array of ``GRP
schemes'', which extend the Godunov scheme. Currently there are extensive studies on the second-order GRP scheme, which proves to be robust and is capable of resolving complex
 multidimensional fluid dynamic problems [M. Ben-Artzi and J. Falcovitz, ``Generalized Riemann Problems in
 Computational Fluid Dynamics", Cambridge University Press, 2003]. More general formulation of the second-order GRP solver can
 be found in [Numer. Math. (2007) 106:369-425], but still confined with a class of ``weakly coupled systems". In this paper, we provide
a unified approach for solving the GRP in the general context of hyperbolic balance laws, {\it without weakly coupled constraint}, towards high order accuracy.
 The derivation of the second-order GRP solver is more concise compared to those in previous works and  the third-order GRP (or quadratic GRP) is resolved for the first time. The latter is shown to be necessary through numerical experiments with strong discontinuities. Our method relies heavily on the new treatment of the rarefaction wave. Indeed, as a main technical step, the ``propagation of singularities" argument for the rarefaction fan, is simplified by deriving the L(Q)-equations,  an ODE system for the ``evolution" of the ``characteristic derivatives" in $x$-$t$ space for generalized Riemann invariants, with aid of the generalized characteristic coordinates. The case of a sonic point is incorporated into a general treatment. The accuracy of the derived GRP solvers are justified and numerical examples are presented for the performance of the resulting scheme.

%

\end{abstract}



\maketitle

\section{Introduction}
\label{sec:1}
\setcounter{equation}{0}
\no
In this paper we consider the {\it generalized Riemann problem} (GRP) for hyperbolic balance laws
\begin{equation}\label{eq:cl}
\f{\p U}{\p t}+\f{\p F(U)}{\p x}=H(x,U),
\end{equation}
where $U=(u_1,\cdots,u_m)$ is the unknown variable with $F=(f_1,\cdots,f_m)$ being the flux functions, and $H(x,U)$ is a
source term resulting from geometrical or physical effects, $x$ is
the spatial variable and $t$ is the time variable. In this study,
 we will concentrate on the numerical aspect of (\ref{eq:cl}), rather
 than important theoretical issues such as well-posedness and solution structures.

In the development of numerical techniques approximating solutions of (\ref{eq:cl}), the finite volume
scheme plays absolutely indispensable role,
wherein one of most crucial ingredients is the construction of numerical fluxes and it boils down to the resolution
 of associated (generalized) Riemann problems
at each computational cell interface. Specifically, we denote by $I_j=[x_{j-1/2}, x_{j+1/2}]$, $\De x=x_{j+1/2}-x_{j-1/2}$, the computational cell numbered $j$,
and by $\{t_n\}_{n=0}^{\infty}$ the sequence of discretized time levels, $\Delta t=t_{n+1}-t_{n}$. The finite volume
scheme is then constructed by integrating the governing equations (\ref{eq:cl}) both in space and time over the control
volume $I_j\times [t_{n+1},t_n]$, yielding
\begin{equation}\label{eq:fv}
U_j^{n+1}=U_j^{n}-\f{\De t}{\De x}\big(F_{j+1/2}^n-F_{j-1/2}^n\big) + \De t H_j^n,
\end{equation}
where
\begin{equation}\label{eq:U_ave}
U_j^{n}=\f{1}{\De x}\int_{x_{j-1/2}}^{x_{j+1/2}} U(x,t_n) \diff x
\end{equation}
is the average of $U(x,t^n)$ over the cell $I_j$.
The remaining terms in (\ref{eq:fv}), $F_{j+1/2}^n$ and $H_j^n$ , are
the temporal average of $F(U(x,t))$ along the interface $x=x_{j+1/2}$ and the space-time integral average of the source $H(x,U)$, i.e.,
\begin{align}\label{eq:F_ave}
&F_{j+1/2}^n=\f{1}{\De t}\int_{t_n}^{t_{n+1}} F(U(x_{j+\f{1}{1}},t)) \diff t,\\[2mm]
\label{eq:S_ave}
&H_j^n=\f{1}{\De t \De x}\int_{t_{n}}^{t_{n+1}}\int_{x_{j-1/2}}^{x_{j+1/2}} H(x,U)\diff x\diff t.
\end{align}

A numerical scheme is obtained if one can supply suitable approximation for $F_{j+1/2}^n$ and $H_j^n$ with given data $U(x,t_n)$ at $t=t_n$.
Formally, for the Godunov-type schemes, this usually consists the following three procedures.
\begin{enumerate}
    \item[A] Data reconstruction: Based on the cell average values $U_j^n$, reconstruct the initial data $U(x,t_n)$ as piece-wise smooth distribution, being constant or polynomial in each cell $I_j$.
    \item[B] Solution evolution: Solve the (generalized) Riemann problem
at each cell interface $x=x_{j+1/2}$ to evolve the solution.
    \item[C] Numerical approximation: Take the numerical integration in (\ref{eq:F_ave}) and (\ref{eq:S_ave}) to get $F_{j+1/2}^n$ and $S_j^n$ under suitable CFL condition.
\end{enumerate}
\begin{figure}
  \begin{center}
    \includegraphics[height=1.6in, width=3.6in, trim=0 0 0 0, clip]{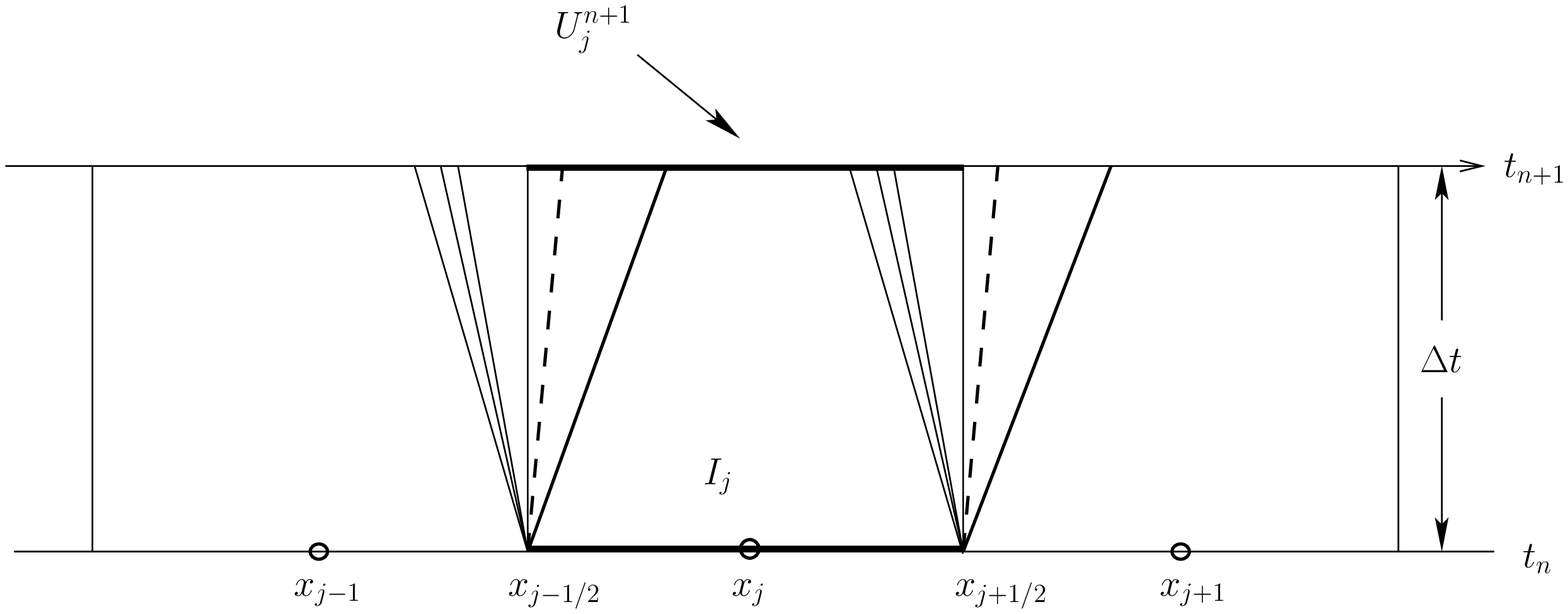}
  \end{center}
  \label{fig:grid}
\end{figure}

For the notable Godunov scheme \cite{god} and higher order schemes using Riemann solvers such as MUSCL \cite{muscl1,muscl2} and TVD \cite{tvd} schemes,
 the classical Riemann problem is solved in each cell interface to evolve the solution from $t^n$ to $t^{n+1}$.
The corresponding initial data are taken as interface limit values of $U(t_n,x)$  in the neighboring cells.
As an extension, the GRP scheme assumes
the piecewise smooth initial data and evolves solutions by analytically solving the generalized Riemann problem at each cell interface with at least second order accuracy.
Currently, the second-order GRP scheme has already been exploited and put into use for several compressible fluid models \cite{ben1,ben2,ben3,ben4,li1,li2,li-sw,yang-tang1,yang-tang2}.
 Let us now outline a standard process for its implementation.
 Assume the data at time $t=t^n$ is piece-wise linear with slope $\si_j^n$, i.e. on $I_j$ we have
\begin{equation}\label{eq:u_g2}
U(x,t^n)=U_j^n+\si_j^n(x-x_j), \quad x\in(x_{j-1/2}, x_{j+1/2}).
\end{equation}
Then the Godnuov-type scheme of second order takes the form
\begin{equation}\label{eq:grp2}
U_j^{n+1}=U_j^{n}-\f{\De t}{\De x}\big(F_{j+1/2}^{n+1/2}-F_{j-1/2}^{n+1/2}\big) + \f{\De t}{2}\big(H_{j+1/2}^{n+1/2}+H_{j-1/2}^{n+1/2}\big),
\end{equation}
where the following notations are used
\begin{equation}\label{eq:grp2_nota}
F_{j+1/2}^{n+{1/2}}=F\big(U_{j+1/2}^{n+1/2}\big),\quad H_{j+1/2}^{n+1/2}=H\big(x_{j+1/2},U_{j+1/2}^{n+1/2}\big),
\end{equation}
 and $U_{j+1/2}^{n+1/2}$ is the mid-point value or the average of $U(x_{j+1/2},t)$ over time interval $[t_n,t_{n+1}]$. For simplicity, the source term is currently discretized with an interface method, which is the trapezoidal rule in space and the mid-point rule in time \cite{ben4,jin} in order to keep second order accuracy. The central issue is how to obtain the mid-point value $U_{j+1/2}^{n+1/2}$, which is formally approximated by the Taylor expansion (ignoring the higher order terms)
\begin{equation}\label{eq:mid_point}
U_{j+1/2}^{n+{1/2}}\cong U_{j+1/2}^n+\f{\De t}{2}\left(\f{\p U}{\p t}\right)_{j+1/2}^n,
\end{equation}
where
\begin{equation}\label{eq:mid_p_lim}
U_{j+1/2}^{n}=\lim_{t\rightarrow t_n+0} U(x_{j+1/2},t),\quad \left(
\f{\p U}{\p t}\r)_{j+1/2}^{n}=\lim_{t\rightarrow t_n+0} \f{\p U}{\p t}(x_{j+1/2},t).
\end{equation}
The value $U_{j+1/2}^{n}$ is obtained by solving the associated Riemann problem for the homogeneous
hyperbolic conservation laws as used in the (first order) Godunov scheme \cite{god}.
 The main ingredient lies upon the calculation of the {\it instantaneous time derivative} $(\f{\p U}{\p t})_{j+1/2}^n$.
Even in the Godunov scheme, the time derivative $(\f{\p U}{\p t})_{j+1/2}^n$ should be properly treated once the source term is present, which makes the solution evolve non-uniformly.

For the solution $U$ being smooth near the grid point $(x_{j+1/2}, t_n)$, it follows directly from (\ref{eq:cl}) that
\begin{equation}\label{eq:U_t_dir}
\l(\f{\p U}{\p t}\r)_{j+1/2}^{n}=-\f{\p F}{\p U}\l(U_{j+1/2}^{n}\r) \l(\f{\p U}{\p x}\r)_{j+1/2}^{n}+H\big(x_{j+1/2}, U_{j+1/2}^n\big).
\end{equation}
However, for the generalized Riemann problem including singularity at grid point $(x_{j+1/2}, t_n)$, (\ref{eq:U_t_dir}) is no longer valid, even for scalar cases, because there exists nonlinear waves (rarefaction
waves or discontinuities) issuing from the singularity point $(x_{j+1/2}, t_n)$.
 Indeed, thinking of the initial data (\ref{eq:u_g2}) with non-zero slopes as a perturbation of piecewise constant Riemann initial data and the source term $S(x,U)$ as a perturbation of the homogenous system of equations, the GRP solution is a perturbation of that of the {\it associated Riemann problem} at least in the neighborhood of the singularity point. It turns out that the GRP solution consists of, for a short time following the ``disintegration" of initial discontinuity, the curvilinear rarefaction wave and the discontinuities (contact discontinuity or shock wave) with time varying speed \cite[Chap. 5]{ben4}.

The solution $U$ together with its derivatives may undergo a jump discontinuity across each wave. Hence, in order to solve the generalized Riemann problem, it requires one to explore
 the mode of the discontinuity for the derivatives coming along with each wave, which is in fact described by a set of linear algebraic equations.
This bears an analogy to the resolution of classical Riemann problem, which involves exploring the relation, usually described by a one parameter curve,
between the two states of $U$ connected by each wave. Indeed, the treatment for capturing the discontinuities of the derivatives
across each type of waves can be sketched out as follows.
\begin{enumerate}
\item[A] Since the generalized Riemann invariants are transported in the transversal direction of the rarefaction fan, it is natural to use them for studying the variation of the derivatives across the rarefaction wave. Actually, the directional (emanating characteristic direction) derivatives of the generalized Riemann invariants are determined by their values on either side of the wave.
\item[B] The generalized Riemann invariants, which remain continuous across corresponding contact discontinuities, are differentiated in the direction of the discontinuity (characteristic).
\item[C] For the shock wave, the identities implied by the Rankine-Hugoniot conditions are differentiated along the shock trajectories.
\end{enumerate}
As indicated in the previous works, the most technical step lies on the treatment for rarefaction fan, which relies on the analysis in term of ``characteristic coordinates''.

The methodology for resolving the generalized Riemann problem is originated in \cite{ben1,ben2,ben3}, wherein the original GRP is designed for the compressible fluid flows with two related Lagrangian and Eulerian versions.
See also the recent textbook \cite{ben4} for detailed discussions. The Eulerian version is always derived by using the Lagrangian case.
 The transformation is quite delicate, particularly for sonic cases, because it becomes singular at sonic points. The direct Eulerian version, more flexible for applications, is developed
recently in the context of shallow water equations \cite{li-sw}, planar compressible fluid and the compressible fluid
flows \cite{li1,li2}. The approach for solving GRP therein, being ready to handle any strict hyperbolic
 system endowed with a coordinate system of Riemann invariants (in particular, the two equations system),
is extended to handle more general weakly coupled systems (in the sense of \cite[Def. 21]{li2}) having only a
 ``partial set" of Riemann invariants.
The common point of the above systems is that the {\it generalized Riemann invariants} (GRI) are coupled in a manner that
enables a ``diagonalized" treatment.
Although many physical systems, including the compressible fluids flow system, belong to such a class of systems, the existing methods for deriving a second-order GRP solver  are rather complicated, which prevents it from practical use in many ways.
For example, the treatment of rarefaction relies heavily on the explicit formulation of the
{\it Asymptotic Characteristic Coordinate} (ACC), which depends on the EOS (equation of state) of the fluid in turn and is sometimes hard to derive. Besides the ACC is not easy to be written out explicitly for higher order GRP solvers that are particularly useful in capturing the propagation of entropy wave \cite{shu2} (see also Fig. \ref{fig:shock-density}). Other closely related efforts can be found in \cite{asym1,asym2} using  the approach of asymptotic analysis  for the resolution of generalized Riemann problems, and  in \cite{eno,tt,ct1} (and the references therein) for approximate Godunov-type high order solvers. The solvers in \cite{tt,ct1} corresponds to the acoustic case and they fails for resolving strong discontinuities, even with very high order accuracy (this point is confirmed through a numerical experiment Fig. \ref{fig:dp100-1000}).
Hence  it is absolutely necessary to develop the high order (at least third order) GRP scheme by resolving nonlinear wave patterns  each computational grid point analytically, in addition to provide an acoustic approximation as the jump there is weak.
\vspace{0.2cm}

Therefore we present a unified approach in this paper, still direct Eulerian, to resolve
 the GRP for general systems of hyperbolic balance laws (\ref{eq:cl}). The weakly coupled constraint is not required here.
The solver for second-order GRP (linear GRP) as well as third-order (quadratic GRP) are derived.
This paper provides a simplified treatment for the main step, resolution of the rarefaction fan.
Indeed, it is carried out by first deriving the system of transport equations for the generalized Riemann invariants and then
deriving the ``evolution" equations, labeled as the L(Q)-equations, for their characteristic derivatives in $x$-$t$ space with aid of the {\it generalized characteristic coordinate} (GCC).
This is based on the following observations. These characteristic derivatives, and hence the resulting ``evolution" equations for them, are independent of the choice of the auxiliary GCC.
Thus the explicit expression of the GCC is not required.
 More importantly, only in the emanating characteristic direction do the derivatives of the GRI (the solution $U$) of any order exist and remain continuous across rarefaction fans.
 No additional assumption is required for the regularity of the solution $U$.
Indeed, the above observations are the reasons why we can derive the third-order (or higher order) GRP solver without many difficulties.
By referring to Section \ref{sec:R-DW} for the resolution of the contact discontinuity and the shock wave, the spatial derivatives of the solution,
 from which the instantaneous time derivatives follow, are obtained by solving a simple system of linear equations in the intermediate regions of the waves. The case of sonic point is handled by supplementing an additional freedom using the differential relation of $U$ along the emanating characteristic direction. A special case frequently occur during the numerical application
 of the GRP scheme is the acoustic case: the initial values of $U$ are continuous at the singularity point. This case
is comparatively easy to handle and requires less computation cost.

 Although this paper focus on exploring solvers for the second-order {\it linear GRP} and the third-order {\it quadratic GRP},
higher order GRP solvers can be derived with the same methodology and a multidimensional extension can be pursued in a forthcoming work \cite{Li-Qian}.
The resulting GRP solvers consist of two steps: (i) the classical Riemann solver; (ii) calculation of instantaneous time derivatives of $U$. As indicated by the solvers, Step (ii) can be straightforward once the full Riemann solution is obtained. Besides, in both steps, only the limiting values of $U$ and its spatial derivatives at two side of the singularity are used, and the resulting linear (resp. quadratic) GRP solver leads to second (resp. third) order accuracy in time approximation to $U$ regardless its initial distribution.

This paper is arranged as follows. In Section 2, a basic setup for the system and the GRP are presented. The resolution of rarefaction wave and discontinuity waves, including the contact discontinuity and shock wave, are detailed in Sections 3 and 4, respectively. We conclude the resolution of GRP in Section 5 and the acoustic approximation in Section 6. As an application example, in Section 7, we derive the GRP solvers for compressible variable duct flow system and show the solvers' accuracy by several tests. Finally, in Section 8, the GRP solvers are used to construct one-step high order numerical scheme and a few 1-D numerical test cases are presented.

\vspace{1cm}

\section{Basic setup for the system and the GRP}
\label{sec:basic}
\setcounter{equation}{0}
As a basic setup, we assume (\ref{eq:cl}) is hyperbolic in the sense that the Jacobian $A(U)=\f{\p F(U)}{\p U}$ of $F(U)$ has $m$ eigenvalues
\begin{equation}\label{eq:eig}
\la_1\le \la_2\le \cdots\le \la_m.
\end{equation}
The set of left (right) eigenvectors $L_k$ ($R_k$) (associated with $\la_k$, $k=1,\cdots,m$) are linearly independent.
The {\it $k$th characteristic field} $\la_k$ can be either genuinely nonlinear in the sense of $\na_{U}\la_k \cdot R_k\neq 0$, or linearly degenerate $\na_{U}\la_k \cdot R_k= 0$.


Now let us state the generalized Riemann problem. It is defined as the initial-value problem for system (\ref{eq:cl}), subject to the initial data
\begin{equation}\label{eq:clini}
U(x,0)=\left\{
\begin{array}{ll}
P_+(x) &\text{if}\ x<0,\\
P_-(x) &\text{if}\ x>0,
\end{array}
\right.
\end{equation}
where $P_\pm(x)$ are vectors, whose components are the smooth functions.
As illustrated in Section \ref{sec:1}, the initial structure of the solution is determined by the {\it associated Riemann problem}:
\begin{equation}\label{eq:rp}
\left\{
\begin{array}{l}
\f{\p U^A}{\p t}+\f{\p F(U^A)}{\p x}=0,\\[2mm]
U^A(x,0)= U_{\pm}, \quad \pm x>0,
\end{array}
\right.
\end{equation}
where $U_{\pm}$ are the limiting values of $P_{\pm}(x)$ at $x=0$, i.e. $U_{\pm}=P_{\pm}(0^\pm)$.
We call the solution of (\ref{eq:rp}) the {\it associated Riemann solution} of (\ref{eq:cl}) and (\ref{eq:clini}).

\begin{assumption}
 The Riemann problem (\ref{eq:rp})
is uniquely solvable, and the solution to (\ref{eq:rp}) consists of $m$ waves $\Gamma_1,\Gamma_2,\cdots,\Gamma_m$. The wave $\Gamma_k$($1\le k\le m$) is an admissible shock, a contact discontinuity, or a rarefaction wave associated with the $k$th characteristic field $\la_k$.
\end{assumption}

Note that the above assumption does not mean we are confined with strict hyperbolic systems that endowed with distinct eigenvalues.\\

Denote by $R^A(x/t,U_-,U_+)$ the Riemann solution of (\ref{eq:rp}). Then we have the following proposition.

\begin{proposition}\label{pro:1}
Let $U(x,t)$ be the solution to the generalized Riemann problem
(\ref{eq:cl}) and (\ref{eq:clini}). Then for every fixed direction $\theta=x/t$,
\begin{equation}
\lim_{t\rightarrow 0}U(\theta t,t)=R^A(\theta,U_-,U_+).
\end{equation}
This implies that the wave configuration for the generalized Riemann problem(\ref{eq:cl}) and (\ref{eq:clini}) is the same as that for the associated Riemann problem (\ref{eq:rp}) around the singularity $(x,t)=(0,0^+)$.
\end{proposition}

\begin{figure}
  \begin{center}
    \subfigure[]{\label{fig:grp-a}%
    \includegraphics[height=1.9in, width=2.8in, trim=0 0 0 0, clip]{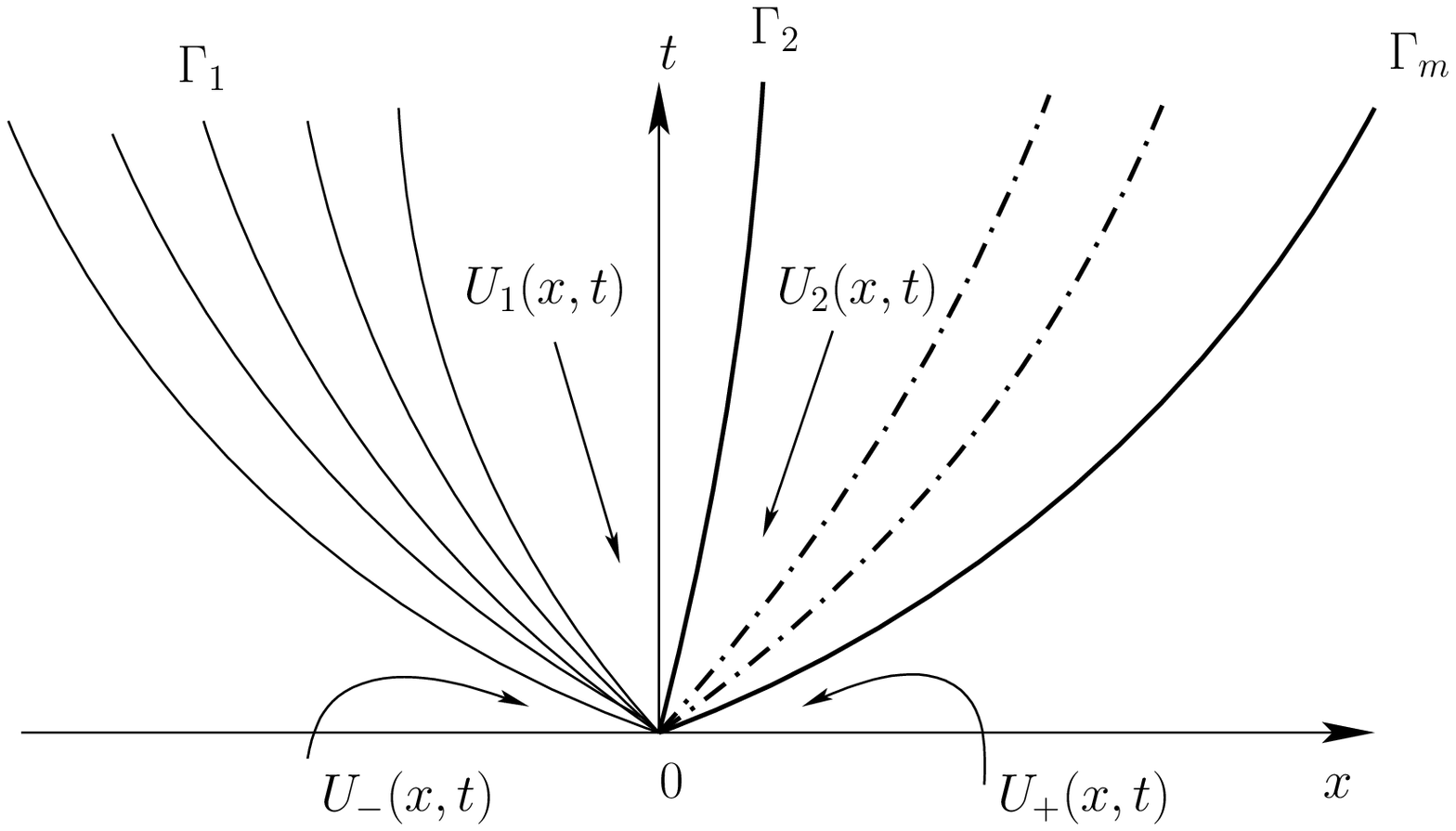} }
    \subfigure[]{\label{fig:grp-b}%
    \includegraphics[height=1.82in, width=2.7in, trim=0 -10 0 0, clip]{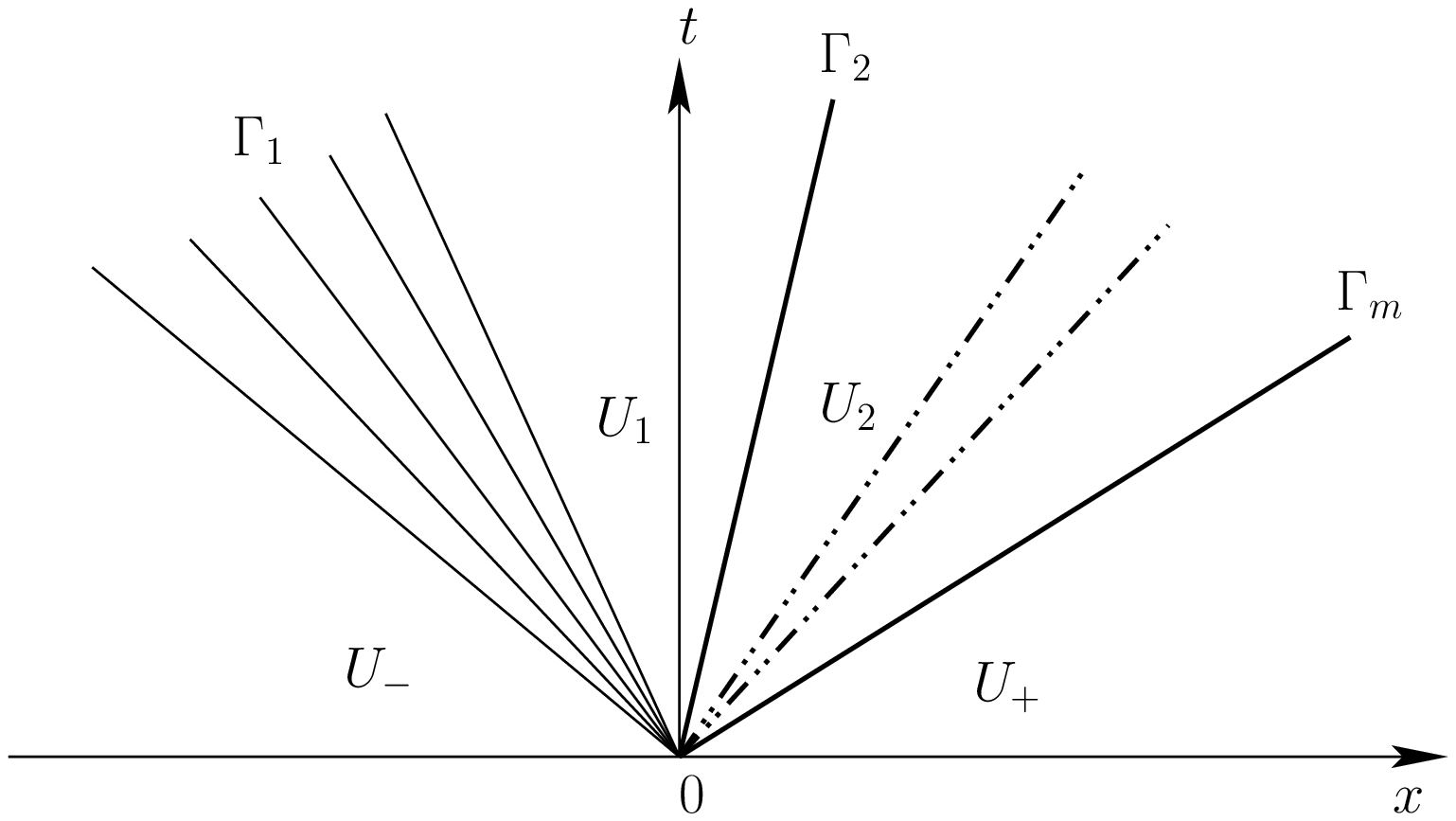} }
  \end{center}
\caption{ Wave configurations: (a) Wave patterns for the GRP with initial data $U(0,x)=P_-(x)$ for $x<0$ and $U(0,x)=P_+(x)$ for $x>0$, $U_\pm=P_\pm(0^\pm)$. (b) Wave patterns for the associated Riemann problem.
}
  \label{fig:grp-waves}
\end{figure}

Proposition \ref{pro:1} is illustrated schematically in Fig.\ref{fig:grp-waves}. The solution of (\ref{eq:rp}) is self-similar, and hence the waves are centered. Correspondingly, the waves for (\ref{eq:cl}) are curved (See \cite{ben4} for more detailed descriptions).

We emphasize that the solution $U$ is smooth in the intermediate regions of these waves and along each emanating characteristic curve in the rarefaction fan (up to the singularity $(0,0^+)$).
To approximate $U$ along $t$-axis with $k$th order accuracy, we can use the Taylor expansion
\begin{equation}\label{eq:tal}
U(x=0,t)=U(0,0^+)+\sum_{\ell=1}^{k}\f{1}{\ell!}\f{\p^{\ell} U}{\p t^\ell}(0,0^+) (t^{k}) +\mathcal{O}(t^{k+1}).
\end{equation}
A solver of the GRP is actually that of evaluating the instantaneous time derivatives
\begin{equation}\label{eq:td}
\f{\p^{\ell} U}{\p t^\ell} (0,0^+)=\lim_{t\rightarrow 0}\f{\p^{\ell} U}{\p t^\ell} (0,t), \quad t>0.
\end{equation}
For convenience,  we label the problem of evaluating (\ref{eq:td})
 with $\ell=1$ (resp. $\ell=1,2$) as the {\it linear GRP} (resp. {\it quadratic GRP}), or {\it LGRP} (reps. {\it QGRP}) for short.
As mentioned in the introduction, this paper concentrates on QGRP.

\section{Resolution of curved rarefaction waves}
\label{rarefaction}
\setcounter{equation}{0}
As pointed out earlier, the main feature of the GRP is the resolution of rarefaction waves and the main ingredients are the Riemann invariants and characteristic coordinates.
Let us consider to first derive the set of transport equations for the (generalized) Riemann invariants in a general setting. For this purpose, we rewrite (\ref{eq:cl}) as a nonconservative form
\begin{equation}\label{eq:ncl}
\f{\p U}{\p t}+A(U)\f{\p U}{\p x}=H(x,U),
\end{equation}
by recalling $A(U)=\p F(U)/\p U$.
Multiplying (\ref{eq:ncl}) by $L=(L_1,\cdots,L_m)$ from the left, it follows that
\begin{equation}\label{eq:lncl}
L\f{\p U}{\p t}+\La L \f{\p U}{\p x}=L H(x,U),
\end{equation}
where $\La=\textrm{diag}(\la_1,\cdots,\la_m)$.
If there exists a set of variables $\v{w}=(w_1,\cdots,w_m)$ satisfying $\f{\p w_k}{\p U} \parallel L_k$ ($k=1,\cdots,m$), then (\ref{eq:lncl}) is equivalent to
\begin{equation}\label{eq:lnclw}
\f{\p \v{w}}{\p t}+\La \f{\p \v{w}}{\p x}=L H(x,U).
\end{equation}
Indeed, $\v{w}$ is a complete set of {\it Riemann invariants}. Unfortunately, most of the systems (\ref{eq:cl}) with $m>2$, including the full system of compressible Euler equations, do not admit such a set of Riemann invariants.
We thus turn to exploring the {\it generalized Riemann invariants} (GRI).

\subsection{The generalized Riemann invariants (GRI)}
Let $\v{w}=(w_1,\cdots,w_{m-1})$ be the generalized Riemann invariants of the $k$th characteristic field. By recalling the definition of GRI \cite{toro}, we have
$$\na_Uw_\ell\cdot R_k=0,\quad \ell=1,\cdots,m-1.
$$
Hence, there exists an invertible $(m-1) \times (m-1)$ matrix $K$, such that
\begin{equation}\label{eq:ms}
K\na_U\v{w}=(L_1,\cdots,L_{k-1},L_{k+1},\cdots,L_m)^T =: L^{(k)}.
\end{equation}
Multiply (\ref{eq:ncl}) by $K\na_U\v{w}$ from the left yields the following proposition.
\begin{proposition}\label{pro:2}
Let $\v{w}=(w_1,\cdots,w_{m-1})$ be the GRI of the $k$th characteristic field. Then in any smooth region of $U$ there holds
\begin{equation}\label{eq:w}
\f{\p \v{w}}{\p t}+B^{(k)}(U) \f{\p \v{w}}{\p x}=L^{(k)}H(x,U),
\end{equation}
where
\begin{equation}\label{eq:coeff_B}
B^{(k)}(U)=K^{-1}\La^{(k)}K,\quad \La^{(k)}=\textrm{diag} (\la_1,\cdots,\la_{k-1},\la_{k+1},\cdots, \la_{m})^T,
\end{equation}
and $K$ is determined by (\ref{eq:ms}).
\end{proposition}

Roughly speaking, Proposition \ref{pro:2} implies that the generalized Riemann invariants of the $k$th characteristic field are transported along the direction different from $\la_k$.  The following useful corollary is straightforward from Proposition \ref{pro:2}  for the resolution of rarefaction wave. See Remark \ref{rem:2} (ii) below.
\begin{corollary}\label{cor:1}
 Let $\Gamma_k$ be a characteristic curve associated with $\la_k$. If $U$ is continuous and piecewise smooth with $\Gamma_k$ being a weak discontinuity curve, then ${\p \v{w}}/{\p t}$ and ${\p \v{w}}/{\p x}$ in Proposition \ref{pro:2} remain continuous across $\Gamma_k$.
\end{corollary}

\subsection{Generalized characteristic coordinates (GCC)}
As mentioned in the introduction, the characteristic coordinates, defined as the integral curves of the characteristic equations, play an important role in the resolution of rarefaction waves. In the region of a rarefaction fan, they work similarly to the usual polar coordinates to single out  singularities.

\begin{figure}
  \begin{center}
    \includegraphics[height=2.2in, width=2.8in, trim=0 0 0 0, clip]{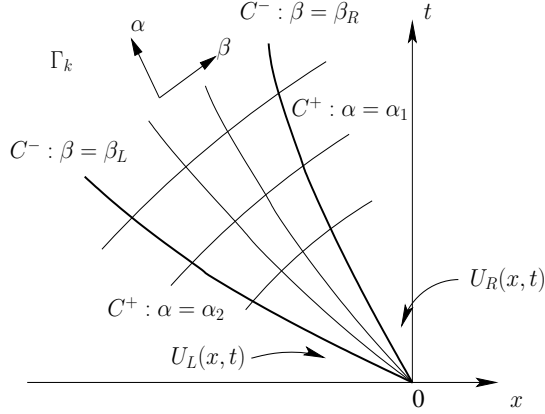}
  \end{center}
\caption{Generalized characteristic coordinates in a rarefaction fan $\Ga_k$.}
  \label{fig:axis}
\end{figure}

Assume that $\Gamma_k$ is a rarefaction wave associated with $\la_k$ and denote by $U_L(x,t)$ (resp. $U_R(x,t)$) the state $U$
 on its left (resp. right) side. See Fig. \ref{fig:axis}. To simplify notations, we
write $\la$ and $B$ below for $\la_k$ and $B^{(k)}$ in Proposition \ref{pro:2}, respectively, for a fixed $k$.

Let $C^-$: $\beta(x,t)=\beta$ and $C^+$: $\al(x,t)=\al$,
 $\be\in[\be_L,\be_R]$, $-\infty\le \al<0$
be the integral curves of the following equations, respectively,
\begin{equation}\label{eq:ceq}
\f{\dif x}{\dif t}=\la,\quad \f{\dif x}{\dif t}=\mu.
\end{equation}
Here, different from the previous works \cite{ben1,ben2,ben3,ben4,li1,li2},  $\mu$ in (\ref{eq:ceq}) is not required to be an eigenvalue of $A(U)$. In fact, it can be the slope of any family of transversal curves different from $\la$. For example, $\mu=-\f{t}{x}$. The GCC used here is a relaxed version of the afore-mentioned ACC.
Moreover, $\be$ and $\al$ are denoted as follows:
$\be$ is the initial value of the slope $\la$ at the singularity $(x,t)=(0,0)$ and $\al$ for the transversal characteristic curves is the $x$-coordinates of the intersection point with the leading $\be$-curve: $\be=\be_L$.

The coordinates $(x,t)$ in the ``triangle" sector of the {\it centered rarefaction wave} shown in Fig. \ref{fig:axis} can be expressed in terms of $\alpha$ and $\beta$,
\beq\label{eq:coor}
x=x(\al,\be),\quad t=t(\al,\be),
\eeq
which satisfy
\beq\label{eq:coorl}
\f{\p x}{\p \alpha}=\la\f{\p t}{\p \alpha},\quad \f{\p x}{\p \be}=\mu\f{\p t}{\p \be}.
\eeq
Denote
\beq
D_\la=\f{\p}{\p t}+\la\f{\p}{\p x}, \quad D_\mu=\f{\p}{\p t}+\mu \f{\p}{\p x}.
\eeq
Then we have
\beq  
\f{\p}{\p \al}=\f{\p t}{\p \al} D_{\la}, \quad \f{\p}{\p \be}=\f{\p t}{\p \be} D_{\mu}.
\eeq
In particular, as $\al=0$, we have
\beq\label{eq:al0}
\f{\p \la}{\p \be}(0,\be)=1,\quad \f{\p t}{\p \be}(0,\be)=0,\quad \be_L\le \be \le \be_R.
\eeq
Here we remind that, as a basic assumption in Section \ref{sec:basic}, the solution $U$ is smooth along each characteristic curve $C^-:\be=\bar{\be}$ inside a rarefaction fan up to the singularity and $D^\ell_\la U$, for any $\ell\geq 1$, takes finite value at $\alpha=0$.

It follows, by differentiating the first equation in (\ref{eq:coorl}) with respect to $\be$, the second with respect to $\al$ and then subtracting, that the function $t=t(\al,\be)$ satisfies
\begin{align}\label{eq:al1}
(\mu-\la)\f{\p^2 t}{\p \al\p \be}=\f{\p \la}{\p \be}\f{\p t}{\p \al}-\f{\p \mu}{\p \al}\f{\p t}{\p \be}.
\end{align}
Setting $\al=0$ and using (\ref{eq:al0}), one obtains
\begin{align}\label{eq:gt1}
\f{\p}{\p \be}\left[\f{\p t}{\p \al}(0,\be)\right]=\f{1}{\mu-\la}\f{\p t}{\p \al}(0,\be).
\end{align}

We continue to make differentiation of (\ref{eq:al1}) with respect to $\al$ to obtain
\begin{align}\label{eq:al2}
(\mu-\la)\f{\p^3 t}{\p \al^2 \p \be}
=&- \f{\p }{\p \al}(\mu-\la) \f{\p^2 t}{\p \al\p \be}+\f{\p^2 \la}{\p\al \p\be}\f{\p t}{\p \al} + \f{\p \la}{\p \be}\f{\p^2 t}{\p \al^2} \nonu\\
&\hspace{1.5cm}-\f{\p^2 \mu}{\p\al^2}\f{\p t}{\p \be}-\f{\p \mu}{\p \al}\f{\p^2 t}{\p \al\p \be}.
\end{align}
Recalling (\ref{eq:al0}) and  (\ref{eq:gt1}) as well as noticing
 \begin{align*}
\f{\p^2}{\p \al \p \be}\la=\f{\p}{\p \be}\l(\f{\p t}{\p \al} D_\la\la\r) =
\f{\p^2 t}{\p \al\p \be} D_\la\la+\f{\p t}{\p \al} \f{\p}{\p \be}\l( D_\la \la \r),
\end{align*}
we take $\al=0$ to obtain
\begin{align}\label{eq:gt2}
\f{\p}{\p \be}\l[\f{\p^2 t}{\p \al^2}(0,\be)\r]
=&
\f{1}{\mu-\la}\f{\p ^2 t}{\p\al^2}(0,\beta)+\f{2 D_\la}{\mu-\la}\l( \f{\p t}{\p\al} \r)^2(0,\beta)\nonu\\
&+\f{1}{\mu-\la} \f{\p}{\p\be}\l(D_\la \la\r) \l(\f{\p t}{\p \al}\r)^2(0,\beta).
\end{align}

The equations \eqref{eq:gt1} and \eqref{eq:gt2}, for ${\p t}/{\p \al}(0,\be)$ and ${\p^2 t}/{\p \al^2}(0,\be)$ repectively, are crucial for deriving the $L(Q)$-equations in next subsection.

\subsection{The $L(Q)$-equations}

In this subsection, we shall derive the linear differential ordinary equations for $D_\la \v{w}(0,\be)$ and $D_\la^2 \v{w}(0,\be) $ with respect to $\be$, namely, the {\it $L$-equations} and {\it $Q$-equations}, respectively, of the GRI.
Precisely, we have the following proposition.
\begin{proposition}[$L(Q)$-equations]\label{pro:LQ}
The $\v{w}$ in Proposition \ref{pro:2} satisfies the {\it $L$-equations}:
\begin{equation}\label{eq:L}
\f{\p}{\p \be} \big[D_\la\v{w}(0,\be) \big]=(\la I-B)^{-1}\big( D_\la \v{w}-L^{(k)}H\big),
\end{equation}
 and the {\it $Q$-equations}:
\begin{align}
\f{\p }{\p \be} \big[D_\la^2 \v{w}(0,\be) \big] =&2 (\la I-B)^{-1} D_\la^2 \v{w} +2D_{\la} (\la I-B)^{-1}D_{\la} \v{w} -2D_{\la}\big[ (\la I-B)^{-1} L^{(k)}H \big] \nonu\\
&\quad
+ \f{\p}{\p \be} (D_\la\la) \big[ (\la I-B)^{-1}D_{\la} \v{w}-(\la I-B)^{-1} L^{(k)}H \big],
\label{eq:Q}
\end{align}
for $\be\in [\be_L, \be_R]$.
\end{proposition}
\begin{remark}\label{rem:2}
(i) Note that (\ref{eq:L}) and (\ref{eq:Q}) are the linear ordinary differential equations for $D_\la \v{w}$ and $D_\la^2 \v{w}$, respectively. Equivalently, by integration, we can formulate (\ref{eq:L}) as
\begin{equation}\label{eq:L-solu}
D_{\la}\v{w}(0,\be)=\mathcal{L}^{(k)} D_{\la}\v{w}(0,\be_L)+\mathcal{SL}^{(k)},
\end{equation}
and (\ref{eq:Q}) as
\begin{equation}\label{eq:Q-solu}
D_{\la}^2\v{w}(0,\be)=\mathcal{Q}^{(k)} D_{\la}^2\v{w}(0,\be_L)+\mathcal{SQ}^{(k)},
\end{equation}
where $\mathcal{L}^{(k)}$ and $\mathcal{Q}^{(k)}$ are both the $(m-1)\times (m-1)$ matrices and $\mathcal{SL}^{(k)}$, $\mathcal{SQ}^{(k)}$ are the $(m-1)$ vectors.

In particular, for most physical systems, the $S$ in (\ref{eq:ms}) is a sub-triangular matrix, so are $B$ and $(\la-B)^{-1}$ in (\ref{eq:L}) and (\ref{eq:Q}).
Hence (\ref{eq:L-solu}) and (\ref{eq:Q-solu}) can be obtained by integrating component by component of $D_\la\v{w}$ in (\ref{eq:L}) and $D_\la^2\v{w}$ in (\ref{eq:Q}), respectively. \\

(ii) Corollary \ref{cor:1} ensures that $\p_x\v{w}$ remains continuous across both the head $\be$-curve: $\be=\be_L$ and the tail $\be$-curve: $\be=\be_R$.
Moreover,
(\ref{eq:L}) is equivalent to the following equation of $\p_x\v{w}$
\begin{equation}\label{eq:L-wx}
\f{\p}{\p \be}\big[\p_x\v{w}(0,\be)\big]=(\la I-B)^{-1}\big[\f{\p}{\p \be} B \p_x\v{w} -\f{\p}{\p \be}(L^{(k)} H) \big],
\end{equation}
which can be formulated as
\begin{equation}\label{eq:L-wx-solu}
\p_x\v{w}(0,\be)=\mathcal{M}^{(k)} \p_x\v{w} (0,\be_L)+ \mathcal{SM}^{(k)}.
\end{equation}
However, we can not derive an equation analogous to (\ref{eq:L-wx}) for $\p_x^2\v{w}(0,\be)$ since $\p_x^2\v{w}(0,\be)$ for $\be\in(\be_L, \be_R)$ does not take a finite value in general.

\end{remark}
{\bf Proof of Proposition \ref{pro:LQ}.} We make use of the regularity of Riemann invariant $\v w$.
Let us first differentiate $\v{w}$ with respect to $\al$ and $\be$ to get
\begin{eqnarray}\label{eq:gw-a1}
\f{\p}{\p \al}\f{\p \v w}{\p \be}=\f{\p}{\p \al}\l(\f{\p t}{\p \be}D_\mu \v{w}\r)
=\f{\p^2 t}{\p \al\p \be}D_\mu\v{w}+\f{\p t}{\p\be}\f{\p}{\p \al}D_\mu\v{w}.
\end{eqnarray}
Similarly, one has
\begin{eqnarray}\label{eq:gw-a2}
\f{\p}{\p \be}\f{\p \v w}{\p \al}=\f{\p}{\p \be}\l(\f{\p t}{\p\al} D_\la\v{w}\r)=\f{\p^2 t}{\p \be\p \al}D_\la\v{w}+\f{\p t}{\p \al}\f{\p}{\p \be}D_\la\v{w}.
\end{eqnarray}
Subtracting these two equations yields
$$
\f{\p t}{\p \al}\f{\p}{\p \be}D_\la\v{w}=
\f{\p t}{\p\be}\f{\p}{\p \al}D_\mu\v{w}+\f{\p^2 t}{\p \al\p \be}(D_\mu\v{w}-D_\la\v{w}).
$$
Using (\ref{eq:al0}) and (\ref{eq:gt1}), one can obtain
\begin{equation}\label{eq:L-pre}
\f{\p}{\p \be} \big[ D_\la\v{w}(0,\be) \big]=\f{1}{\mu-\la}\l(D_\mu-D_\la\r)\v{w}=\p_x\v{w}.
\end{equation}
Recall Proposition \ref{pro:2}. Then we arrive at (\ref{eq:L}).

We proceed by differentiating (\ref{eq:gw-a1}) and (\ref{eq:gw-a2}) with respect to $\al$ to obtain
\begin{eqnarray*}
\f{\p^3}{\p \al^2 \p \be}\v{w}=\f{\p^2}{\p \al^2}\l(\f{\p t}{\p \be}D_\mu \v{w}\r)
=\f{\p^3}{\p \al^2 \p \be}D_\mu\v{w}+2\f{\p^2 t}{\p \al\p\be}\f{\p}{\p \al}D_\mu\v{w}+ \f{\p t}{\p\be}\f{\p^2}{\p \al^2}D_\mu\v{w},
\end{eqnarray*}
and
\begin{align*}
\f{\p^3}{\p \be \p \al^2}\v{w}&=\f{\p^2}{\p \be \p \al}\l(\f{\p t}{\p\al} D_\la\v{w}\r)=\f{\p}{\p \be}\l( \f{\p^2 t}{\p \al^2}D_\la \v{w} +\f{\p t}{\p \al} \f{\p }{\p \al}D_\la \v{w}\r)\nonu\\
&=\f{\p^3 t}{\p \be \p \al^2} D_\la \v{w}+
\f{\p^2 t}{\p \al^2}\f{\p}{\p \be} D_\la\v{w}+2\f{\p^2 t}{\p \al\p \be}\f{\p t}{\p \al}D_\la^2 \v{w}+ \l(\f{\p t}{\p \al}\r)^2\f{\p}{\p \be}D_\la^2\v{w}.
\end{align*}
Subtract the above two equations and then set $\al=0$ to yield (using (\ref{eq:al0}) and (\ref{eq:gt1}) again)
\begin{align*}
\l(\f{\p t}{\p \al}\r)^2\f{\p }{\p \be} \big[ D_\la^2 \v{w}(0,\be) \big]&=\f{2}{\mu-\la}\l(\f{\p t}{\p \al}\r)^2 D_\la\big((\mu-\la)\p_x \v{w} \big)\\
&\quad +\f{\p^3 t}{\p \be \p \al^2}(\mu-\la)\p_x\v{w}-\f{\p^2t}{\p \al^2} \f{\p}{\p \be}D_\la \v{w}.
\end{align*}
Recalling (\ref{eq:L-pre}), it follows that
\begin{align*}
\l(\f{\p t}{\p \al}\r)^2\f{\p }{\p \be} \big[ D_\la^2 \v{w}(0,\be) \big]&=2\l(\f{\p t}{\p \al}\r)^2 D_\la (\p_x\v{w}) +\f{2}{\mu-\la}\l(\f{\p t}{\p \al}\r)^2 D_\la(\mu-\la) \p_x\v{w} \\
&\quad+ \l[\f{\p^3 t}{\p \be \p \al^2}(\mu-\la)- \f{\p^2t}{\p \al^2} \r] \p_x \v{w}.
\end{align*}
Inserting (\ref{eq:gt2}) into the last term of the above equation, we can obtain (after suitable reduction)
\begin{align}\label{eq:Q-pre}
\f{\p }{\p \be}\big[ D_\la^2 \v{w}(0,\be) \big]=2 D_\la (\p_x\v{w}) +\f{\p}{\p \be} (D_\la\la) \p_x\v{w}.
\end{align}
Recalling Proposition \ref{pro:2} once more, we  obtain (\ref{eq:Q}).
\ef

\section{Resolution of curved discontinuities }
\label{sec:R-DW}
\setcounter{equation}{0}
In this section, we resolve the curved discontinuity wave, which can be a contact discontinuity or a shock wave.
Let $\Gamma_k$ be the discontinuity wave and denote by $U_L(x,t)$ (resp. $U_R(x,t)$) the state $U$ on its left (resp. right) side.

\subsection{The contact discontinuity}
Assume for the present that the  $\Gamma_k$ is a curved contact discontinuity. We use the same notations as in the previous section. The propagation speed of $\Gamma_k$ is $\la$ by suppressing the subscript for simplicity in notations.

A significant feature of contact discontinuity is that the generalized Riemann invariant remains continuous across the wave. Thus we take differentiations of $w$ along the trajectory of $\Gamma_k$ to obtain
\beq\label{eq:ct}
D_\la^\ell\big( \v{w}(U_R) \big)= D_\la^\ell \big( \v{w}(U_L) \big),
\eeq
for $\ell=1, 2$.

By recalling (\ref{eq:w}), we have
\begin{equation}\label{eq:d1w}
D_\la\big( \v{w}(U) \big)=
[(\la I- B ) \na_{U} \v{w}] (\p_x U)+L^{(k)} H.
\end{equation}
Thus, while $\ell=1$, (\ref{eq:ct}) is equivalent to
\begin{equation}\label{eq:ct-1}
[(\la I- B) \na_{U} \v{w}]_R (\p_x U)_R-[(\la I- B ) \na_{U} \v{w}]_L (\p_x U)_L= -(L^{(k)} H)_R +(L^{(k)} H)_L.
\end{equation}

\subsection{The shock wave}
Now, let us assume $\Gamma_k$ is a curved shock wave with propagation speed denoted by $\si$.  Then along the shock trajectory, the Rankine-Hugoniot relation reads
\beq\label{eq:RH}
F(U_R)-F(U_L)=\si (U_R-U_L).
\eeq
Denote $D_\si=\f{\p}{\p t}+\si\f{\p}{\p x}$. By taking the directional derivative of (\ref{eq:RH}) along the shock trajectory $\Gamma_k$, one can get
\beq\label{eq:RH-0}
D_\si^\ell \big(F(U_R)-\si U_R \big)  =
D_\si^\ell \big(F(U_L)-\si U_L \big),
\eeq
for $\ell=1,2$.

While $\ell=1$, by noting that
\begin{align}
D_\si \big(F(U)-\si U \big) &=(A-\si I) D_\si U-D_\si \si U\nonumber\\
&=-(A-\si I)^2(\p_x U)-(A-\si I)(L^{(k)}H)-D_\si \si U,\label{eq:d1RH}
\end{align}
 (\ref{eq:RH-0}) is equivalent to
\begin{align}\label{eq:RH-1}
&-(A_{R}-\si I)^2(\p_x U)_{R}+(A_{L}-\si I)^2(\p_x U)_{L} -D_{\si} \si (U_{R}-U_{L})\nonu\\
&\hspace{4cm}=-(A_{R}-\si I)(L^{(k)}H)_{R}+(A_{L}-\si I)(L^{(k)}H)_{L}.
\end{align}

An alternative approach for resolving the shock wave is by using the $m-1$ Rankine-Hugoniot relations in the form
\begin{align}\label{eq:RH-11}
\Psi(U_L, U_R)=0, \quad \Psi=(\Psi^1,\cdots, \Psi^{(m-1)}),
\end{align}
which is equivalent to (\ref{eq:RH}). By differentiating (\ref{eq:RH-11}) along the shock wave, the relation equations for $(D_\si^kU)_R$ and $(D_\si^kU)_L$($k= 1, 2$) can be obtained directly.
This later approach is usually more efficient for practical use, since $\si$ does not appear in (\ref{eq:RH-11}) now. However, for a general purpose, we shall use the former approach in the following discussion.

\section{The GRP solvers}
\label{GRP-alg}
\setcounter{equation}{0}
In this section, we will present the full solver for the linear GRP and quadratic GRP. Indeed,
since the solution $U$ is smooth in the region on the left (resp. right) of $\Gamma_1$ (resp. $\Gamma_m$), the time derivatives of $U$ are thus determined by (\ref{eq:cl}) and the initial data (\ref{eq:clini}).

For the nonsonic case, it suffices for us to determine the spatial derivatives $\p_x U$ and $\p_x^2U$ in the intermediate regions of $\Gamma_k$ ($k=1,\cdots,m$), since the times derivatives $\p_t U$ and $\p_t^2 U$ follows directly from (\ref{eq:cl}). For the sonic case that $t$-axis lies inside the rarefaction fan, we need to give an independent treatment.

\subsection{The nonsonic case}
The nonsonic case refers to the case that  the $t$-axis is located in the intermediate regions of $\Gamma_k$ ($k=1,\cdots,m$).
The $m$ waves $\Gamma_k$, $k=1,\cdots, m$ separates the half space $t>0$ into $m+1$ regions. The region on the left (right) of $\Gamma_k$ is labeled as $\Omega_{k-1/2}$ ($\Omega_{k+1/2}$). The associated state of $U$ in $\Omega_{k-1/2}$ is labeled as $U_{k-1/2}$. The same notation apply for the derivatives of $U$, such as $(\p_x U)_{k-1/2}$.

Now let us summarize the resolution of the linear GRP for the nonsonic case in the following proposition.
\begin{proposition}[Linear GRP: Nonsonic case]\label{pro:nonsonic}
Assume that the solution of problem (\ref{eq:cl}) and (\ref{eq:clini}) consists of $m$ waves $\Ga_k$, $k=1,\cdots, m$.
Then the $(m-1)\times m$ unknowns $(\p_x U)_{k-1/2}$, ($k=2,\cdots,m$) in the intermediate regions of $\Gamma_K$ and the number $D_{\si_k}\si_k$ are determined by the following linear algebraic system
\beq\label{eq:sys-LGRP}
\left\{\begin{array}{l}
(\na_U\v{w})_{k+1/2} (\p_x U)_{k+1/2}-\mathcal{M}^{(k)} (\na_U\v{w})_{k-1/2} (\p_x U)_{k-1/2}=\mathcal{SM}^{(k)}, \\[2mm]
\hspace{2.5cm}\textrm{if}\  \Gamma_k\  \textrm{is a rarefaction wave;}\\[3mm]
[(\la_k I-B^{(k)})\na_U\v{w}]_{k+1/2} (\p_x U)_{k+1/2}-[(\la_k I-B^{(k)})\na_U\v{w}]_{k-1/2} (\p_x U)_{k-1/2}\\
\hspace{3.5cm}
=-(L^{(k)} H)_{k+1/2}+(L^{(k)}H)_{k-1/2}, \\[2mm]
\hspace{2.5cm}\textrm{if}\  \Gamma_k\  \textrm{is a contact discontinuity wave;}\\[3mm]
-(A_{k+1/2}-\si_k I)^2(\p_x U)_{k+1/2}+(A_{k-1/2}-\si_k I)^2(\p_x U)_{k-1/2} -D_{\si_k} \si_k (U_{k+1/2}-U_{k-1/2})\\
\hspace{3.5cm}=-(A_{k+1/2}-\si_k I)(L^{(k)}H)_{k+1/2}+(A_{k-1/2}-\si_k)(L^{(k)}H)_{k-1/2}\\[2mm]
\hspace{2.5cm}\textrm{if}\ \Gamma_k\ \textrm{is a shock wave.}\\
\end{array}
\right.
\eeq
 Here, the $\mathcal{M}^{(k)}$ and $\mathcal{SM}^{(k)}$ are as in (\ref{eq:L-wx-solu}). $(\p_x U)_{1/2}= (\p_x U)_L$, $(\p_x U)_{m+1/2}= (\p_x U)_R$ and the $U_{k+1/2}$ in the coefficients are determined by $R^A(\theta, U_-, U_+)$.
Having solved $(\p_x U)_{k-1/2}$, the time derivatives $(\p_t U)_{k-1/2}$ ($k=1,\cdots,m+1$) are determined by
\begin{equation}\label{eq:Ut-det}
(\p_t U)_{k-1/2}=-A(U_{k-1/2})(\p_x U)_{k-1/2}+H(x,U_{k-1/2}).
\end{equation}
\end{proposition}

\begin{remark}
To solve (\ref{eq:sys-LGRP}), we suggest that the unknowns be ordered as
$$(\cdots,D_{\si_k}\si_k,U_{k-1/2},U_{k+1/2},\cdots)
$$
if $\Gamma_k$ is a shock wave, and use {\it Gauss-Jordan elimination with rows partial pivoting}.
\end{remark}

{\bf Proof of Proposition \ref{pro:nonsonic}} As illustrated previously, the solution $U$ of (\ref{eq:cl}) and (\ref{eq:clini})
is smooth in the regions $\Omega_{k-1/2}$, $k=1,\cdots,m+1$. In the regions $\Omega_{1/2}$ and $\Omega_{m+1/2}$, the spatial derivatives $\p_x U$ are determined by the initial data $(\p_x U)_L$ and $(\p_x U)_R$, respectively. As indicated by the resolution of rarefaction wave and discontinuous waves in Sections \ref{rarefaction} and \ref{sec:R-DW}, the relations between $(\p_x U)_{k+1/2}$ and $(\p_x U)_{k-1/2}$ are described by a set of linear algebraic equations. (\ref{eq:sys-LGRP}) is obtained by combining (\ref{eq:L-wx-solu}), (\ref{eq:ct-1}) and (\ref{eq:RH-1}) and (\ref{eq:Ut-det}) follows directly from (\ref{eq:ncl}).
\ef

To present the quadratic GRP solver, we need to give a few formulations.
In the regions where the flow is smooth, by applying $\p_x$ and $\p_t$ to (\ref{eq:cl}), we have
\begin{align}
\label{eq:u-tx}
&\p_{t}(\p_xU)=-A \p_{x}^2U-\p_{x}A \p_{x}U+\p_{x}H,\\
\label{eq:u-tt}
&\p_{t}^2U=-A \p_{t}(\p_xU)-\p_{t}A \p_{x}U+\p_{t}H.
\end{align}
Inserting (\ref{eq:cl}) and (\ref{eq:u-tx}) into (\ref{eq:u-tt}), $\p_{t}^2 U$ can be expressed as a function of $U, \p_x U, \p_{x}^2U$:
\begin{equation}\label{eq:u-ex}
\p_{t}^2U=\mathcal{A}_Q(U, \p_{x}U, \p_{x}^2U).
\end{equation}

For the GRI $\v{w}$ of $k$th characteristic fields, by noticing that
\begin{align*}
&\p_{x}^2\v{w}=\p_x(\na_U \v{w}) \p_xU+\na_U \v{w} \p_{x}^2 U\\[2mm]
&D_{\la_k} (\p_x\v{w})=(\la_k I-B^{(k)})\p_x^2\v{w}-\p_xB^{(k)}\p_x\v{w}+\p_x (L^{(k)} H),\\
&D_{\la_k}^2\v{w}=D_{\la_k}(\la_k I-B^{(k)})\p_x\v{w}+(\la_k I-B^{(k)})D_{\la_k}(\p_x\v{w})+D_{\la_k} (L^{(k)} H),
\end{align*}
we can get
\begin{equation*}\label{eq:exp-d2w}
 D_{\la_k}^2\v{w}=M_{r}^{(k)}(U) \p_x^2 U+B_r^{(k)}(U,\p_x U),
\end{equation*}
with
\begin{equation}\label{eq:mrbr}
\begin{split}
& M_{r}^{(k)}(U)=(\la_k I-B^{(k)} )^2\na_U \v{w},\\[2mm]
& B_r^{(k)}(U,\p_xU)=[D_{\la_k}(\la_k I -B^{(k)})-(\la_k-B^{(k)}) \p_x B^{(k)}] \p_x\v{w}\\
& \hspace{1cm}+(\la_k I-B^{(k)})^2 \p_x(\na_U \v{w}) \p_x U+(\la_kI-B^{(k)}) \p_x(L^{(k)} H) +D_{\la_k}(L^{(k)} H).
\end{split}
\end{equation}

To resolve the shock wave, we shall use
\begin{align*}
D_\si^2 \big(F(U)-\si U \big)&= D_\si \big( (A-\si I) D_\si U- D_\si \si U \big)\\
&=D_\si \big( -(A-\si I)^2 \p_xU+ (A-\si I) H - D_\si \si U \big)\\
&=M_s(U,\si)\p_x^2 U- D_\si^2 \si U+ B_s(U, \p_x U, \si, D_\si \si),
\end{align*}
with
\begin{equation}\label{eq:msbs}
\begin{split}
& M_s(U,\si)=(A-\si I)^3,\\[2mm]
& B_s(U, \p_x U, \si, D_\si \si)=(D_\si A-2 D_\si \si I ) D_\si U  +(A-\si I)^2 (\p_x A\p_x U-\p_x H)\\
&\hspace{3cm}+(A-\si I)[-D_\si (A-\si I) \p_xU +D_\si H].
\end{split}
\end{equation}

Similar to Proposition \ref{pro:nonsonic}, by combining (\ref{eq:Q}) and (\ref{eq:ct}), (\ref{eq:RH-0}) with $\ell=2$, we have the following proposition for the quadratic GRP solver in nonsonic case.

\begin{proposition}[Quadratic GRP: Nonsonic case]\label{pro:QGRP}
Assume that the solution of problem (\ref{eq:cl}) and (\ref{eq:clini}) consists of $m$ waves $\Ga_k$, $k=1,\cdots, m$.
Then
the $(m-1)\times m$ unknowns $(\p_x^2 U)_{k-1/2}$ ($k=2,\cdots,m$) in the intermediate regions of $\Gamma_k$ and the number $D_{\si_k}^2\si_k$ are determined by the following linear algebraic system
\beq\label{eq:sys-QGRP}
\left\{\begin{array}{l}
M^{(k)}_r (\p_x^2 U)_{k+1/2}- \mathcal{Q}^{(k)} M_r^{(k)}(U_{k-1/2}) (\p_x^2 U)_{k-1/2}=\\
 \hspace{2.4cm}-B_r^{(k)}(U_{k+1/2},(\p_x U)_{k+1/2} )+\mathcal{Q}^{(k)} B_r^{(k)}(U_{k-1/2},(\p_x U)_{k-1/2} ), \\[2mm]
\hspace{4cm}\textrm{if}\  \Gamma_k\  \textrm{is a rarefaction wave;}\\[3mm]
M^{(k)}_r (\p_x^2 U)_{k+1/2}- M_r^{(k)}(U_{k-1/2}) (\p_x^2 U)_{k-1/2}=\\
 \hspace{2.4cm}-B_r^{(k)}(U_{k+1/2},(\p_x U)_{k+1/2} )+ B_r^{(k)}(U_{k-1/2},(\p_x U)_{k-1/2} ), \\[2mm]
\hspace{4cm}\textrm{if}\  \Gamma_k\  \textrm{is a contact discontinuity wave;}\\[3mm]
 M_s(U_{k+1/2},\si)(\p_x^2 U)_{k+1/2} - M_s(U_{k-1/2},\si)(\p_x^2U)_{k-1/2}- D_{\si_k}^2 \si_k (U_{k+1/2}-U_{k-1/2})= \\
 \hspace{1.8cm}-B_s(U_{k+1/2},(\p_x U)_{k+1/2}, \si_k, D_{\si_k} \si_k)+ B_s(U_{k+1/2},(\p_x U)_{k+1/2}, \si_k, D_{\si_k} \si_k), \\[2mm]
\hspace{4cm}\textrm{if}\ \Gamma_k\ \textrm{is a shock wave.}\\
\end{array}
\right.
\eeq
Here, the $ M_{r}^{(k)}(U)$, $B_r^{(k)}(U,\p_xU)$, $M_s(U,\si)$ and $B_s(U, \p_x U, \si, D_\si \si)$ are as in (\ref{eq:mrbr}) and (\ref{eq:msbs}).  $(\p_x^2U)_{1/2}= (\p_x^2U)_L$ and $(\p_x^2U)_{m+1/2}= (\p_x^2U)_R$. The $U_{k+1/2}$ and $(\p_xU)_{k+1/2}$ in the coefficients are determined by Proposition \ref{pro:nonsonic}.
Having solved $(\p_x^2 U)_{k-1/2}$, the time derivatives $(\p_t^2 U)_{k-1/2}$ ($k=1,\cdots,m+1$) can be obtained by using
\begin{equation*}
(\p_t^2 U)_{k-1/2}=\mathcal{A}_Q(U_{k-1/2}, (\p_{x}U)_{k-1/2}, (\p_{x}^2U)_{k-1/2}).
\end{equation*}

\end{proposition}

\subsection{The sonic case}
\label{sec:sonic}
As far as the sonic case is concerned, the $t$-axis is located inside the rarefaction wave, $\Gamma_k$ for instance, and is in fact tangential to the $\la_k$-characteristic curve. Thus, for this case, we need to solve $D_{\la_k} U$ and $D_{\la_k}^2 U$.
Moreover, the explicit expression of $D_{\la_k} U$ (with respect to $\be$) is required when solving the Q-equations (\ref{eq:Q}), as is the main step of the QGRP solver.

Although  $D_{\la_k} \v{w}$ and $D_{\la_k}^2 \v{w}$ are readily obtained from (\ref{eq:L-solu}) and (\ref{eq:Q-solu}), we still need to make up an additional freedom. In fact, (\ref{eq:lncl}) implies
the differential relations of $U$ along the $\la_k$ characteristic curve
\begin{equation}\label{eq:lncl-k}
L_k  D_{\la_k} U =L_k H.
\end{equation}
 Combining (\ref{eq:lncl-k}) with
\beq\label{eq:w-u}
\na_U \v{w}D_{\la_k} U=D_{\la_k}\v{w},
\eeq
$D_{\la_k} U$ can be determined.

Furthermore, applying $D_{\la_k}$ to (\ref{eq:lncl-k}) and (\ref{eq:w-u}) yields
\begin{equation}\label{eq:lncl-k-2}
L_k D_{\la_k}^2 U = -D_{\la_k}L_k  D_{\la_k}U+ D_{\la_k}(L_k H),
\end{equation}
and
\beq\label{eq:w-u-2}
\na_U \v{w}D_{\la_k}^2U=-D_{\la_k} (\na_U \v{w}) D_{\la_k}U+ D_{\la_k}^2\v{w}.
\eeq
Then $D_{\la_k}^2 U$ can be solved by combining (\ref{eq:lncl-k-2}) and (\ref{eq:w-u-2}).

In the sonic case, where $\la_k=0$, we use the following observation for the LGRP,
\begin{align}\label{eq:sonic-ut}
U_t(0,0^+)=D_{\la_k} U(0,0^+).
\end{align}
Also, by taking $\la_k=0$, we have
\begin{align}\label{eq:sonic-utt}
\p_t^2U(0,t)=D_{\la_k}^2 U(0,t)-D_{\la_k} \la_k \p_xU(0,t),
\end{align}
for $t>0$. However, the above observation can not be used to calculate $\p_t^2U(0,0^+)$, since generally neither $\p_t^2U(0,0^+)$ nor $\p_xU(0,0^+)$ takes finite value inside the rarefaction wave fan, expect for $U$ being the GRI $\v{w}$. See Remark \ref{rem:2}.

\begin{figure}
  \begin{center}
    \includegraphics[height=2.0in, width=3.2in, trim=0 0 0 0, clip]{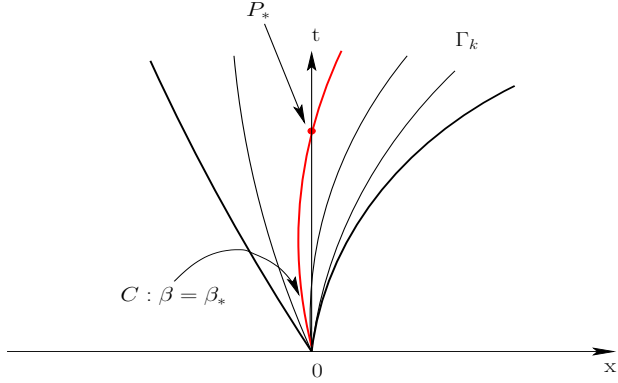}
  \end{center}
\caption{The characteristic curve $\be=\be_*$, red line.}
  \label{fig:sonic}
\end{figure}

For the QGRP, we shall use the following method to give a second order in time approximation of $U$ in $t$-axis.
For any point $P_*=(0, \Delta t)$ with $\Delta t$ being small, to evaluate $U(P_*)$, we need to find the initial slope $\be_0$ of the characteristic curve $C$ which emanates from the singularity and goes though $P_*$. See Fig. \ref{fig:sonic}.  Since for any $(x(t),t)\in C$,
\begin{equation}\label{eq:curve-app}
x(t)=\int_{0}^{t} \la \diff s=\int_{0}^{t} \big(\la_k(0) +D_{\la_k}\la_k(0) s +O(s^2) \big)\diff s,
\end{equation}
the initial slope $\be_*=\la_k(0)$ can be approximated by solving
\begin{equation}\label{eq:be-app}
\be_*+D_{\la_k} \la_k (\be_*)\f{\Delta t}{2}=0,
\end{equation}
for which, we can use the Newton iteration with initial guess $\be_*=0$.

Having determined $\be_*$, $U(P_*)$ can be evaluated as
\begin{equation}\label{eq:U-sonic-app}
U(P_*)\approx U(\be_*)+D_{\la_k} U (\be_*) \Delta t+D_{\la_k}^2 U(\be_*) \f{\Delta t^2}{2}.
\end{equation}

\section{The acoustic approximation}

\setcounter{equation}{0}
As $U_-=U_+$ and $\p_x^\ell P_-(0^-)\neq \p_x^\ell P_+(0^+)$, we refer it to as the acoustic case and all waves $\Gamma_k$ are acoustic.  Fixed $\la_k$, the wave $\Gamma_k$ degenerates to a characteristic curve and  the states $U_L$, $U_R$ on both sides of $\Gamma_k$ are the same. In particular, as the initial data has a small jump  $\|U_-- U_+\|\ll 1$, we adopt the {\em acoustic approximation} in the sense that $U_-$ and $U_+$ are regarded as the same approximately.

Now let us look at the acoustic wave $\Gamma_k$.
We use the continuity property of $U$  and make differentiation along $\Gamma_k$ to obtain  $D_{\la_k}U_L=D_{\la_k}U_R$. Then we proceed to use \eqref{eq:ncl} to get
\begin{equation}
(\la_k I-A)_L (\p_xU)_L +H_L =(\la_k I-A)_R (\p_xU)_R+ H_R.
\label{eq:ac-0}
\end{equation}
Note that $U_L=U_R$ and recall the notation $L^{(k)} =(L_1,\cdots, L_{k-1},L_{k+1},\cdots, L_m)^T$ in \eqref{eq:ms}. Then we find that the \eqref{eq:ac-0} is equivalent to
\begin{equation}\label{eq:acous-1}
L^{(k)} (\p_xU)_L=L^{(k)} (\p_xU)_R.
\end{equation}
Moreover, applying $D_{\la_k}$ to (\ref{eq:acous-1}) yields
\begin{align}\label{eq:acous-2}
(\la_k I- \La^{(k)}) L^{(k)} & \big( (\p_x^2U)_R- (\p_x^2U)_L \big)=\nonu\\
& D_{\la_k} L^{(k)}\big((\p_xU)_L-(\p_xU)_R\big) -L^{(k)}\big( (\p_xA \p_xU)_L-(\p_xA \p_xU)_R \big).
\end{align}
In addition, if $(\p_xU)_L= (\p_xU)_R$, then (\ref{eq:acous-2}) is reduced to
\begin{equation}\label{eq:acous-3}
L^{(k)} (\p_x^2U)_L=L^{(k)} (\p_x^2U)_R.
\end{equation}

Interestingly, in the course of acoustic approximation,  we can obtain $\p^\ell_x U$ equivalently by solving  linear  classical Riemann-type  problems
\begin{equation}\label{eq:acous-RP}
\left\{\begin{array}{l}
\f{\p }{\p t} \l(\p_x^\ell U \r) +A(U_*) \f{\p }{\p x} \l(\p_x^\ell U\r)=0,\\[2mm]
 \p_x^\ell U (x,0)=
\left\{\begin{array}{ll}
\p_x^\ell P_-(0) & \textrm{if}\ x<0,\\
\p_x^\ell P_+(0)  & \textrm{if}\ x>0.
\end{array}
\right.
\end{array}
\right.
\end{equation}
with  $U_*=(U_-+U_+)/2$. Note that the components of $L^{(k)} \p_x U$ in (\ref{eq:acous-1}) are nothing but the $m-1$ generalized Riemann invariants of system (\ref{eq:acous-RP}) associated with $\la_k$.
We also note that, as indicated by (\ref{eq:acous-1}) or (\ref{eq:acous-RP}), the spatial derivatives $\p_x U$ are independent of the source term $H$.
In addition, if $ \p_x P_- (0)  \approx  \p_x P_+(0)$, from (\ref{eq:acous-3}), we see that $\p_x^2U$ can also be approximated by solving (\ref{eq:acous-RP}) with $\ell=2$.
In general, we have the following proposition.

\begin{proposition}\label{pro:acoustic}
  For any $k\geq 1$, assume that we have $\p_x^\ell P_- (0)  =  \p_x^\ell P_+(0)$ ($0\leq \ell \leq k-1$ with $\p_x^0 U$ stands for $U$). Then $\p_x^k U$ are determined by the linear system (\ref{eq:acous-RP}) with $\ell=k$.

\end{proposition}
\begin{remark}
(i) Since $U$ is analytical in the regions $\Omega_{j-1/2}$ ($j=1,\cdots,m+1$), the corresponding time derivatives $\p_t^\ell U$, $0 \leq \ell \leq k$ follow from the Cauchy-Kowalewski procedure as illustrated in \cite{toro}. \\

(ii) We note here that,
under the acoustic assumption in Proposition \ref{pro:acoustic}, all the approximate DRP solvers proposed in \cite{eno,tt,ct1} are valid and are actually equivalent to the present acoustic GRP solvers.
\end{remark}

As for the resolution of GRP (\ref{eq:cl})-(\ref{eq:clini}), if the initial data (\ref{eq:clini}) has a jump discontinuity, we can derive the solvers analytically as in Section \ref{GRP-alg} to calculate the time derivatives of $U$, with possible acoustic approximation for a partial set of waves. This leads to the solver which we label as the {\it LGRP$_\infty$ (QGRP$_\infty$) solver}.
While the jump $U_+-U_-$ of $U$ is very small, we can use (\ref{eq:acous-RP}) (or possibly (\ref{eq:acous-2})) to calculate the space derivatives approximately. The resulting LGRP (QGRP) solver is labeled as the {\it LGRP$_1$ (QGRP$_1$) solver}.

\section{An example: The variable area duct flow system}
\label{sec:exm}
In this section, we will take the system of variable area duct flow as an example to test the GRP solvers in the previous section. The flow system is
\begin{align}\label{eq:euler}
&\hspace{1.5cm} \f{\p }{\p t}U+A(x)^{-1}\f{\p}{\p x}\big[ A(x) F(U) \big]+\f{\p}{\p x}G(U)=0,\\
\nonumber
&U=\begin{pmatrix}\rho \\
                           \rho u\\
                           E
                  \end{pmatrix},
\quad
F(U)=\begin{pmatrix}\rho \\
                           \rho u^2+p\\
                           (E+p) u
                  \end{pmatrix},
\quad
G(U)=\begin{pmatrix} 0 \\
                     p \\
                     0
                  \end{pmatrix}.
\end{align}
Here, $\rho$, $u$, $e$ are the density, velocity and internal energy, respectively. $p=p(\rho, e)$ is the pressure, $E=\rho( e+{1}/{2} u^2)$ is the total energy and the function $A(x)$ is the area of the duct. When $A(x)=1$, the system (\ref{eq:euler}) represents the planar compressible Euler equations. We discuss the case of polytropic gases, for which $p=(\gamma-1) \rho e$, where $\gamma$ is the ratio of specific heats.

\subsection{Formulation of the GRP Solvers}
\label{sec:solver-form}
In terms of the primitive variables $Q=(\rho, u, p)$, system (\ref{eq:euler}) can be written, for smooth flow, as
\begin{equation}\label{eq:euler-pv}
\f{\p Q }{\p t}+J \f{\p Q }{\p x}=H,
\quad
J=\left(\begin{array}{ccc} u & \rho & 0\\
                           0 & u & \f{1}{\rho}\\
                            0 & \rho c^2 & u\\
            \end{array}\right),
\quad
H=\begin{pmatrix} -\f{A'(x)}{A(x)}\rho u \\
                     0 \\
                     -\f{A'(x)}{A(x)} \rho c^2 u
                  \end{pmatrix}.
\end{equation}
Here, $c$ is the local speed of sound, given by $c^2=\f{\gamma p}{\rho}$.

The system (\ref{eq:euler}), or equivalently (\ref{eq:euler-pv}), possesses three eigenvalues
$$
\la_-=u-c, \quad \la_0=u, \quad \la_+=u+c.
$$
The three pairs
$$
\v{w}_-=(S,\psi),\quad \v{w}_0=(u,p), \quad \v{w}_+=(S,\phi)
$$
are
the generalized Riemann invariants associated with $\la_-$, $\la_0$, $\la_+$. Here, $S=p \rho^{-\gamma}$ is the the entropy, and the two varibles $\psi$, $\phi$ are
$$
\psi=u+ \f{2}{\gamma-1} c \quad \phi=u- \f{2}{\gamma-1} c.
$$

\begin{figure}
  \begin{center}
    \includegraphics[height=1.8in, width=2.8in, trim=0 0 0 0, clip]{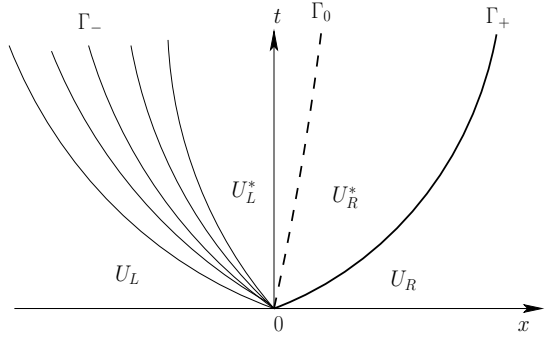}
  \end{center}
\caption{ Typical wave configuration of the variable area duct flow system.}
\label{fig:euler}
\end{figure}

We now start to resolve the generalized Riemann problem for (\ref{eq:euler}) subject to initial data (\ref{eq:clini}).
Assume that the configuration is as shown in Fig.\ref{fig:euler}: a rarefaction wave $\Ga_-$ associated with $\la_-$ moves to the left, a shock wave $\Ga_+$ associated with $\la_+$ moves to the right.
For the variable $U$, let us denote by $U_{L}$ and $U_{R}$ its values on the left-hand side and right-hand side of the three waves, respectively. Similarly, the values of $U$ on the two side of $\Ga_0$ are denoted by $U_{L}^*$ and $U_{R}^*$, as is illustrated by Fig. \ref{fig:euler}. Similar notations will be applied to other variables. For example, $(\p_x U)_{L}^*$ is the value of $\p_x U$ on the left-hand side domain of the contact discontinuity.

As stated in Section \ref{GRP-alg}, to resolve rarefaction wave associated with $\la_-$, we need to use the associated GRI $\v{w}_-$. Indeed, in view of Proposition \ref{pro:2}, we have the following equation for $\v{w}_-$,
\begin{equation}\label{eq:euler-ri}
\f{\p \v{w}_- }{\p t} + B_-
\f{\p \v{w}_-}{\p x}
=H_-,
\quad
B_-=\l( \begin{array}{cc} u & 0\\
       -\f{1}{\gamma-1} \rho^{\gamma-1} & u+c
            \end{array}
\r),
\quad
H_-=\begin{pmatrix}   0 \\
                     -\f{A'(x)}{A(x)} c u
                  \end{pmatrix}.
\end{equation}

By recalling $S(0,\be)=S_L$, $\psi(0,\be)=\phi_L$, $\la_-(0,\be)=\be$ in the rarefaction wave fan and using the L(Q)-equations in Proposition \ref{pro:LQ}, we can obtain the following proposition. The coefficients $A_1,\cdots, A_{19}$, $B_1,\cdots, B_{12}$,  $Z_1(\be)$ and $Z_2(\be)$ are given in \ref{app:1}.
\begin{proposition}[L(Q)-equations: $\gamma\neq 3, 5/3$]\label{prop:euler-lq}
Let $\Ga_-$ be the rarefaction wave as in Fig. \ref{fig:euler}. Then for $\gamma\neq 5/3, 3$, we have
\begin{equation}\label{eq:euler-L}
\begin{split}
& D_{\la_-} S(0,\be) = A_1 (\psi_L -\be )^{\f{\gamma+1}{\gamma-1}},\\
& D_{\la_-} \psi(0,\be)=A_2 (\psi_L -\be )^{\f{\gamma+1}{2(\gamma-1)}}+A_3 (\psi_L -\be )^{\f{2\gamma}{\gamma-1}}+B_1(\psi_L-\be)+B_2(\psi_L-\be)^2.
\end{split}
\end{equation}
and
\begin{equation}\label{eq:euler-Q}
\begin{split}
& D_{\la_-}^2 S(0,\be)=\big[(\psi_L -\be_L)^{-\f{2(\gamma+1)}{\gamma-1}} D_{\la_-} S(0,\be_L)+Z_1(\be)-Z_1(\be_L) \big](\psi_L -\be)^{\f{2(\gamma+1)}{\gamma-1}},\\
& D_{\la_-} \psi(0,\be)=\big[(\psi_L -\be_L)^{-\f{\gamma+1}{\gamma-1}} D_{\la_-} \psi(0,\be_L)+Z_2(\be)-Z_2(\be_L) \big] (\psi_L -\be)^{\f{\gamma+1}{\gamma-1}},
\end{split}
\end{equation}
 for $\be\in [\be_L, \be_{L}^*]$.
\end{proposition}
\begin{remark}
(i) For the cases of $\gamma=5/3, 3$, the L-equations of $\v{w}_-$ are given in \ref{app:1}. The corresponding Q-equations are omitted here.
Besides, the formulae for  $\p_x S(0,\be)$ and $\p_x \Psi(0,\be)$, which can be obtained by integrating (\ref{eq:L-wx}) directly, are also omitted here.

(ii) For $\ga\neq 5/3, 3$,
the differential relation in the rarefaction wave fan associated with $\la_-$
\begin{align}\label{eq:diff-1}
D_{\la_-} \phi=-\f{1}{\gamma(\gamma-1)} \f{c}{ S} D_{\la_-} S+\f{A'(x)}{A(x)} cu
\end{align}
leads to
\begin{align}\label{eq:diff-11}
D_{\la_-} \phi= A_4 (\psi_L -\be )^{\f{2\gamma}{\gamma-1}} +B_3(\psi_L-\be)+B_4(\psi_L-\be)^2.
\end{align}
By noticing $\la_-=\f{3-\gamma}{4} \psi+\f{1+\gamma}{4} \phi$,

we have
\begin{align}\label{eq:lala}
\f{\p}{\p \be}( D_{\la_-} \la_-)
= A_7 (\psi_L -\be )^{\f{\gamma+1}{2(\gamma-1)}}+A_8 (\psi_L -\be )^{\f{2\gamma}{\gamma-1}}+B_7(\psi_L-\be)+B_8(\psi_L-\be)^2.
\end{align}

Moreover, $D_{\la_-}^2 \phi$ can be determined by
\begin{align}\label{eq:diff-2}
D_{\la_-}^2 \phi=-\f{1}{\gamma (\gamma-1)} \l(\f{c}{ S} D_{\la_-}^2 S+D _{\la_-}\l(\f{c}{ S}\r) D_{\la_-} S \r)+D_{\la_-}\l(\f{A'(x)}{A(x)} c u\r),
\end{align}
which follows from (\ref{eq:diff-1}).
The above formulas will be used to resolve the sonic case. See Proposition \ref{prop:euler-sonic}.
\end{remark}

We also note that, in order to resolve the contact discontinuity wave $\Ga_0$, the following equation of $\v{w}_0$ will be used,
\begin{equation}\label{eq:euler-w0}
\f{\p \v{w}_0 }{\p t} + B_0
\f{\p \v{w}_0}{\p x}
=H_0,
\quad
B_0=\l( \begin{array}{cc} u & \f{1}{\rho}\\
       \rho c^2 & u
            \end{array}
\r),
\quad
H_0=\begin{pmatrix}   0 \\
                     -\f{A'(x)}{A(x)} \rho c^2 u
                  \end{pmatrix}.
\end{equation}

Now let us present the LGRP$_\infty$ solver for problem (\ref{eq:euler})-(\ref{eq:clini}) in the following proposition, corresponding to the wave configuration in Fig. \ref{fig:euler}.

\begin{proposition}[Linear GRP]\label{prop:QGRP}
Assume a typical wave configuration for the
generalized Riemann problem of (\ref{eq:euler}) and (\ref{eq:clini}) as shown in Fig.\ref{fig:euler}. Then $(\p_x Q)_L^*$ and $(\p_x Q)_R^*$ are determined by the set of linear equations
\begin{equation}\label{eq:lgrp}
\begin{split}
& [(\la_- I-B_-)\na_Q \v{w}_-]_L^* (\p_x Q)_L^*= D_{\la_-}\v{w}_- (U_L^*)-H_-(U_L^*),\\[2mm]
&(\la_0 I-B_0)_L^*(\p_x \v{w}_0 )_L^*- (\la_0 I-B_0)_R^*(\p_x \v{w}_0)_R^*=-H_0(U_L^*)+H_0(U_R^*),\\[2mm]
& [(\na_Q F-\si \na_Q U)(\si I- J)]_R^* (\p_x Q)_R^*- D_\si \si (U_R^*-U_R)
=\\[2mm]
& [(\na_Q F-\si \na_Q U)(\si I- J)]_L^* (\p_x Q)_R -[(\na_Q F-\si\na_Q U) H]_R^*+[(\na_Q F-\si\na_Q U) H]_R.
\end{split}
\end{equation}
Here, $D_{\la_-} \v{w}_- (U_L^*)$ is determined by (\ref{eq:euler-L}), $\si$ is the speed of the shock associated with $\la_+=u+c$.

Moreover, for the sonic case where the $t$-axis is located in the rarefaction, the $\p_t Q$ at $t$-axis ($\be=0$) are determined by
\begin{equation}\label{eq:sonic-exm}
\begin{split}
& [(\la_- I-B_-)\na_Q \v{w}_-] (\p_t Q)= D_{\la_-}\v{w}_--H_-,\\[2mm]
& \p_t u+ \f{1}{\rho c} \p_t p=\f{A'(x)}{A(x)} c u,
\end{split}
\end{equation}
where the $U$ in the coefficients takes value at $t$-axis where $\la_-=0$.
\end{proposition}

{\bf Proof.} The linear system (\ref{eq:lgrp}) for $(\p_x^* U)_L^*$ and $(\p_x^* U)_R^*$ can be obtained by combining (\ref{eq:euler-L}), (\ref{eq:ct}) and (\ref{eq:RH-0}) with $l=1$ and using the following expressions
\begin{align*}
& D_{\la_-} \l(\v{w}_-(U) \r)= (\la_- I-B_-)\na_Q \v{w}_- \p_x Q+ H_-,\\[3mm]
& D_{\la_0} \l( \v{w}_0(U) \r)=(\la_0 I-B_0) \p_x \v{w}_0+ H_0,\\[3mm]
& D_\si [(\na_Q F-\si \na_Q U)]= (\na_Q F-\si \na_Q U)[(\si I- J)(\p_x Q) + H]  - D_\si^2 \si U.
\end{align*}

To resolve the sonic case, we shall use the differential relation along the $\la_-$ characteristic curve:
\begin{equation}\label{eq:diff-form2}
D_{\la_-} u+ \f{1}{\rho c} D_{\la_-} p=\f{A'(x)}{A(x)} c u,
\end{equation}
which is equivalent to (\ref{eq:diff-1}). (\ref{eq:sonic-exm}) follows from (\ref{eq:euler-L}) and (\ref{eq:diff-form2}) by setting $\la_-=0$.
\ef.

The QGRP$_\infty$ solver for the nonsonic case and sonic case are presented in the following two propositions, respectively.

\begin{proposition}[Quadratic GRP: Nonsonic case]\label{prop:QGRP}
Assume a typical wave configuration for the
generalized Riemann problem of (\ref{eq:euler}) and (\ref{eq:clini}) as shown in Fig. \ref{fig:euler}. Then $(\p_x^2 Q)_L^*$ and $(\p_x^2 Q)_R^*$ are determined by the set of linear equations
\begin{equation}\label{eq:qgrp}
\begin{split}
& M_r(U_L^*) (\p_x^2 Q)_L^*= D_{\la_-}^2\v{w}_- (U_L^*)-B_r(U_L^*, (\p_x U)_L^* )\\[2mm]
& M_c(U_L^*)(\p_x^2 \v{w}_0 )_L^*- M_c(U_R^*)(\p_x^2 \v{w}_0)_R^*
=-B_c\big(U_L^*,(\p_x U)_L^*\big)+B_c\big(U_R^*,(\p_x U)_R^*\big),\\[2mm]
& M_s(U_R^*,\si)(\p_x^2 Q)_R^*- D_\si^2 \si (U_R^*-U_R)
=  -B_s\big(U_R^*,(\p_x U)_R^*, \si, D_\si \si \big)  \\
& \hspace{4.5cm}+ M_s(U_R, \si)(\p_x^2 Q)_R+B_s\big(U_R,(\p_x U)_R,\si, D_\si \si \big).
\end{split}
\end{equation}
Here, $D_{\la_-}^2 \v{w}_- (U_L^*)$ is determined by (\ref{eq:euler-Q}), $\si$ is the speed of the shock associated with $\la_+=u+c$,
\begin{align*}
& M_r(U)=(\la_-I-B_-)^2 \na_Q \v{w}_{-},\\[2mm]
& B_r(U, \p_x U)=[D_{\la_-}(\la_- I-B_-)-(\la_- I-B_-)\p_x B_-]\p_x\v{w}_- \\
& \hspace{3cm} +(\la_- I-B_-)^2\p_x(\na_Q \v{w}_-)  \p_x Q +(\la_- I- B_-)\p_x H_- +D_{\ka_-} H_-, \\[2mm]
& M_c(U)=(\la_0 I-B_0)^2,\\[2mm]
& B_c(U,\p_x U)= D_u(\la_0 I-B_0) \p_x\v{w}_0-(\la_0 I-B_0) \p_x B_0 \p_x \v{w}_0 +(\la_0 I- B_0)\p_x H_0 +D_{\la_0} H_0 ,\\[2mm]
& M_s(U,\si)=(\na_QF-\si \na_QU) (\si I- J)^2,\\[2mm]
& B_s(U,\p_x U, \si, D_\si \si)=(D_\si (\na_QF)-2 D_\si \si \na_QU -\si D_\si (\na_QU) ) D_\si Q\\
&\hspace{4.5cm} +(\na_QF- \si \na_QU)[D_\si (\si I-J) -(\si I-J) \p_x J ] \p_x Q\\
&\hspace{4.5cm} +(\na_QF- \si \na_QU)[(\si I-J) \p_x H +D_{\si} H ].
\end{align*}
\end{proposition}

{\bf Proof.} This proposition can be proved by combining (\ref{eq:euler-Q}), (\ref{eq:ct}) and (\ref{eq:RH-0}) with $l=2$ and using the following expressions
\begin{align*}
& D_{\la_-}^2 \l(\v{w}_-(U) \r)= M_r(U) \p_x^2 Q+ B_r(U,\p_x U),\\[3mm]
& D_{\la_0}^2 \l( \v{w}_0(U) \r)= M_c(U) \p_x^2 \v{w}_0+ B_c(U,\p_x U),
\end{align*}
and
\begin{align*}
D_\si^2 \big(F(U)-\si U \big)=
M_s(U,\si)\p_x^2 Q- D_\si^2 \si U+ B_s(U,\p_x U, \si, D_\si \si).
\end{align*}
\ef.

\begin{proposition}[Quadratic GRP: Sonic case]\label{prop:euler-sonic}
Assume that the $t$-axis is located inside the rarefaction wave associated with $\la_-$. Denoting by $\Phi=(S, \phi, \psi)$, then for any point $P_*=(0, \Delta t)$ with $\Delta t$ being small, we have
\begin{equation}\label{eq:euler-sonic-app}
\Phi(P_*)= \Phi(0)+D_{\la_-} \Phi(\be_*) \Delta t+D_{\la_-}^2 \Phi(\be_*) \f{\Delta t^2}{2}+O(\Delta t^3),
\end{equation}
where $\be_*$ is the root of
\begin{equation}\label{eq:euler-sonic-root}
\be + D_{\la_-} \la_-(\be) \f{\Delta t}{2}=0,
\end{equation}
and
$D_{\la_-}^\ell \Phi$, $\ell=1,2$ are determined by (\ref{eq:euler-L})-(\ref{eq:diff-2}).
\end{proposition}

The other wave configurations can be treated similarly. In particular, if a $\la_+$-rarefaction wave is involved, in order to get the linear equation for $ (\p_x^2 Q)_R^*$ analogous to the first equation of (\ref{eq:qgrp}), it requires one to derive the L(Q)-equation for $\v{w}_+$. However,
a better choice for us is to use the following property of system (\ref{eq:euler}): (\ref{eq:euler}) holds true under the transformation $\mathscr{T}:(\rho, u, p, A)(x,t)\rightarrow (\rho,-u,p, A)(-x,t)$. In fact, if we denote by $\tilde{Q}=\mathscr{T} (Q)$,
 then $\tilde{U}_L^*=\mathscr{T}(U_R^*) $, $\tilde{U}_L^*=\mathscr{T}(U_R^*)$, $\la_-(\tilde{U}_R^*)=\la_+(U_L^*)$.
By expressing $D_{\la_-}^2 \v{w}_-(\be)$ as a function of $U_L,(\p_xU)_L,(\p_x^2U)_L, A'(x), A''(x)$ and $\be$: $D_{\la_-}^2 \v{w}_-(\be)= \mathscr{W} \big(U_L,(\p_xU)_L,(\p_x^2U)_L, A'(x), A''(x), \be \big)$, we have
\begin{align}\label{eq:la3}
 M_r(\tilde{U}_L^*, -A'(x), A''(x)) (\p_x^2 \tilde{Q})_L^*=&\mathscr{W} \big( \tilde{U}_L,\p_x \tilde{U}_L,(\p_x^2 \tilde{U})_L, -A'(x), A''(x), \la_-(\tilde{U}_L^*) \big)  \nonumber\\
&\quad -B_r\big(\tilde{U}_L^*, (\p_x \tilde{U})_L^* , -A'(x), A''(x)\big).
\end{align}
Noting that $(\p_x^2 {Q})_R^* =((\p_x^2 \tilde \rho)_L^*,-(\p_x^2 \tilde u)_L^*, (\p_x^2 \tilde p)_L^*)$, (\ref{eq:la3}) is indeed the derived linear equation for $ (\p_x^2 Q)_R^*$. The sonic case corresponding to the $\la_+$-rarefaction wave can be resolved using the same technique.

\subsection{Tests for the GRP solvers}

In this section we assess the performance of the GRP solvers for the compressible Euler equations system, i.e. (\ref{eq:euler}) with $A(x)=1$. The aim is to show, via several test problems, the accuracy and behavior of the present solvers. As tests, we use the generalized Riemann problems proposed by \cite{ct1} and construct new ones with large jumps in pressure.
The first test has no jump discontinuities in the state variables but admits discontinuities in derivatives at the interface.
The more demanding test problems are constructed from the first one, by adding a discontinuity in pressure.
Six new cases are thus generated by varying the strength of the initial pressure jump $\Delta p=(p_L-p_R)/p_R$ at the interface, namely $\Delta p=10^k$, $k=-2,\cdots,3$.
The last test problem for the sonic case are constructed by adding $\Delta u=28$ to the initial flow velocity of the test case corresponding to  $\Delta p=100$.

In \cite{ct1}, the authors test the first five problems using three type of DRP solvers with only partial success.
Since no exact solutions are known, the reference solutions are obtained numerically, by solving the test problems on very fine mesh on the interval $[-1,1]\times [0,t_0]$ of $(x,t)$.
To do this, the authors of \cite{ct1,ct2} suggest using the Random Choice Method or Weighted Average Flux method to avoid the large nonphysical oscillations of early time solution.
To this aspect, a detailed description can be found in \cite{ct1,ct2}, which is beyond the scope of this work. Here, 
our numerical reference solutions are obtained simply by using the Godunov flux in the context of finite volume method and then correcting their values on the early time interval $[0, t_0/20]$ using an interpolation method. Such a measure does not affect  our accuracy tests.

For each of these tests, we will compare the GRP solver based solution at the interface $x=0$, as a function of time determined by (\ref{eq:tal}), with the {\it reference numerical solution}. As will be shown, the present GRP solvers is truly accurate, having the expected accuracy not only for all the test cases in \cite{ct1}, but also for cases with much larger initial jump in state variables.

\subsubsection{Acoustic case with continuous state and jump in derivatives}

\begin{figure}
  \begin{center}
    \includegraphics[height=1.7in, width=2.0in, trim=0 0 0 0, clip]{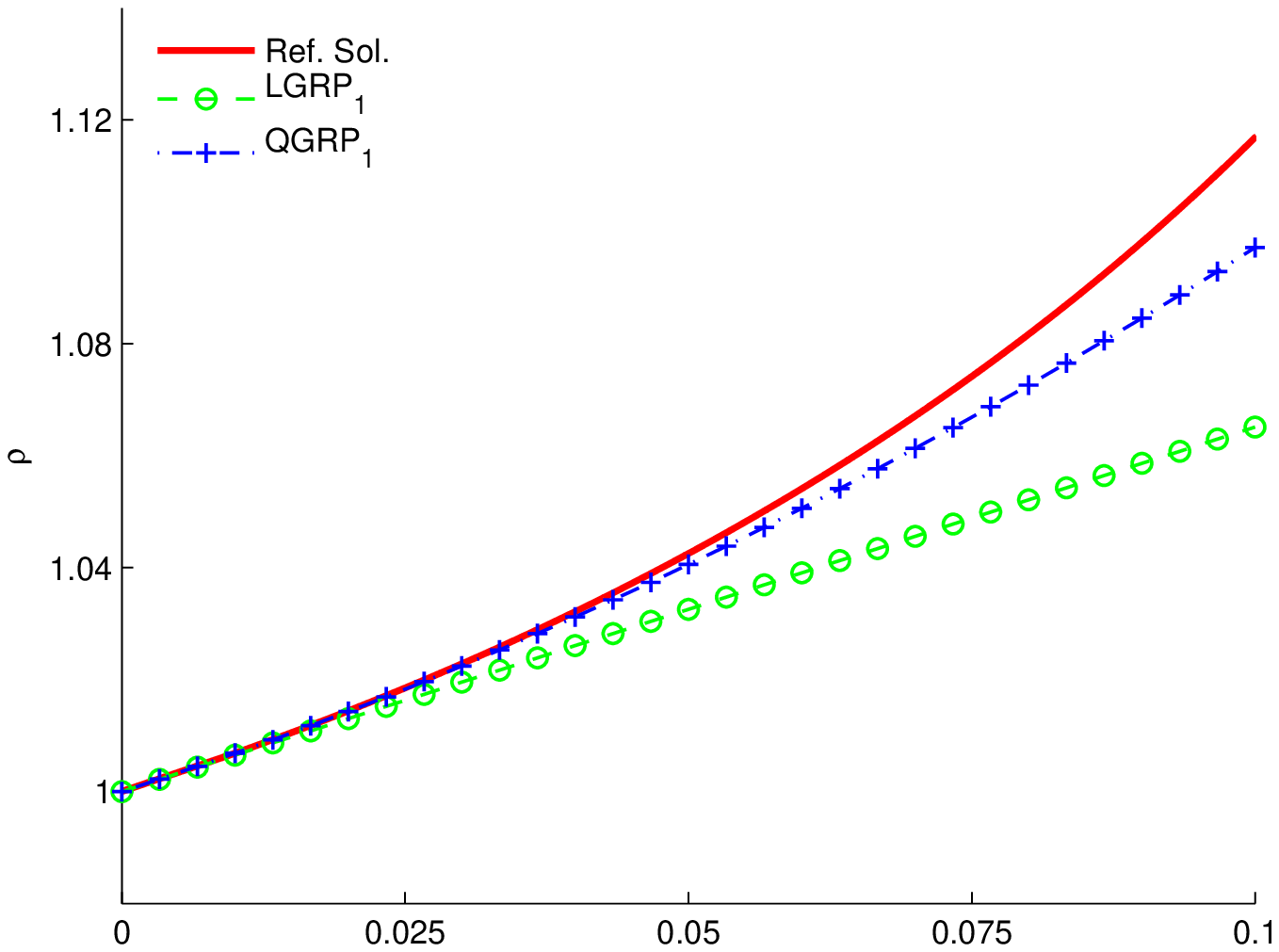}
    \includegraphics[height=1.7in, width=2.0in, trim=0 0 0 0, clip]{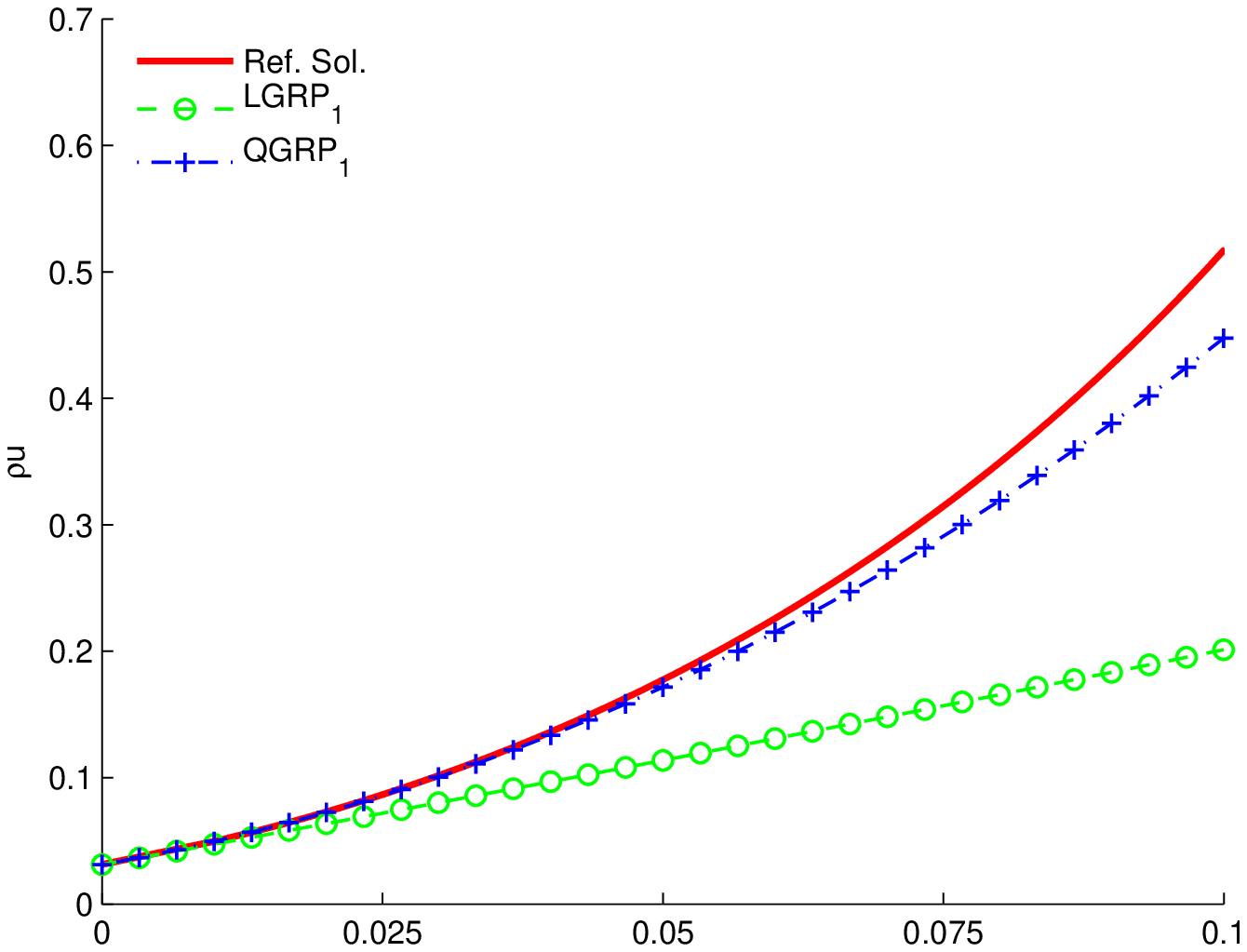}
    \includegraphics[height=1.7in, width=2.0in, trim=0 0 0 0, clip]{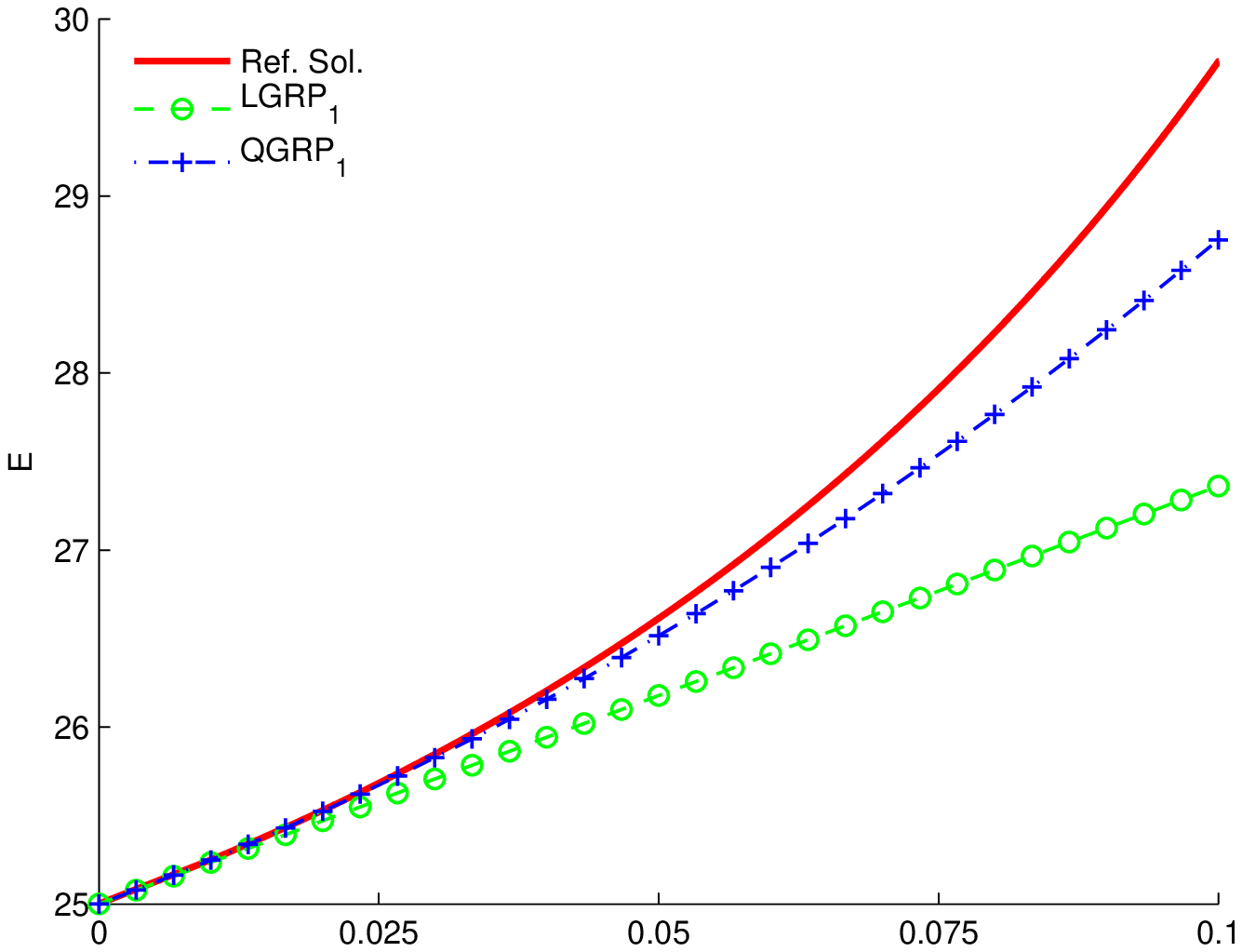}
  \caption{Acoustic case: Reference solution and GRP solvers based solutions.}
  \label{fig:dp0}
  \end{center}
\end{figure}

This test corresponds to the following initial condition
\begin{equation}\label{eq:test1}
\begin{split}
& \rho_L(x,0)=1+0.56431 x+ 2.62896 x^2,  \\[2mm]
& u_L(x, 0)=0.03125-1.024 x+1.92 x^2,\\[2mm]
&p_L(x,0)=10 -0.216 x+1.08 x^2,\\[2mm]
& \rho_R(x,0)=1+ 2.04204 x,  \\[2mm]
& u_R(x,0)=0.03125-0.25 x+0.75 x^2,\\[2mm]
& p_R(x,0)=10,
\end{split}
\end{equation}
which is indeed the initial conditions used as Test 2 in \cite{ct1} with slight modification, keeping only the same leading terms up to second order at $x=0$.
The initial condition (\ref{eq:test1}) has a continuous state but with discontinuous derivatives at $x=0$.
The solution for this problem contains a left-going and a right-going acoustic waves. The $t$-axis is located in the intermediate region of the two acoustic waves. This test aims at testing the accuracy of the acoustic GRP solvers.
Fig. \ref{fig:dp0} shows the solution of LGRP$_1$ solver and QGRP$_1$ solver, for each component of $U$, and the errors measured in $L_\infty$ with the rate of convergence are displayed in Table \ref{tab:acoustic}.

\begin{table}
\caption{The $L_\infty$ error of $U$ and convergence rate of the acoustic GRP solvers \hfill  }
\label{tab:acoustic}
\centering
\begin{tabular}{lllllllll}
\hline
 & \multicolumn{2}{c}{$t=0.1$} & \multicolumn{2}{c}{$t=0.05$} & \multicolumn{2}{c}{$t=0.025$}  & \multicolumn{2}{c}{$t=0.0125$} \\
Solver  &  \makebox[0pt][l]{\rule[2.0ex]{30mm}{0.5pt}}Error   & Order & \makebox[0pt][l]{\rule[2.0ex]{30mm}{0.5pt}}Error   & Order & \makebox[0pt][l]{\rule[2.0ex]{30mm}{0.5pt}}Error   & Order & \makebox[0pt][l]{\rule[2.0ex]{30mm}{0.5pt}}Error  & Order    \\

\hline

 LGRP$_{1}$ & 2.420e+0 & -- & 4.407e-1 & 2.46 & 9.592e-2 & 2.20 & 2.251e-02 & 2.09  \\

 QGRP$_{1}$ & 1.011e+0 & -- & 9.861e-2 & 3.36 & 1.127e-2 & 3.13  & 1.439e-3 & 2.97  \\

                               \hline

\end{tabular}
\end{table}

\subsubsection{Nonsonic case with jump in initial state}

In this subsection, we test the GRP solvers in the nonsonic case with initial conditions having jump in state variables. The initial conditions  are generated from  (\ref{eq:test1}) by adding a term in $p_L$ and thus generating a jump $\Delta p=(P_L(0,0),P_R(0,0))/ P_R(0,0)$ in pressure at $x=0$.

\begin{figure}
  \begin{center}
    \includegraphics[height=2.0in, width=2.5in, trim=0 0 0 0, clip]{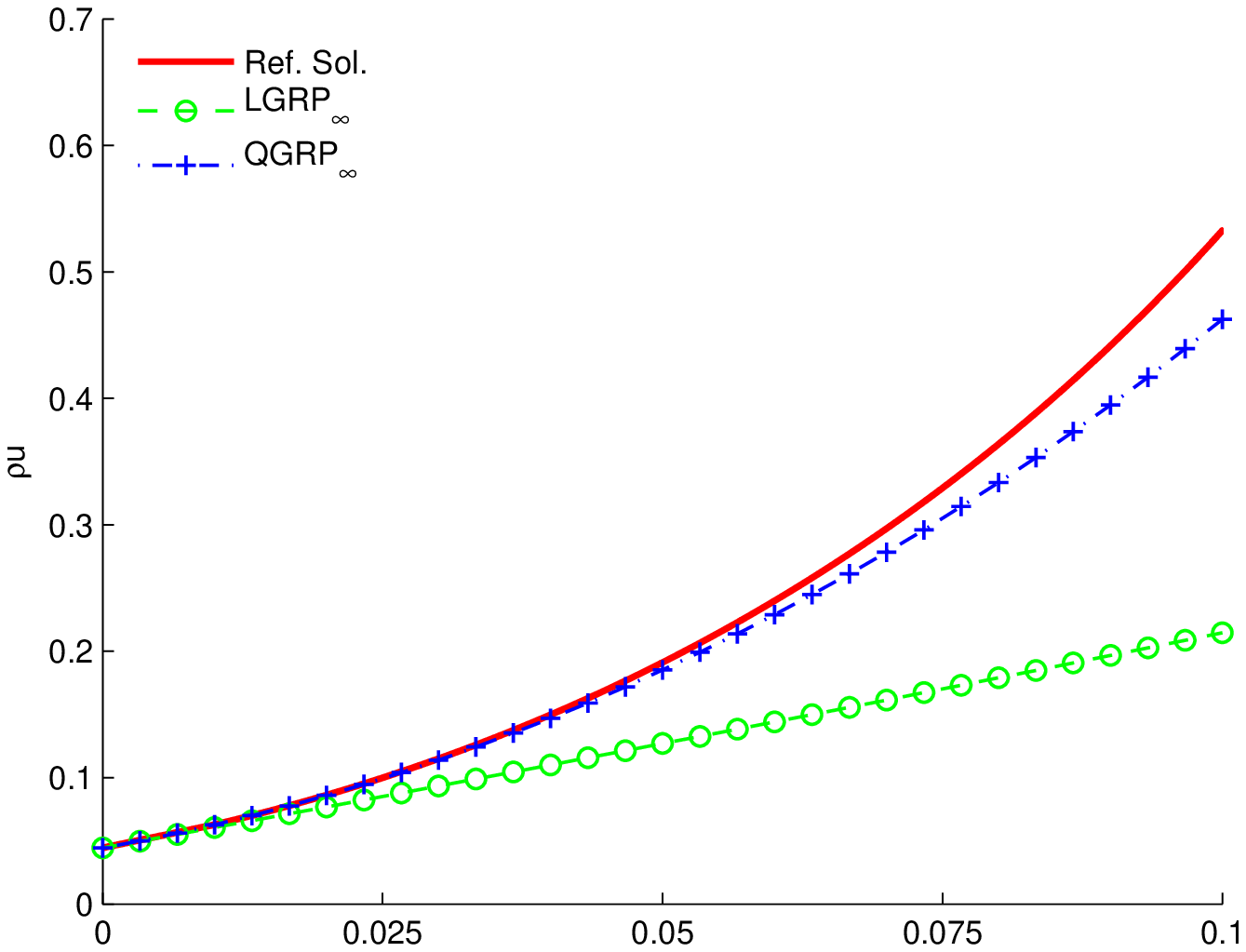}
    \includegraphics[height=2.0in, width=2.5in, trim=0 0 0 0, clip]{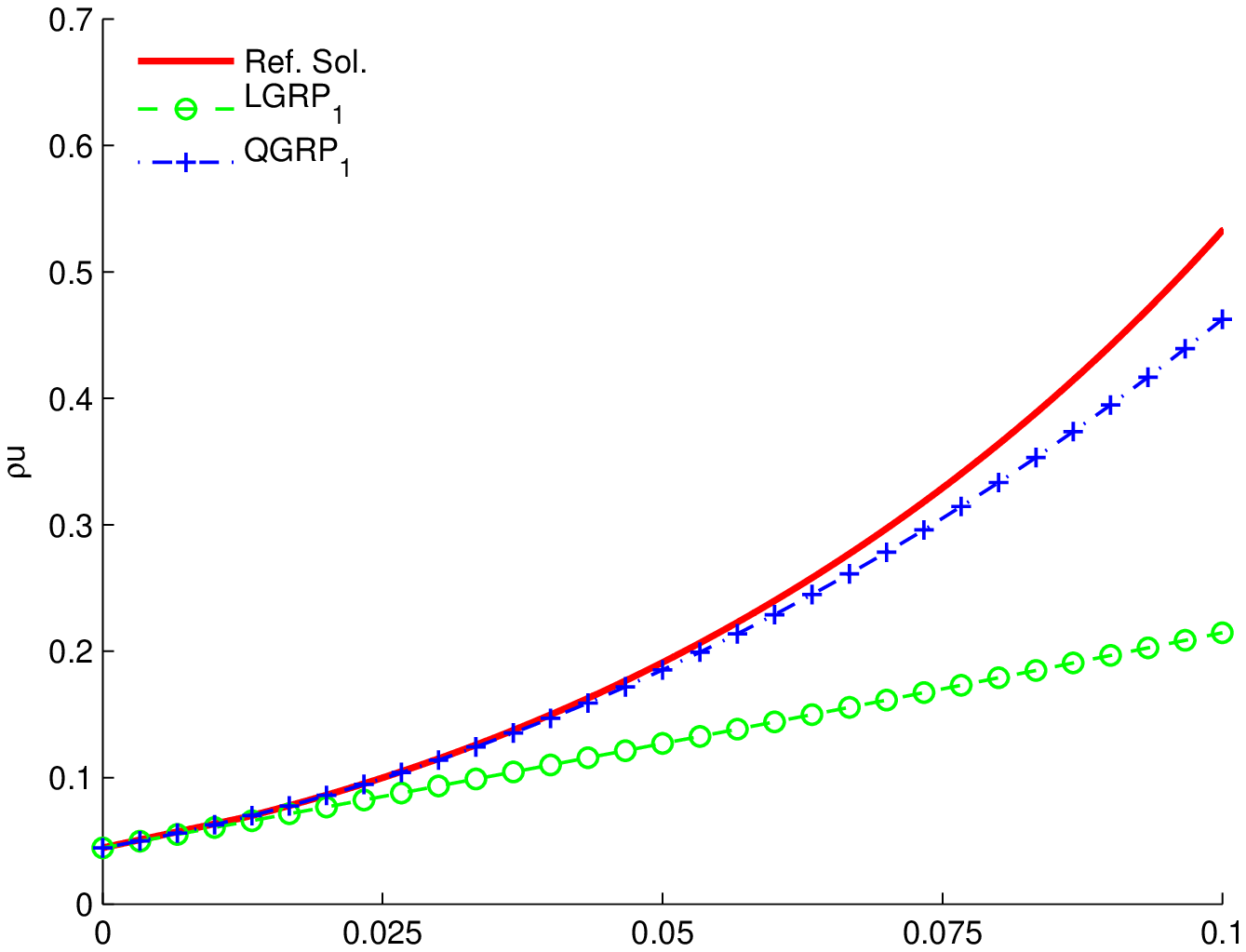}\\

    \includegraphics[height=2.0in, width=2.5in, trim=0 0 0 0, clip]{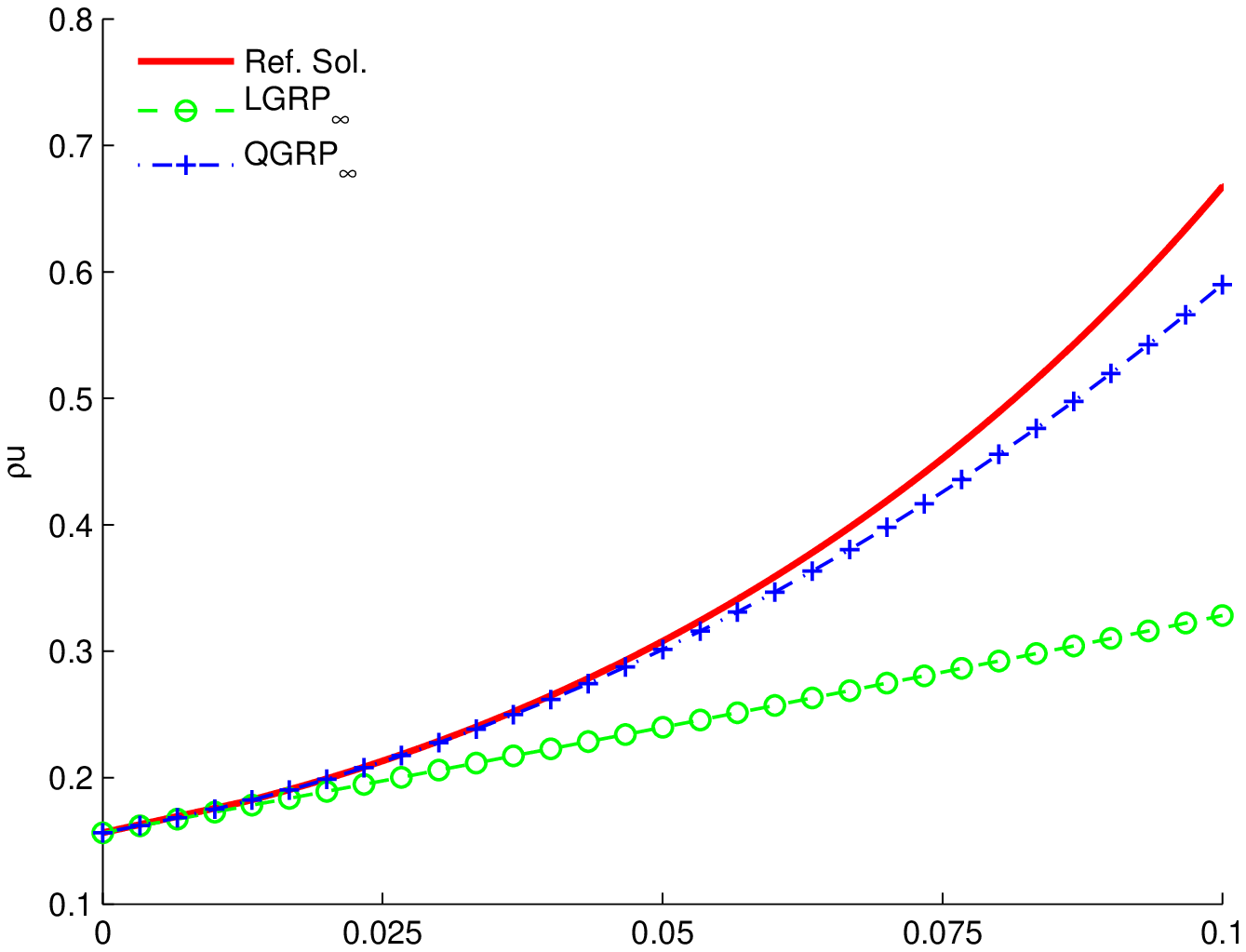}
    \includegraphics[height=2.0in, width=2.5in, trim=0 0 0 0, clip]{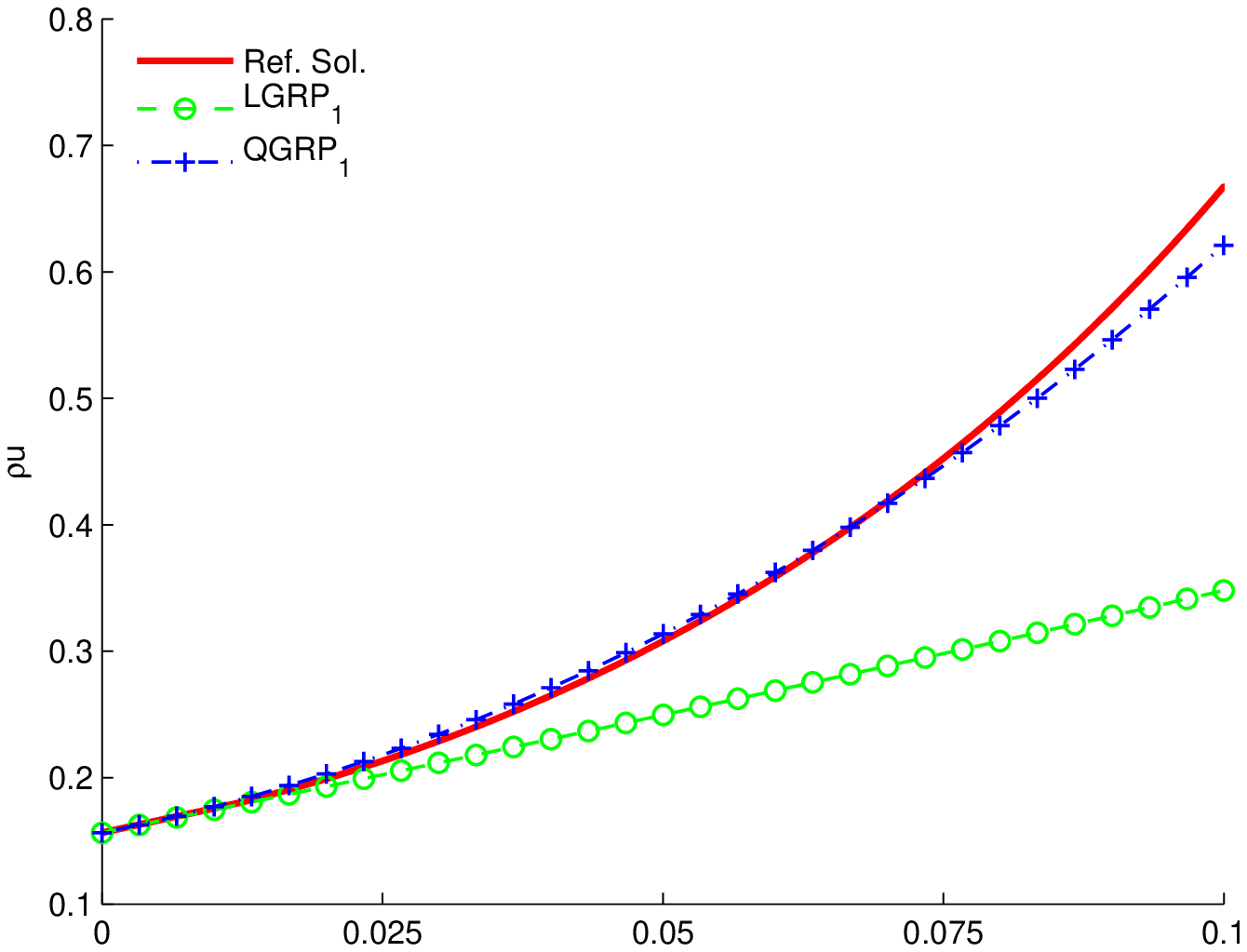} \\

    \includegraphics[height=2.0in, width=2.5in, trim=0 0 0 0, clip]{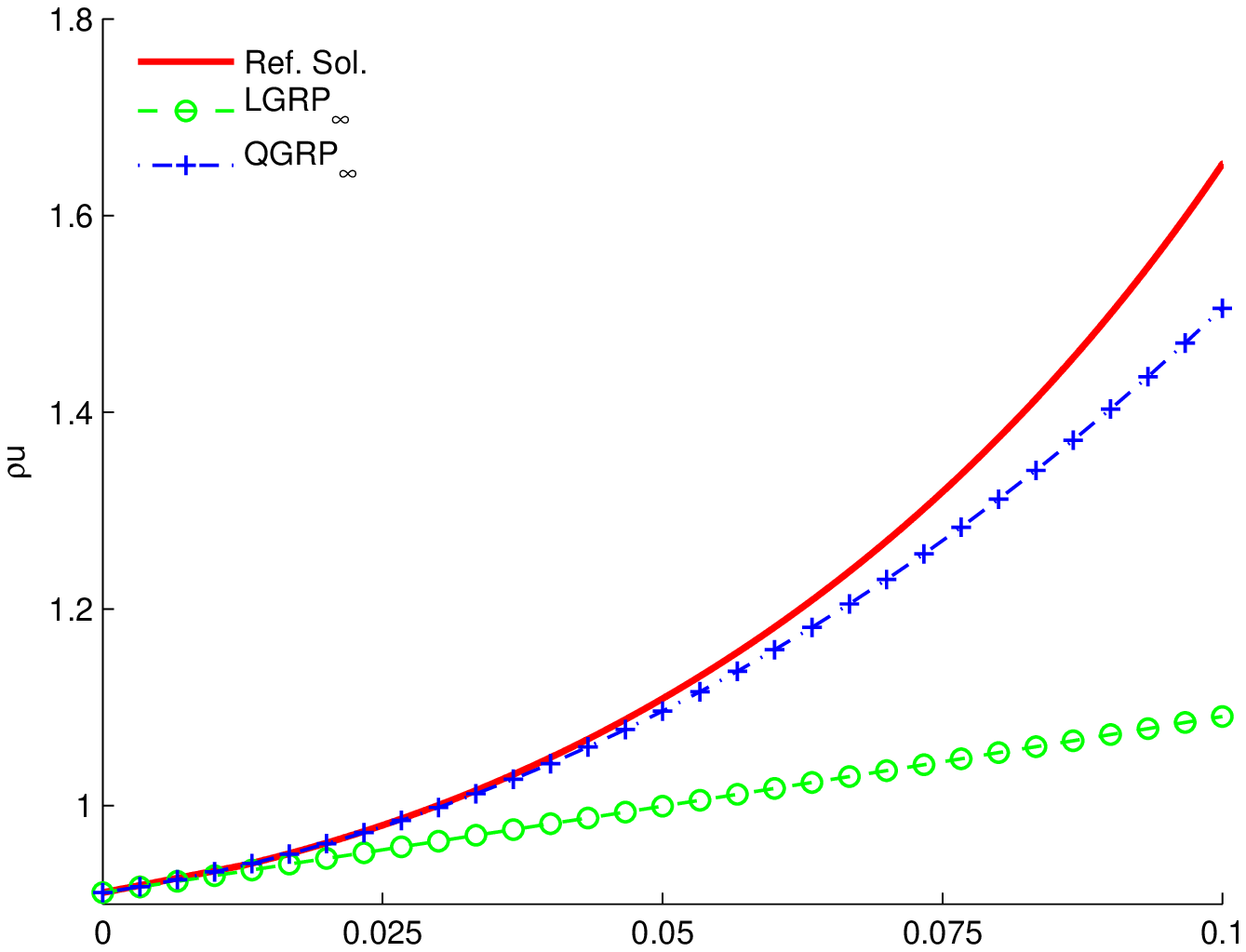}
    \includegraphics[height=2.0in, width=2.5in, trim=0 0 0 0, clip]{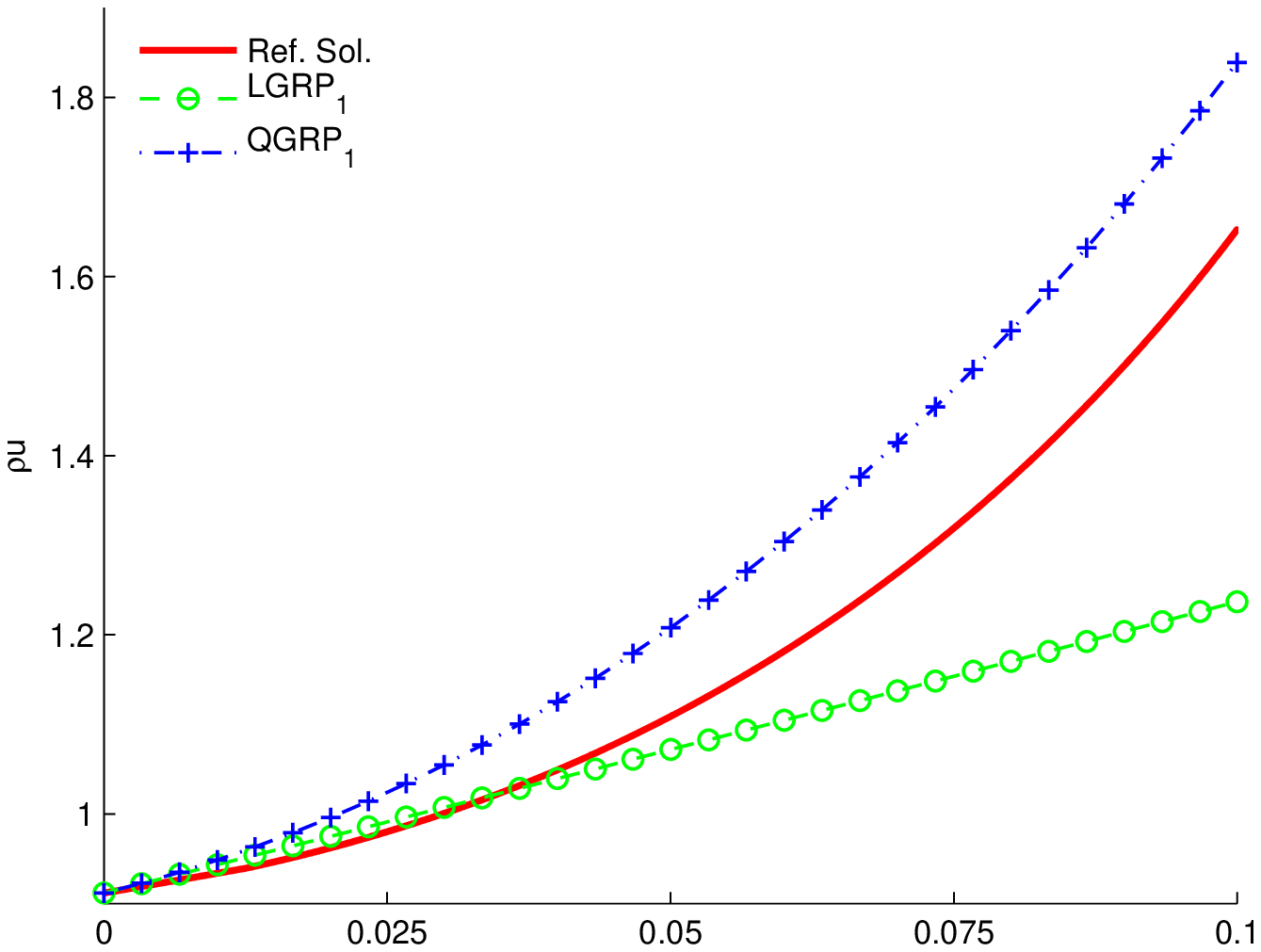}\\

    \includegraphics[height=2.0in, width=2.5in, trim=0 0 0 0, clip]{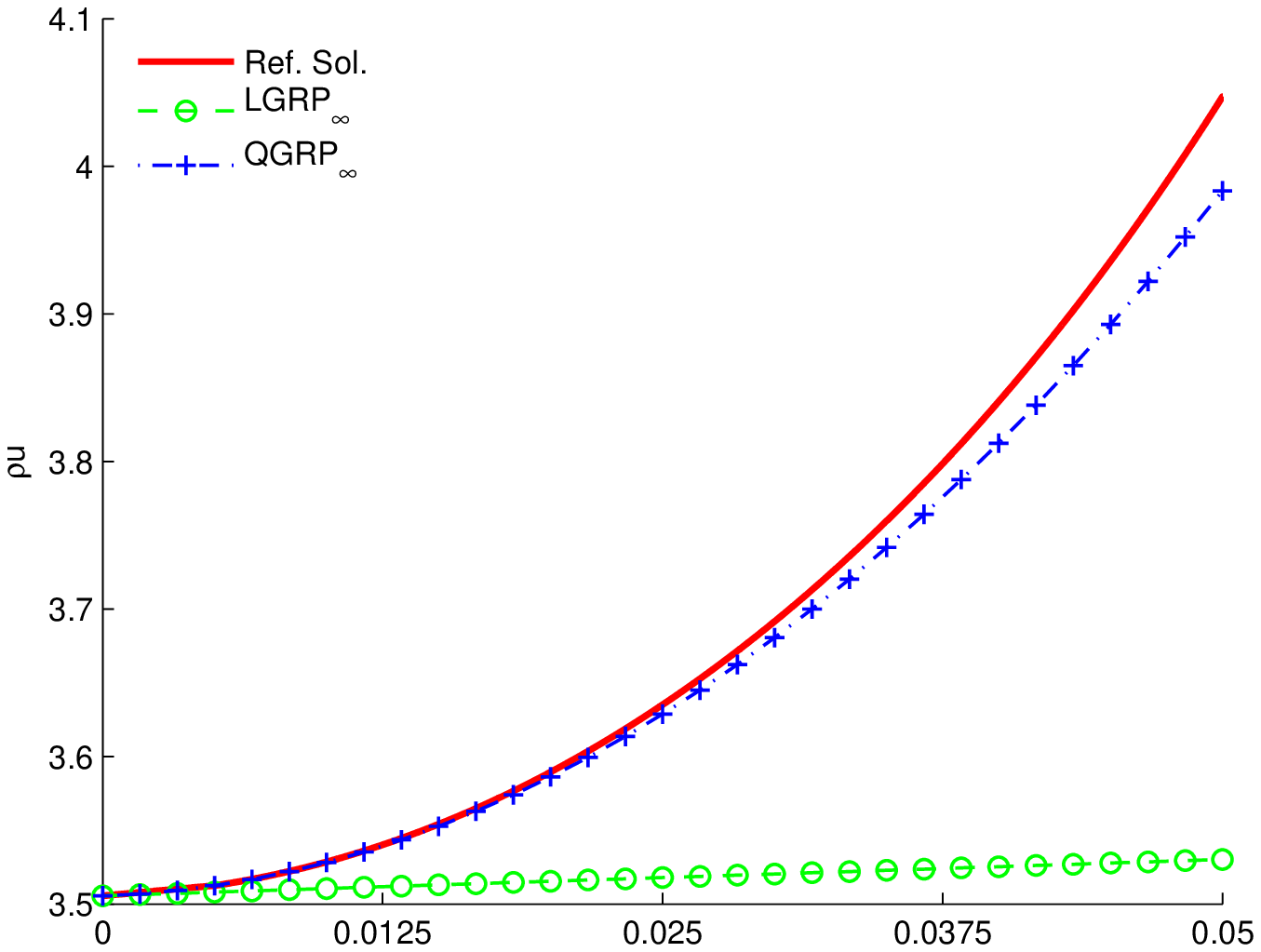}
    \includegraphics[height=2.0in, width=2.5in, trim=0 0 0 0, clip]{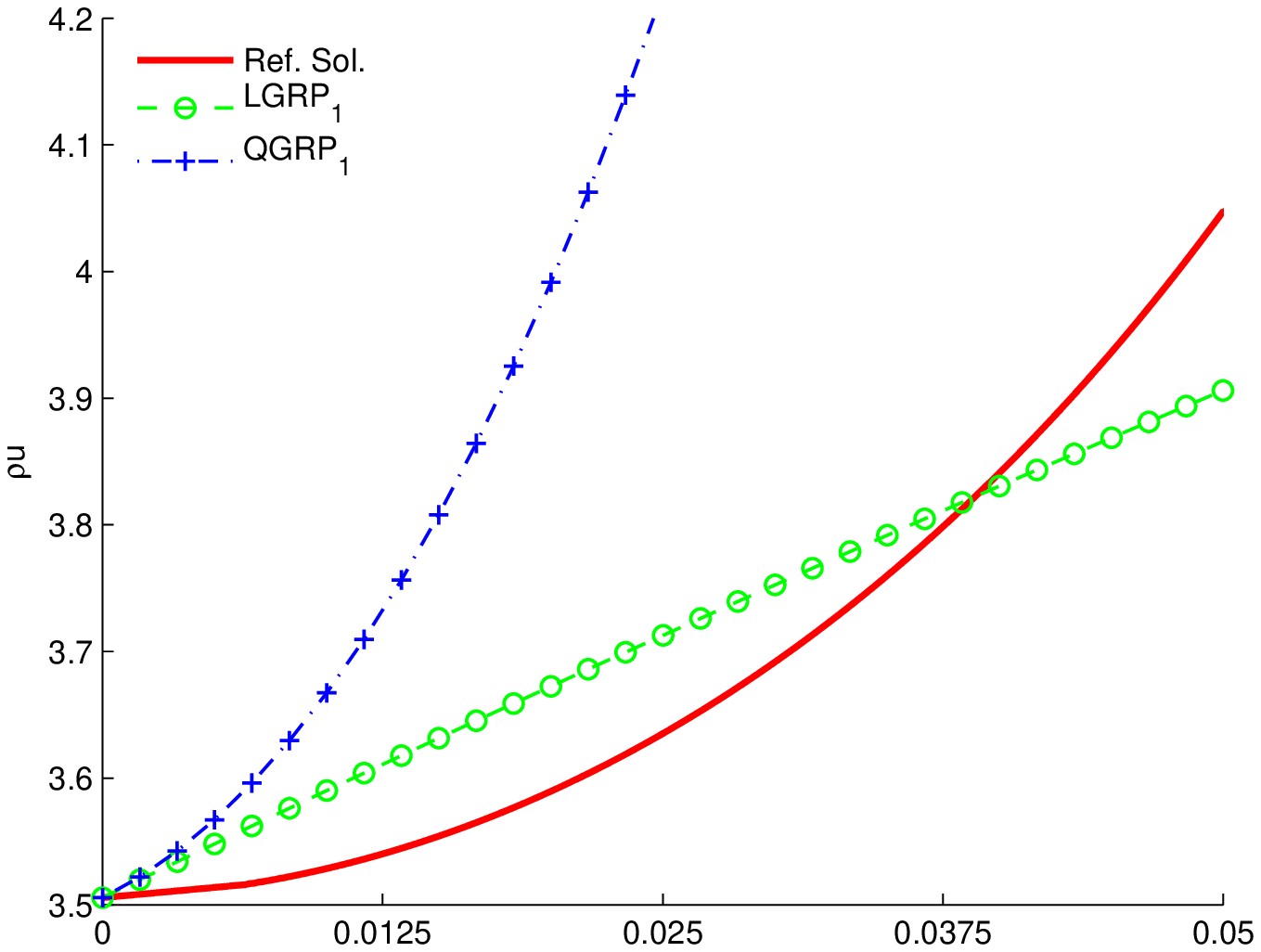}

  \caption{ Nonsonic case: Reference solution and GRP solvers based solutions. Left: LGRP$_\infty$ and LGRP$_\infty$ solver; Right: LGRP$_1$ and LGRP$_1$ solver. From top to bottom: $\Delta p=0.01$, $\Delta p=0.1$, $\Delta p=1$, $\Delta p=10$. }
  \label{fig:dp0.01-10}
  \end{center}
\end{figure}

\begin{figure}
  \begin{center}
    \includegraphics[height=2.0in, width=2.5in, trim=0 0 0 0, clip]{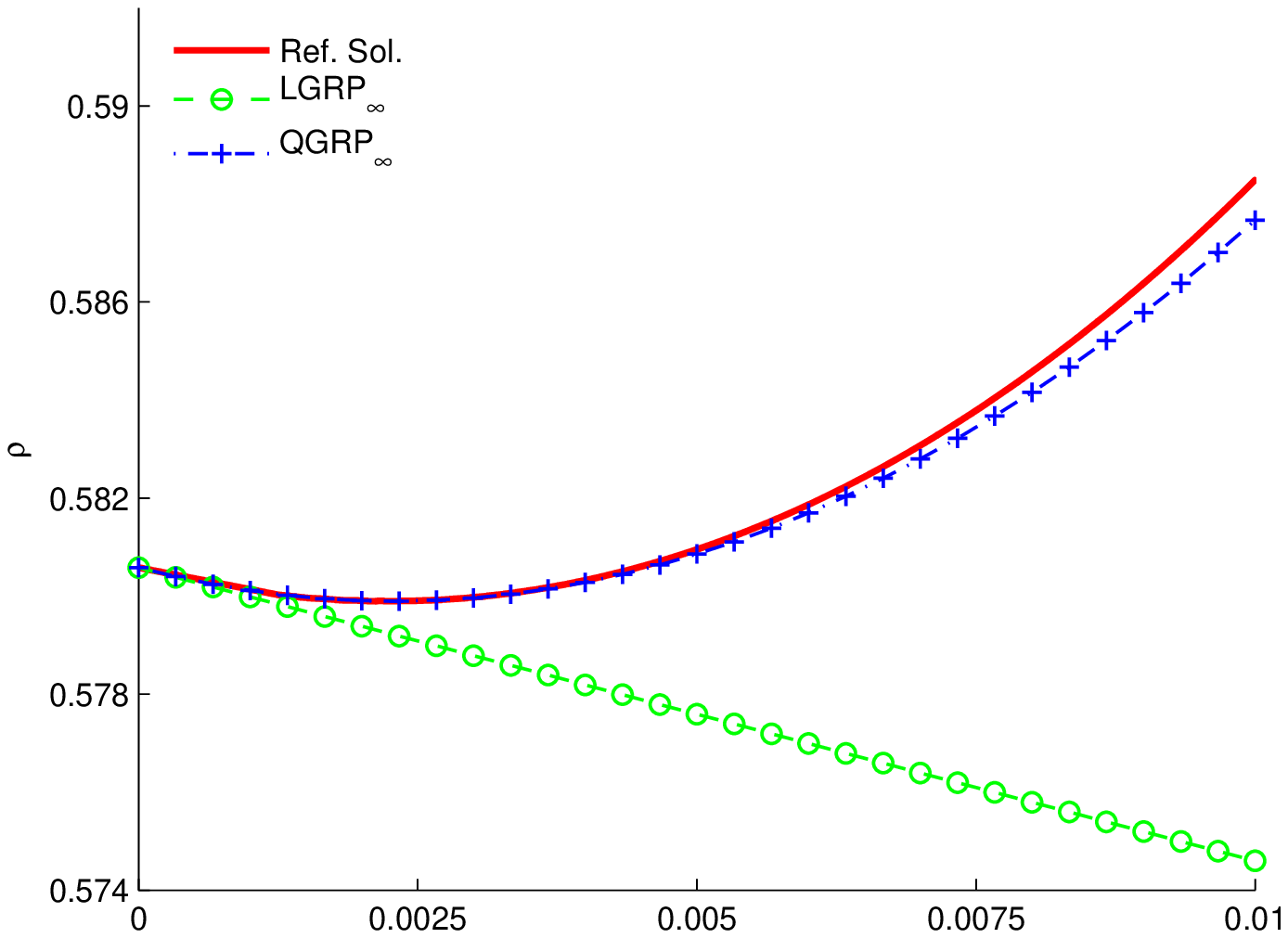}
    \includegraphics[height=2.0in, width=2.5in, trim=0 0 0 0, clip]{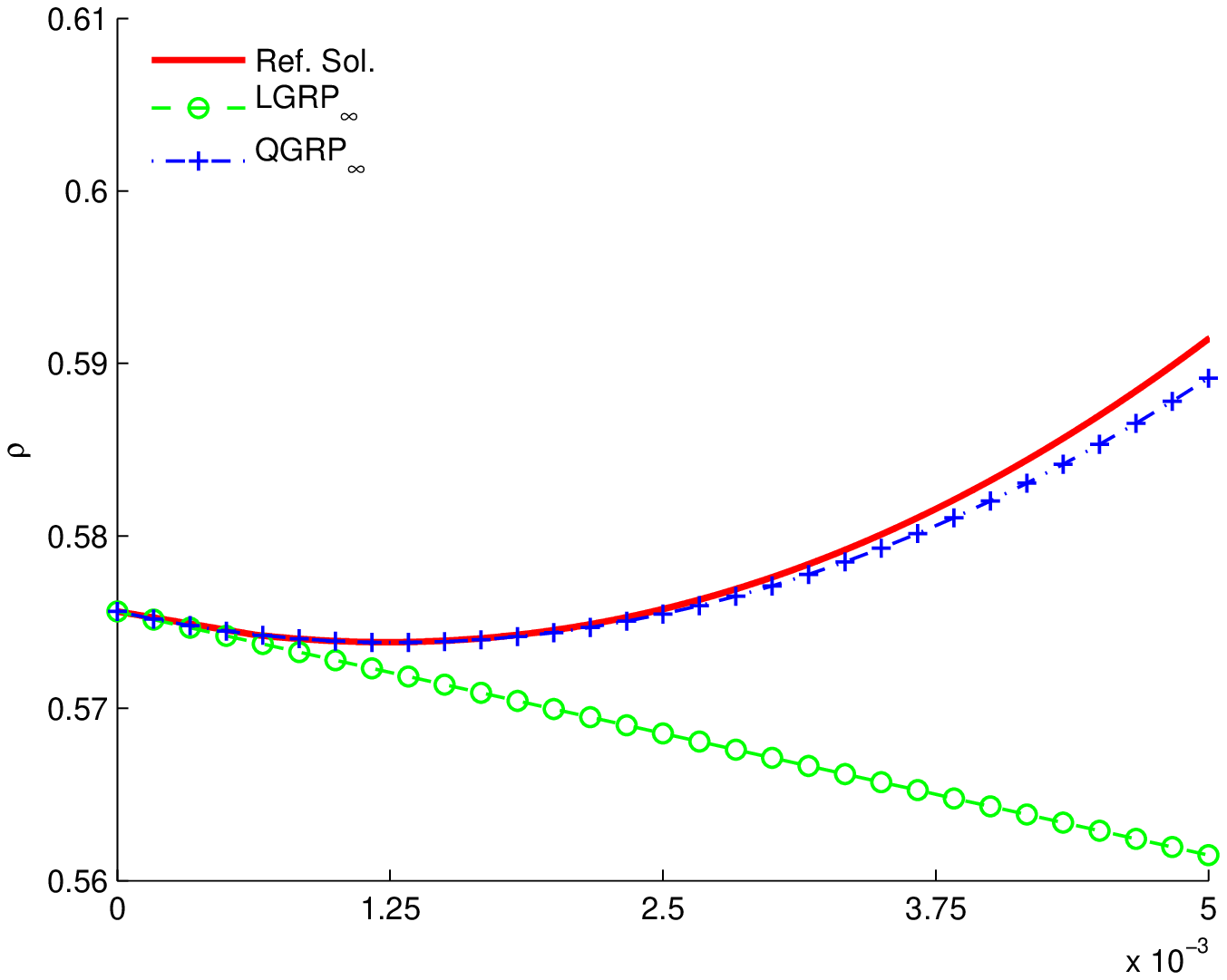}\\

    \includegraphics[height=2.0in, width=2.5in, trim=0 0 0 0, clip]{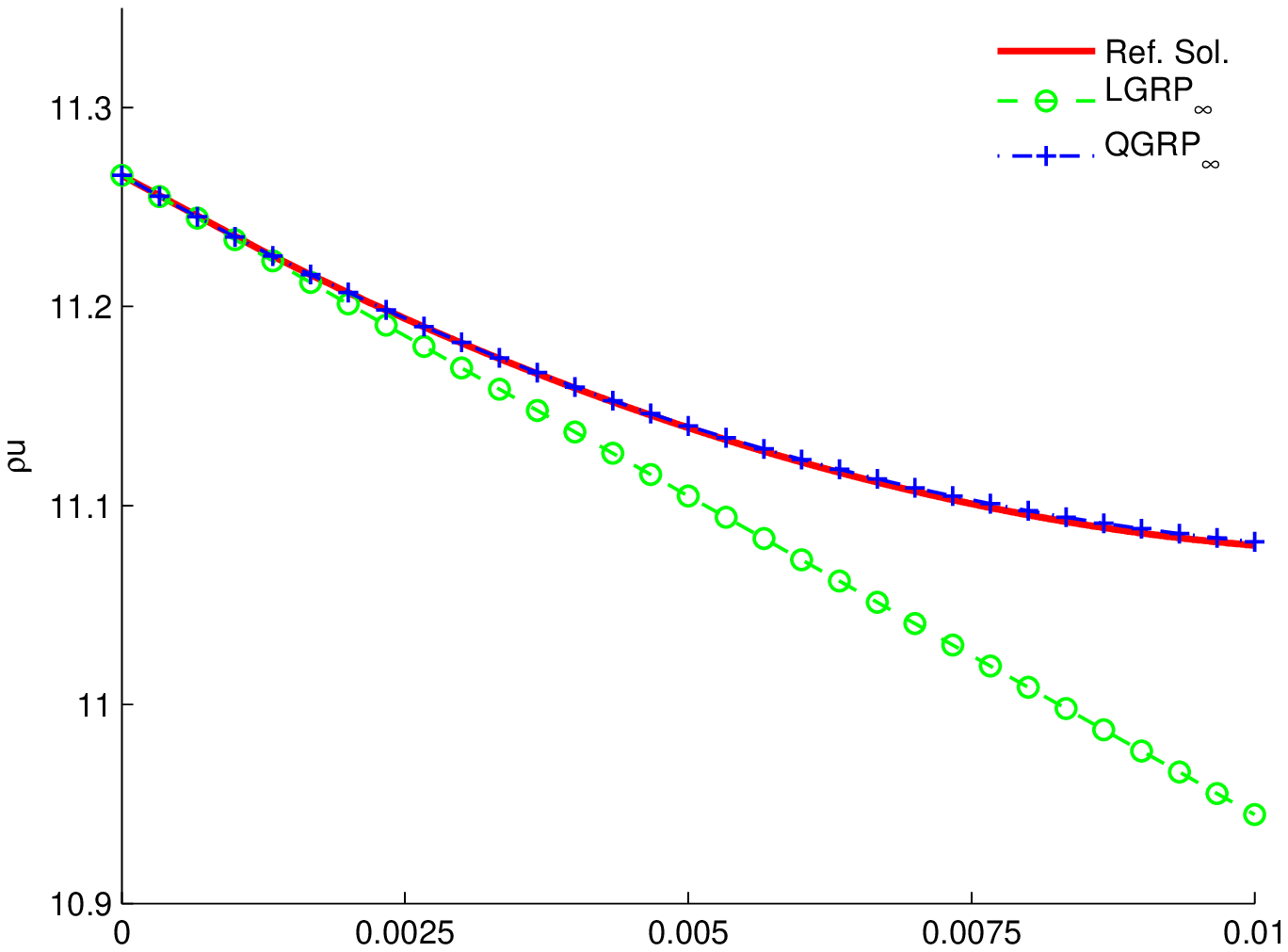}
    \includegraphics[height=2.0in, width=2.5in, trim=0 0 0 0, clip]{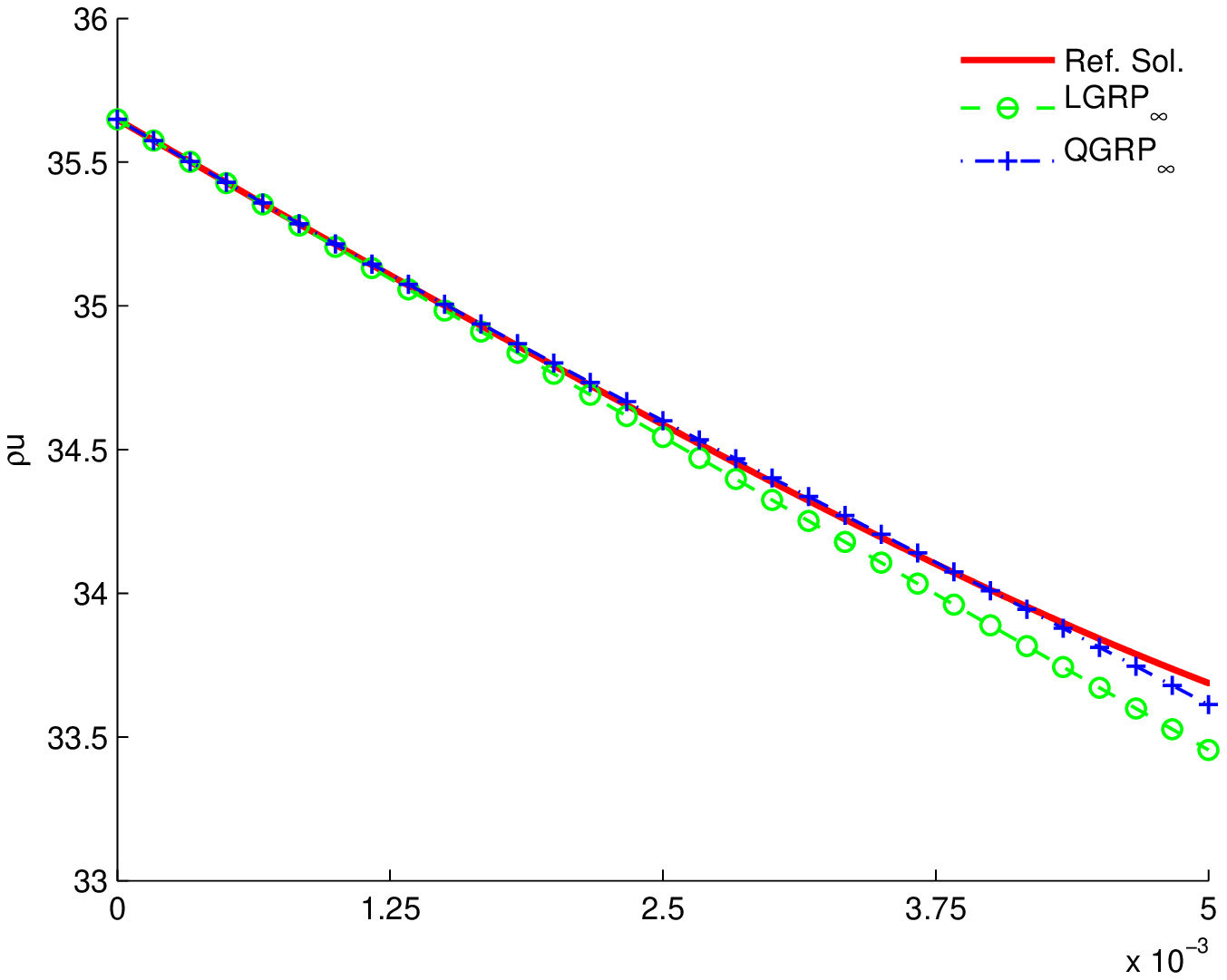}\\

    \includegraphics[height=2.0in, width=2.5in, trim=0 0 0 0, clip]{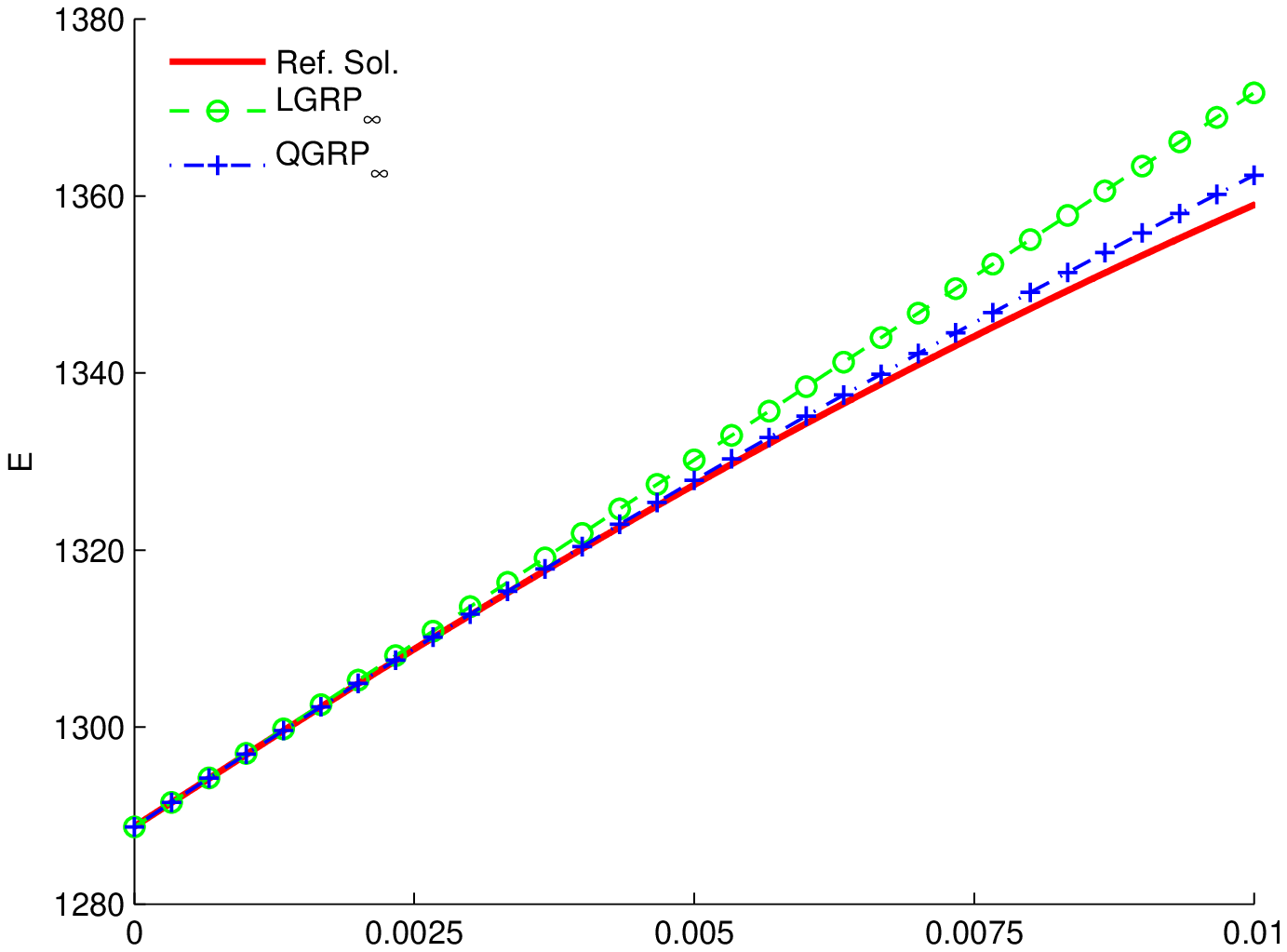}
    \includegraphics[height=2.0in, width=2.5in, trim=0 0 0 0, clip]{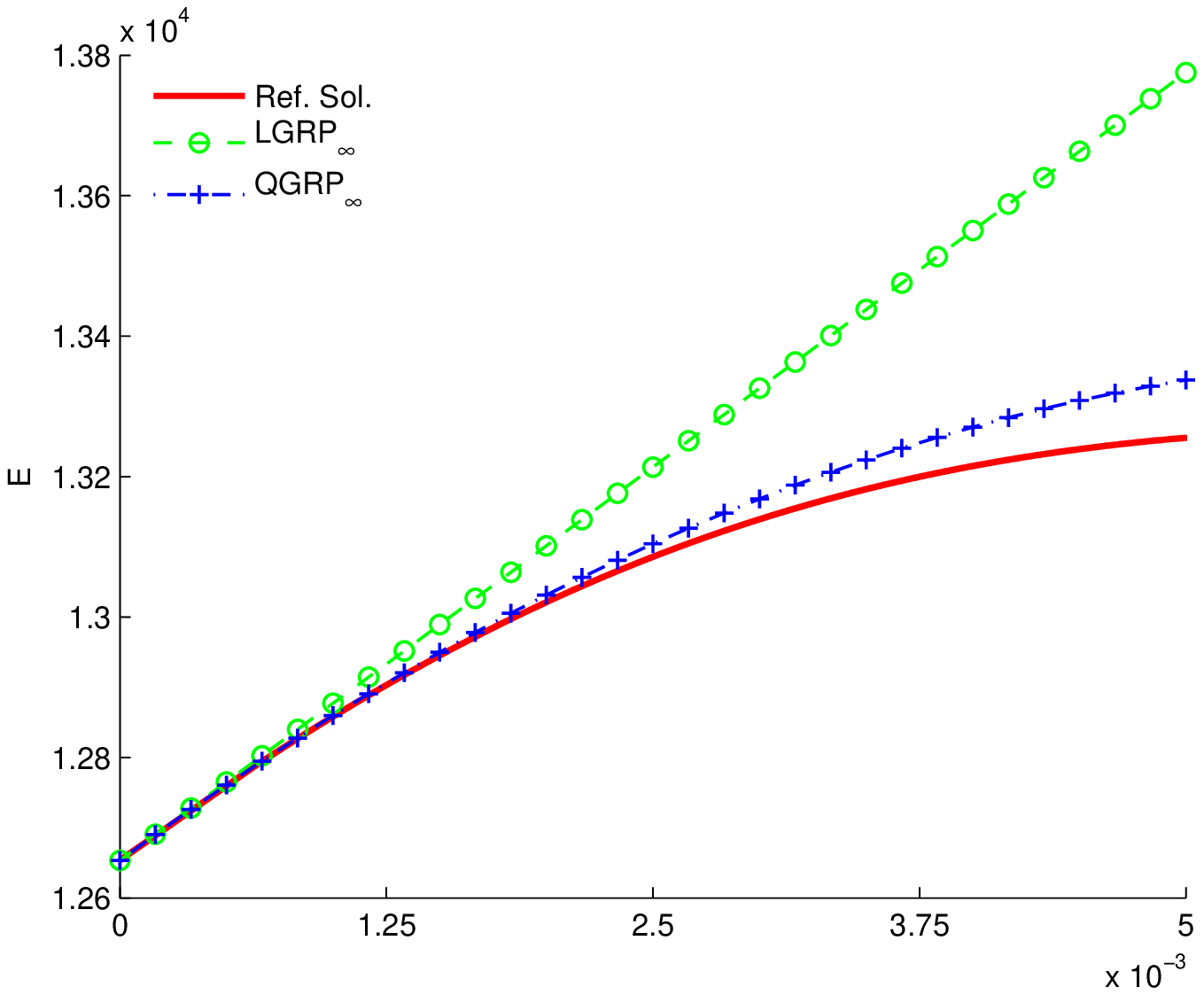}
  \caption{Nonsonic case: Reference solution and GRP solvers based solutions. Left: $\Delta p=100$; Right: $\Delta p=1000$. }
  \label{fig:dp100-1000}
  \end{center}
\end{figure}

The $L_\infty$ error of vector $U$ with the convergence rate for the LGRP$_\infty$ and QGRP$_\infty$ solver are tabulated in Tables \ref{tab:lgrp} and \ref{tab:qgrp}, respectively.
We can see that for all cases the LGRP$_\infty$ attains order two and QGPR$_\infty$ solver is essentially third order. For the QGRP solver, the decay of accuracy in some cases may be caused by the limited resolution of reference solution. Fig. \ref{fig:dp0.01-10} show the results of the acoustic LGRP$_1$ (resp. QGRP$_1$) solver in comparison with that of the LGRP$_\infty$ (resp. QGRP$_\infty$) solver with $\Delta p$ arranges from $0.01$ to $10$. When $\Delta p$ is small, the acoustic solvers do serve as good approximations of their counterpart ones.
However, as the jump $\Delta p$ increases, the performance of the acoustic solvers becomes worse. As $\Delta p=10$, the acoustic solvers give absolutely wrong initial slopes. Indeed, the behaviors of the acoustic solvers are essentially the same with that of the approximate solvers studied in \cite{ct1, ct2}.
For even larger pressure jump cases $\Delta p=100$ and $\Delta p=1000$, the solution profiles are shown in Fig. \ref{fig:dp100-1000}. We can see that, for all these test cases, the LGRP$_\infty$ and QGRP$_\infty$ solver based solutions agree well with the reference solutions.

\begin{table}
\caption{The $L_\infty$ error of $U$ and convergence rate of the LGRP$_\infty$ solvers: Nonsonic case \hfill }
 \label{tab:lgrp}
\centering
\begin{tabular}{rlllllllll}
\hline
& & \multicolumn{2}{c}{$t=t_0$} & \multicolumn{2}{c}{$t=t_0/2$} & \multicolumn{2}{c}{$t=t_0/4$}  & \multicolumn{2}{c}{$t=t_0/8$} \\
$\Delta p$ & $t_0$  &  \makebox[0pt][l]{\rule[1.8ex]{30mm}{0.5pt}}error   & Order & \makebox[0pt][l]{\rule[1.8ex]{30mm}{0.5pt}}error   & Order & \makebox[0pt][l]{\rule[1.8ex]{30mm}{0.5pt}}error   & Order & \makebox[0pt][l]{\rule[1.8ex]{30mm}{0.5pt}}error  & Order    \\
\hline

 0.01 & 0.1 & 2.456e+0  & -- & 4.478e-1 & 2.46 & 9.762e-2 & 2.20 & 2.294e-2 & 2.09  \\

 0.1 & 0.1   & 2.782e+0  & --   & 5.100e-1 & 2.45  & 1.115e-1 &  2.19  & 2.630e-2 & 2.08    \\

 1 & 0.1   &  6.012e+0 & --  & 1.128e+0 & 2.41  & 2.507e-1 & 2.17 & 5.881e-2 & 2.09  \\

 10 & 0.05   & 5.823e+0 & --  & 1.406e+0 & 2.05 & 3.515e-1 & 2.00  & 8.778e-2 & 2.00    \\

 100 & 0.01   &  1.265e+1 & --   & 2.810e+0 & 2.17 & 6.501e-1 & 2.11 & 1.600e-1 & 2.02    \\

 1000 & 0.005   & 5.201e+2 & --   & 1.277e+2 & 2.03   & 3.010e+1 & 2.08  & 7.300e+0 & 2.04  \\

                               \hline
\end{tabular}
\end{table}

\begin{table}
\caption{The $L_\infty$ error of $U$ and convergence rate of the QGRP$_\infty$ solvers: Nonsonic case \hfill }
\label{tab:qgrp}
\centering
\begin{tabular}{rlllllllll}
\hline
& & \multicolumn{2}{c}{$t=t_0$} & \multicolumn{2}{c}{$t=t_0/2$} & \multicolumn{2}{c}{$t=t_0/4$}  & \multicolumn{2}{c}{$t=t_0/8$} \\
$\Delta p$ & $t_0$  &  \makebox[0pt][l]{\rule[1.8ex]{30mm}{0.5pt}}error   & Order & \makebox[0pt][l]{\rule[1.8ex]{30mm}{0.5pt}}error   & Order & \makebox[0pt][l]{\rule[1.8ex]{30mm}{0.5pt}}error   & Order & \makebox[0pt][l]{\rule[1.8ex]{30mm}{0.5pt}}error  & Order    \\
\hline

 0.01 & 0.1 & 1.024e+0 & -- & 9.997e-2 & 3.32 & 1.157e-2 & 3.11 & 1.517e-3 & 2.93  \\

 0.1 & 0.1  & 1.141e+0 & --   & 1.113e-1 & 3.36 & 1.278e-2 & 3.12 & 1.721e-3 & 2.89    \\

 1 & 0.1   & 2.289e+0 & -- & 2.250e-1 & 3.35 & 2.714e-2 & 3.05   & 3.105e-3 &  3.13  \\

 10 & 0.05   & 1.729e-1 & --  & 2.035e-2 & 3.09 & 3.411e-3 & 2.58  & 7.979e-4 &  2.10  \\

 100 & 0.01 & 3.350e+0 & -- & 4.900e-1  & 2.77 &  7.000e-2  & 2.81 &  1.000e-2  & 2.81   \\

 1000 & 0.005  & 8.220e+1 & -- & 1.850e+1 & 2.15 & 2.800e+0 & 2.72 &  4.000e-1  & 2.81    \\
                               \hline
\end{tabular}
\end{table}

\subsubsection{Sonic case}

\begin{figure}
  \begin{center}
    \includegraphics[height=2.0in, width=2.5in, trim=0 0 0 0, clip]{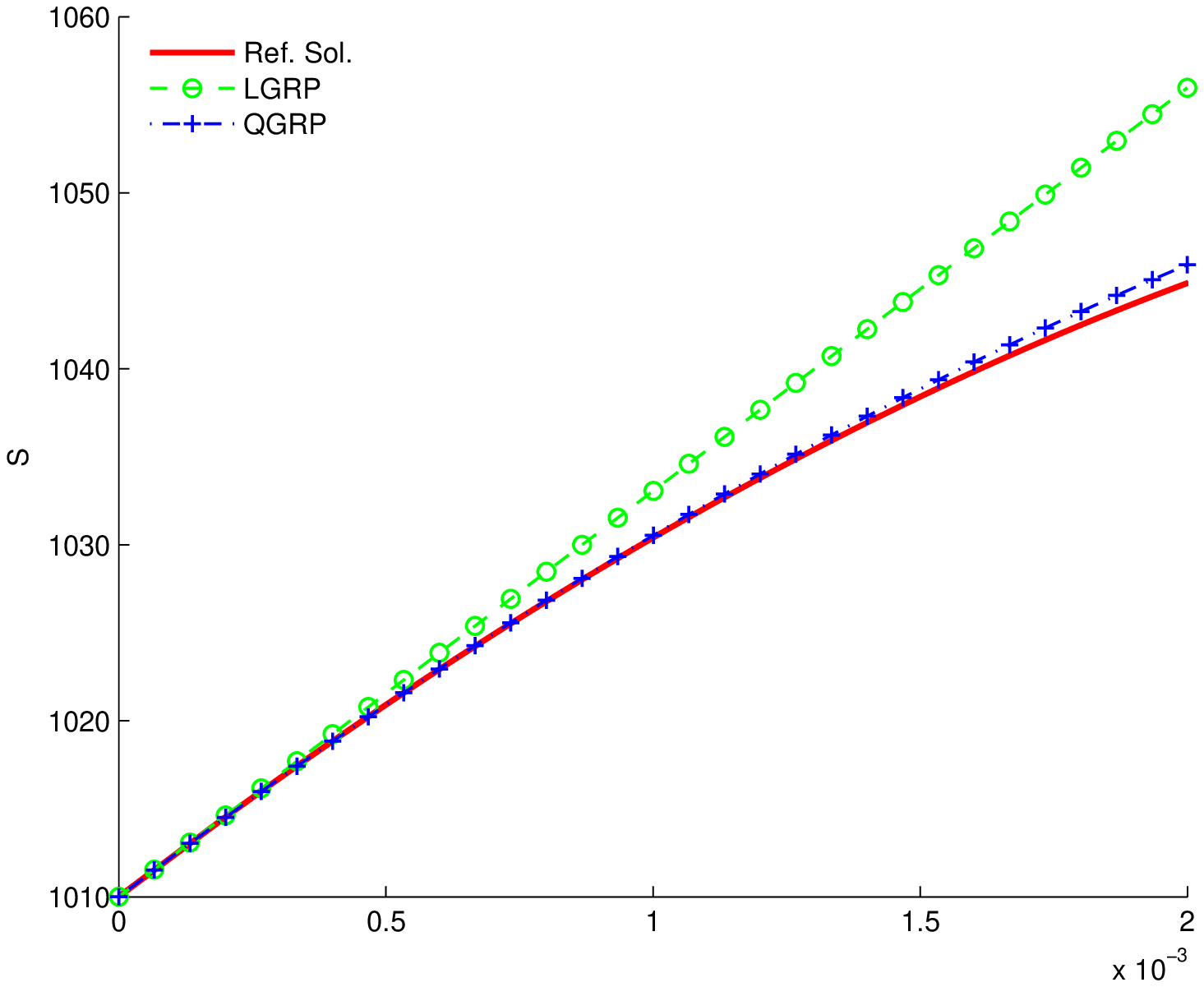}
    \includegraphics[height=2.0in, width=2.5in, trim=0 0 0 0, clip]{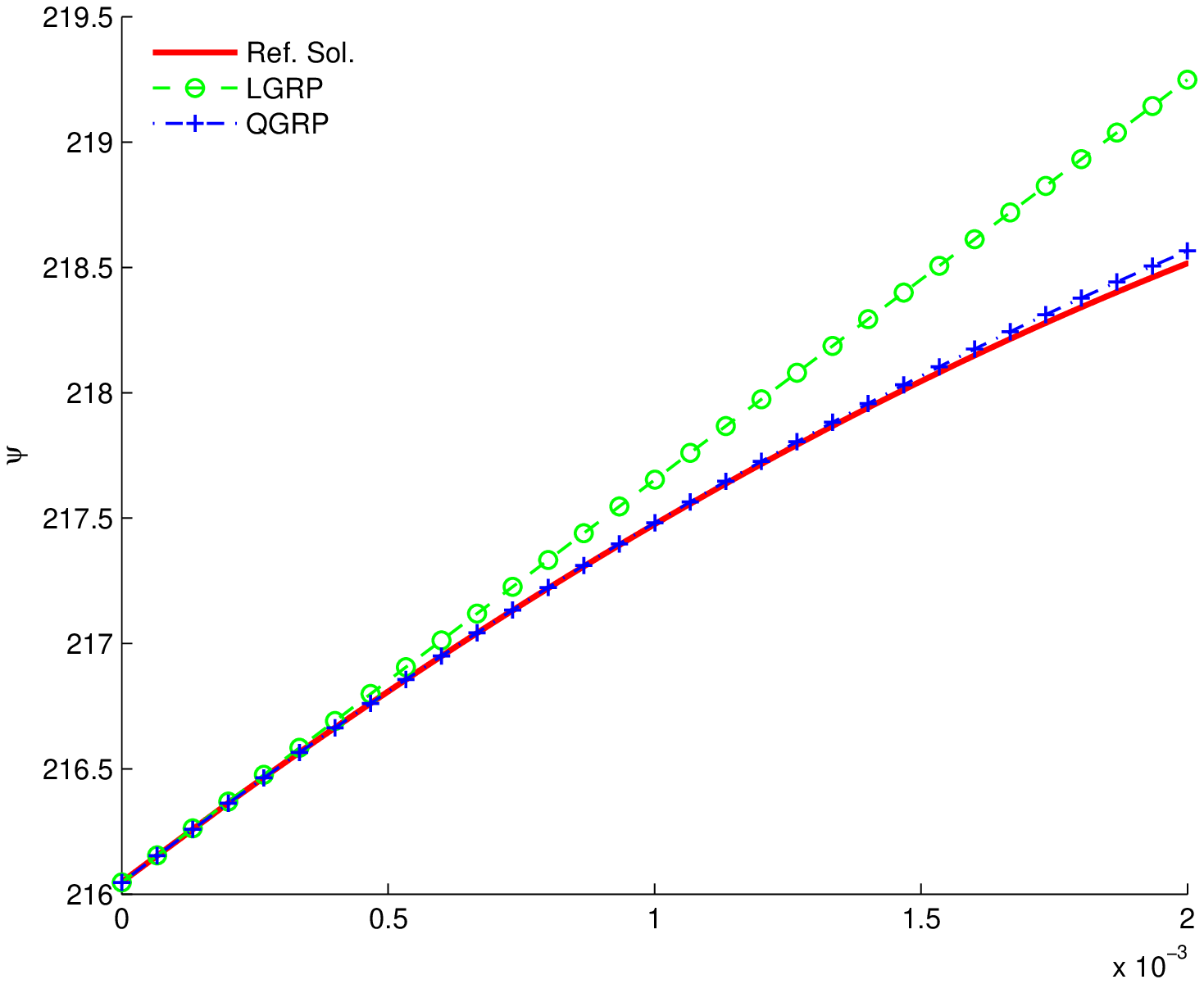}\\

    \includegraphics[height=2.0in, width=2.5in, trim=0 0 0 0, clip]{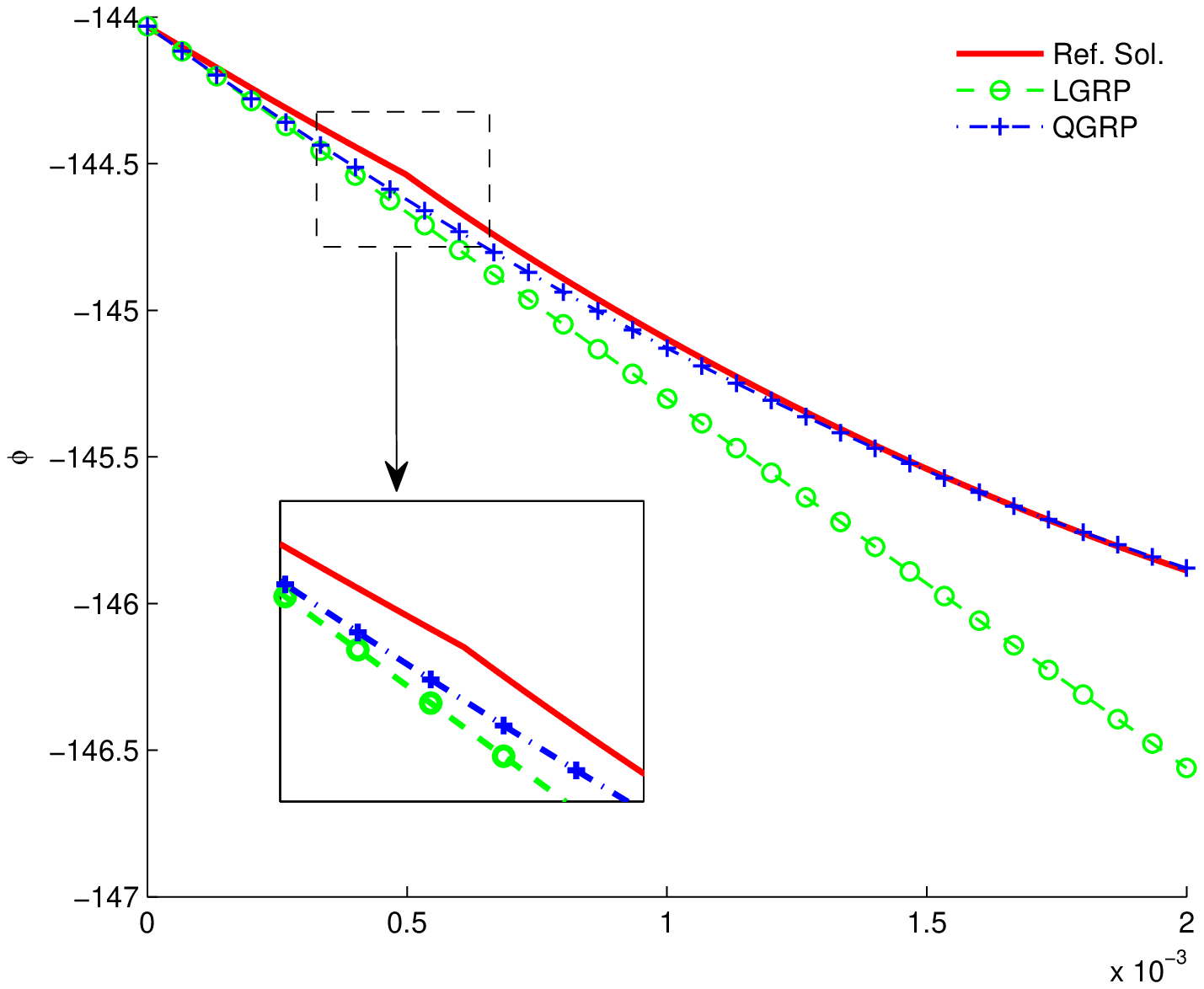}
    \includegraphics[height=2.0in, width=2.5in, trim=0 0 0 0, clip]{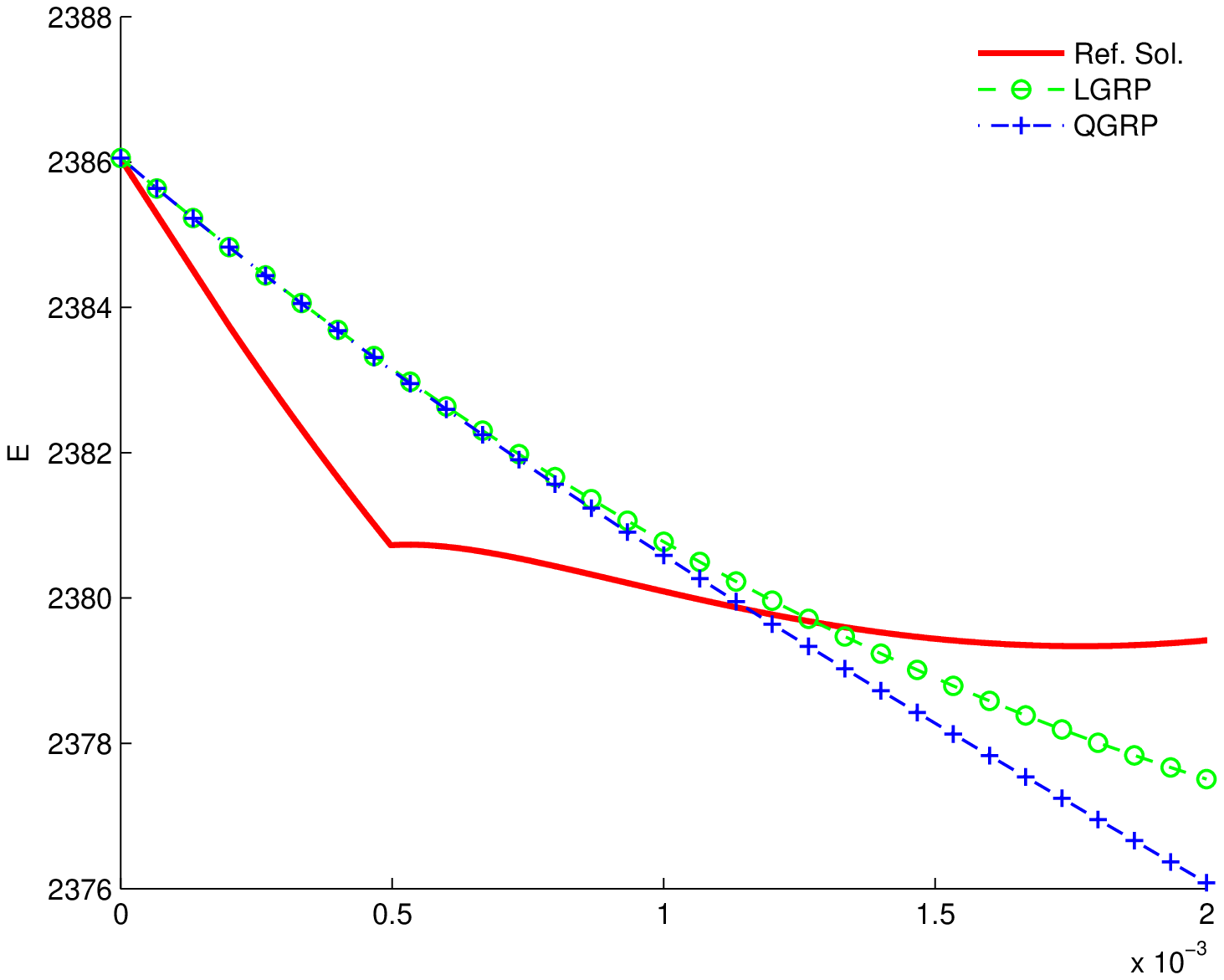}
  \caption{Sonic case: Reference solution and GRP solvers based solutions. Uniform mesh of 2.0e-4 cell size are used for computing reference solution.
}
  \label{fig:sonic-1w}
 \end{center}
\end{figure}

For the sonic case, the test problem is generated by adding $\Delta u= 28$ on the initial velocity $u(x,0)$ of the generalized Riemann problem in previous section corresponding to $p=100$.
Compared to the previous tests, it is more difficult to compute the reference solution for this case and we need to use a finer mesh with a smaller time interval. The reason is twofold. For the first, the solution is singular in the rarefaction fan, and for the second, we have observed an {\it aberration phenomenon} when computing reference solution.
The aberration phenomenon is illustrated by Fig. \ref{fig:sonic-1w}: the computed reference of $\phi$ (or $E$) exhibits a {\it weak discontinuity point} and an {\it aberration region}. However, for the GRI, the $S$ and $\psi$, such a phenomenon is not observed.
This phenomenon is different from the afore-mentioned early-time oscillation \cite{ct1}, since it is GRI-dependent.
As the mesh is refined, the weak discontinuous point converges to the singularity $(0, 0^+)$ and the numerical solution converges.

This phenomenon can be viewed as a numerical justification of the fact that the second time derivative of a variable, expect for the GRI, takes infinite value at the singularity. See Section \ref{sec:sonic}.

The errors in terms of the vector $\Phi=(S, \psi, \phi)$ and the convergence rates for the GRP solvers are displayed in Table \ref{tab:sonic}.
As suggested in Section \ref{sec:sonic}, for resolving the sonic case, we use the Newton iteration method with initial gauss $\be=0$ to solve (\ref{eq:euler-sonic-root}). Here, for the tolerance $TOL = 1.0e -7$, the number of iterations required for convergence is no more than three.

\begin{table}
\caption{The $L_\infty$ error of $\Phi$ and convergence rate of the GRP$_\infty$ solvers: Sonic case \hfill }
 \label{tab:sonic}
\centering
\begin{tabular}{lllllllll}
\hline
             & \multicolumn{2}{c}{$t=t_0$} & \multicolumn{2}{c}{$t=2/3 t_0$} & \multicolumn{2}{c}{$t=t_0/2$}  & \multicolumn{2}{c}{$t=1/3 t_0$} \\
Solver  &  \makebox[0pt][l]{\rule[2.0ex]{30mm}{0.5pt}}Error   & Order & \makebox[0pt][l]{\rule[2.0ex]{30mm}{0.5pt}}Error   & Order & \makebox[0pt][l]{\rule[2.0ex]{30mm}{0.5pt}}Error   & Order & \makebox[0pt][l]{\rule[2.0ex]{30mm}{0.5pt}}Error  & Order    \\

\hline

 LGRP$_{\infty}$ & 1.114e+1 & -- & 4.809e+0 & 2.06 & 2.671e+0 & 2.06 & 1.165 & 2.05  \\

 QGRP$_{\infty}$ & 1.052e+0 & -- & 3.208e-1 & 3.13 & 1.304e-1 & 3.13  & 4.006e-2 & 2.91  \\

                               \hline

\end{tabular}
\end{table}

\begin{figure}
  \begin{center}
    \includegraphics[height=2.0in, width=2.5in, trim=0 0 0 0, clip]{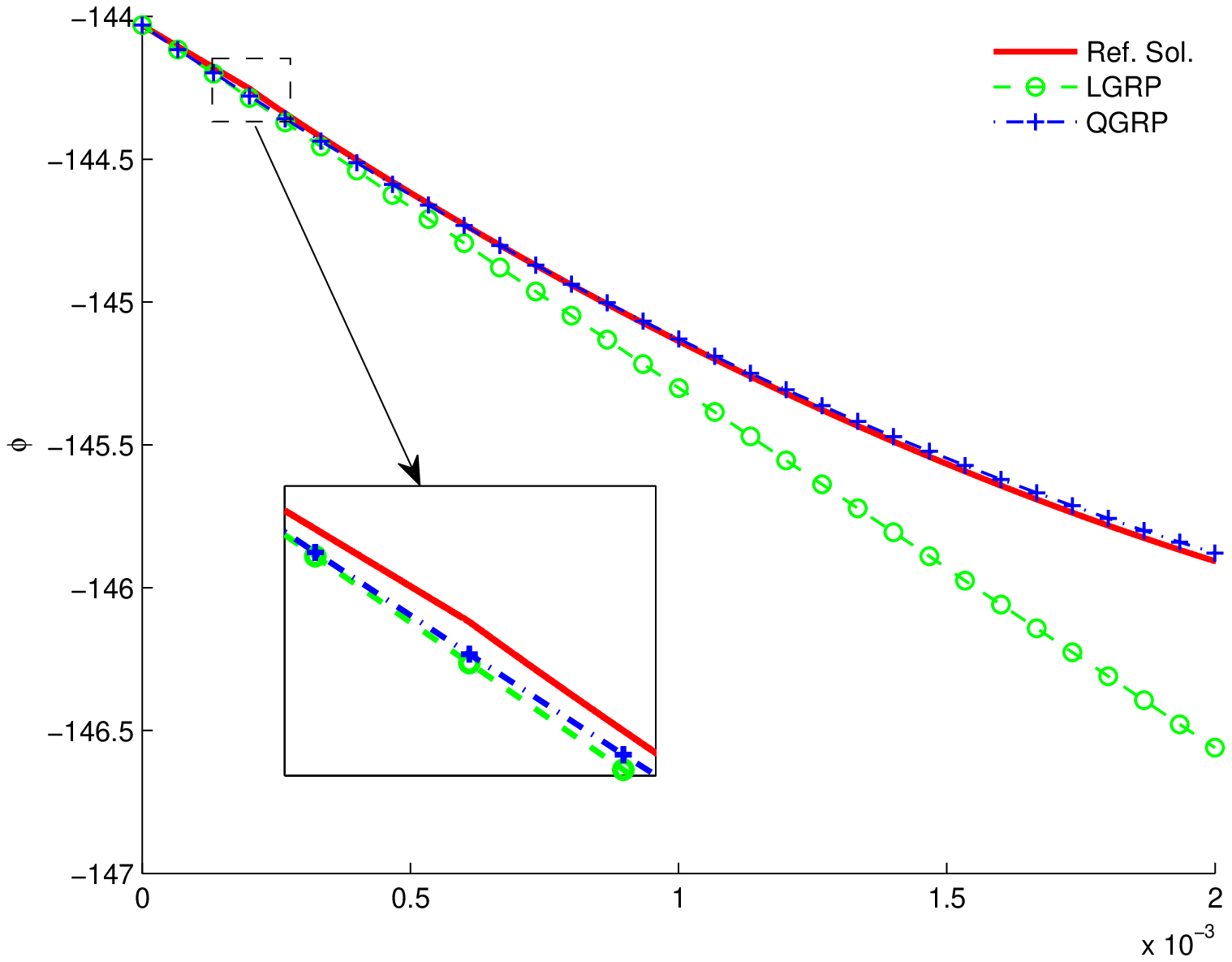}
    \includegraphics[height=2.0in, width=2.5in, trim=0 0 0 0, clip]{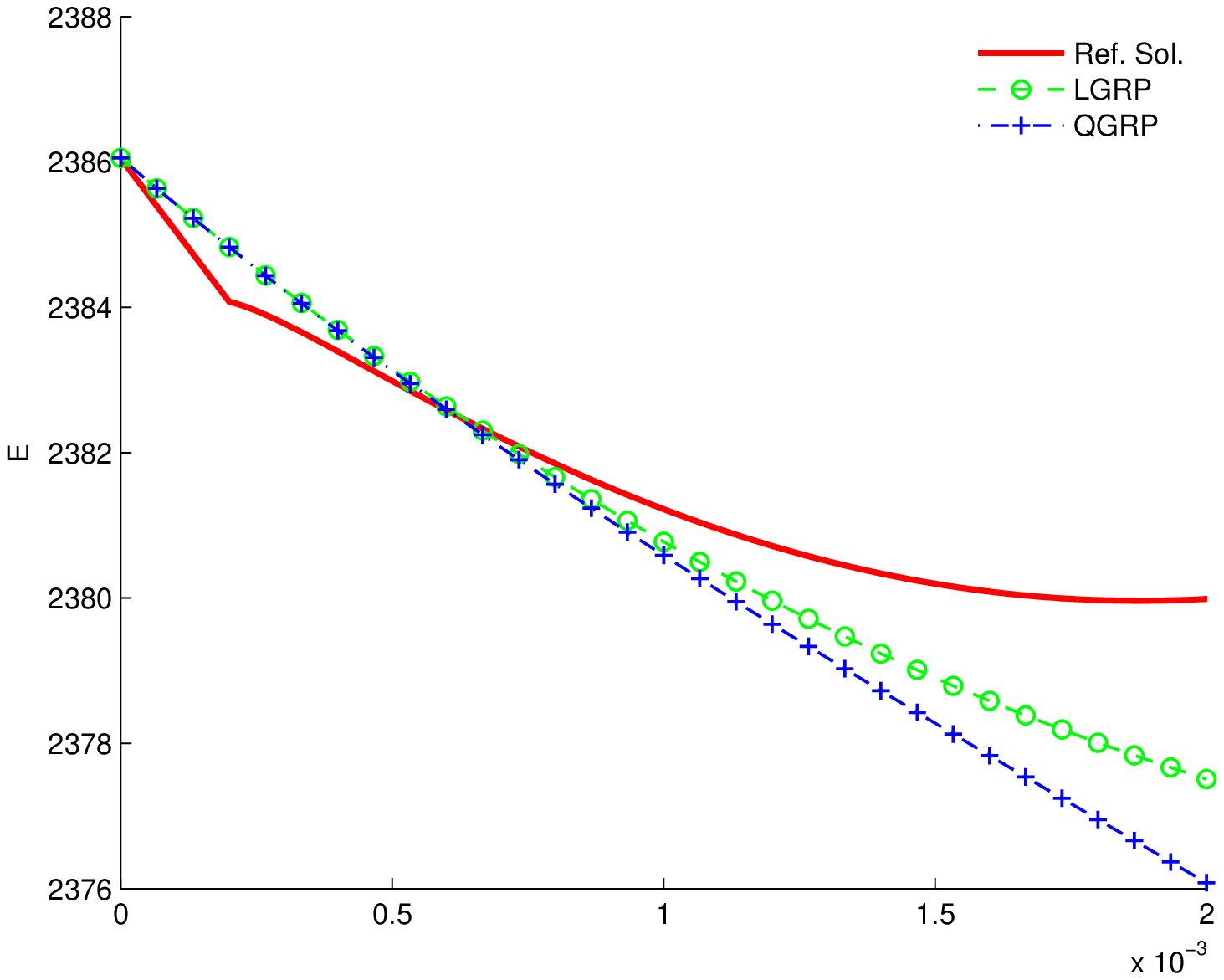}
  \caption{ Sonic case: Reference solution and GRP solvers based solutions. Uniform mesh of 2.5e-5 cell size are used for computing reference solution.
}
  \label{fig:sonic-8w}
  \end{center}
\end{figure}

\section{Numerical schemes}

In this section, we turn using GRP solvers to construct one step high order numerical schemes, namely, the GRP schemes.
In the introduction, we have described the process of implementing the second-order numerical scheme, where the LGRP solver provides a second-order approximation of the flux function from a piecewise linear discontinuous initial data.
The process of implementing the QGRP solver based third-order numerical scheme is similar. The differences is that we need to provide a third-order subcell data reconstruction on each time step and use two point quadrature for the integral of (\ref{eq:F_ave}) to compute the numerical flux, i.e.
\begin{align}\label{eq:F_ave2}
F_{j+1/2}=\omega_1F(U(x_{j+1/2}, \tau_1))+\omega_2F(x_{j+1/2}, \tau_2).
\end{align}
On each quadrature points $(x_{j+1/2}, \tau_i)$, the vector $U$ are calculated through (\ref{eq:tal}), wherein the $ U(0, 0^+)$, $\p_tU(0, 0^+)$ and $\p_t^2U(0, 0^+)$ are determined by solving a generalized Riemann problem on the cell interface using the QGRP solver.

In the following, we present several one-dimensional examples to test the performance of our schemes.
The uniform size meshes are used for all the test cases. For the second-order scheme, the van Leer limiter \cite{muscl2} is used to perform the linear reconstruction. For the third-order scheme, we use the same reconstruction method as in \cite{luo-xu}. In fact, we use the 5rd order WENO technique to reconstruct pointwise variables of $U$ at each cell interface.
Then based on the cell interface values and the cell averages of $U_j^n$, a third-order polynomial is constructed as the subcell flow distributions at time $t^n$.
 In the following numerical examples, the WENO reconstruction is carried out based on the characteristic decomposition \cite{shu1} and the CFL number is set to be 0.5.

For all the problems, the GRP solutions are plotted against the exact solutions. The solid lines represent the exact solution, the circles show the second-order scheme solution, while the crosses stand for the third-order scheme solution.

\subsection{Sod problem}

\begin{figure}
  \begin{center}
    \includegraphics[height=2.0in, width=2.6in, trim=0 0 0 0, clip]{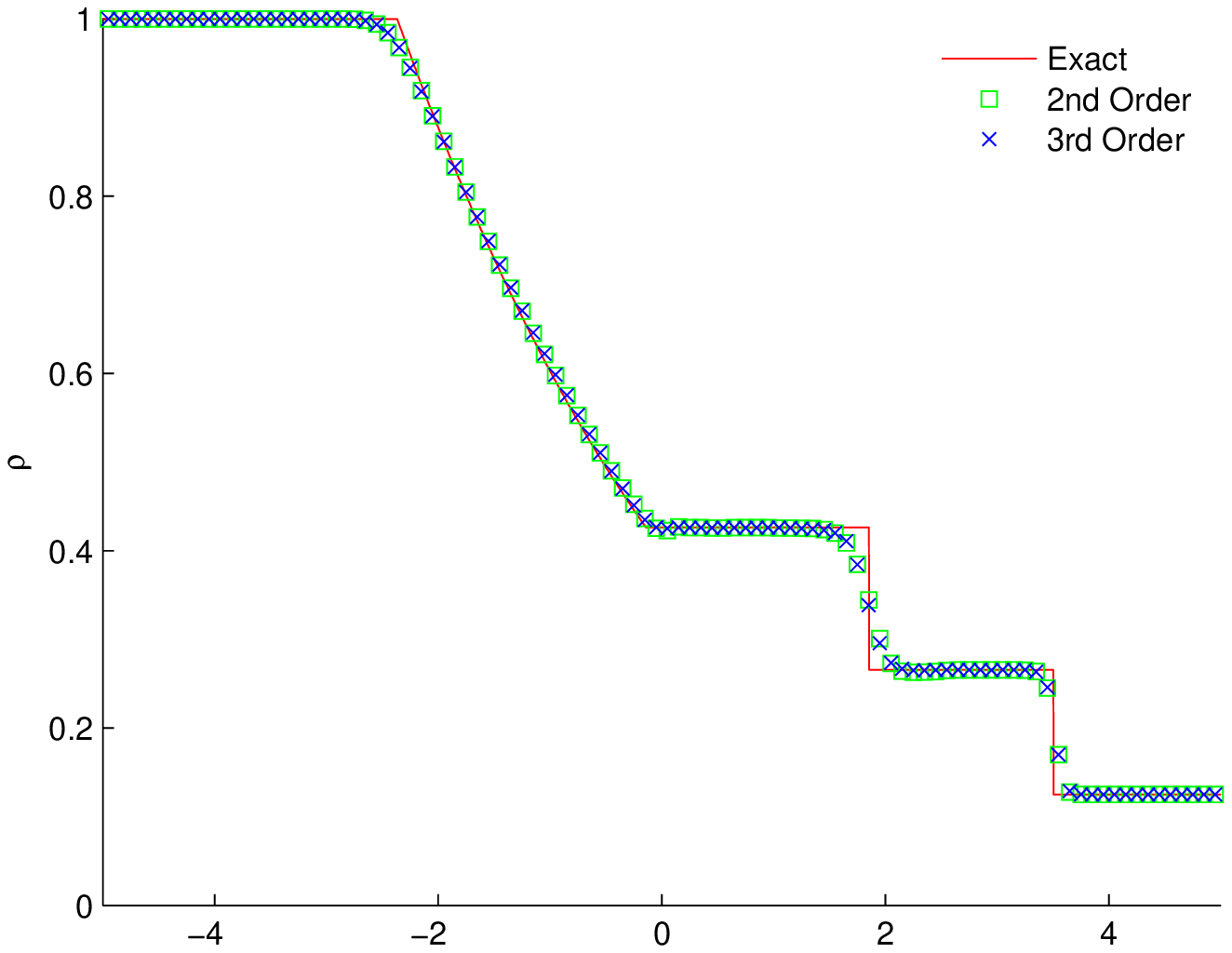}
    \includegraphics[height=2.0in, width=2.6in, trim=0 0 0 0, clip]{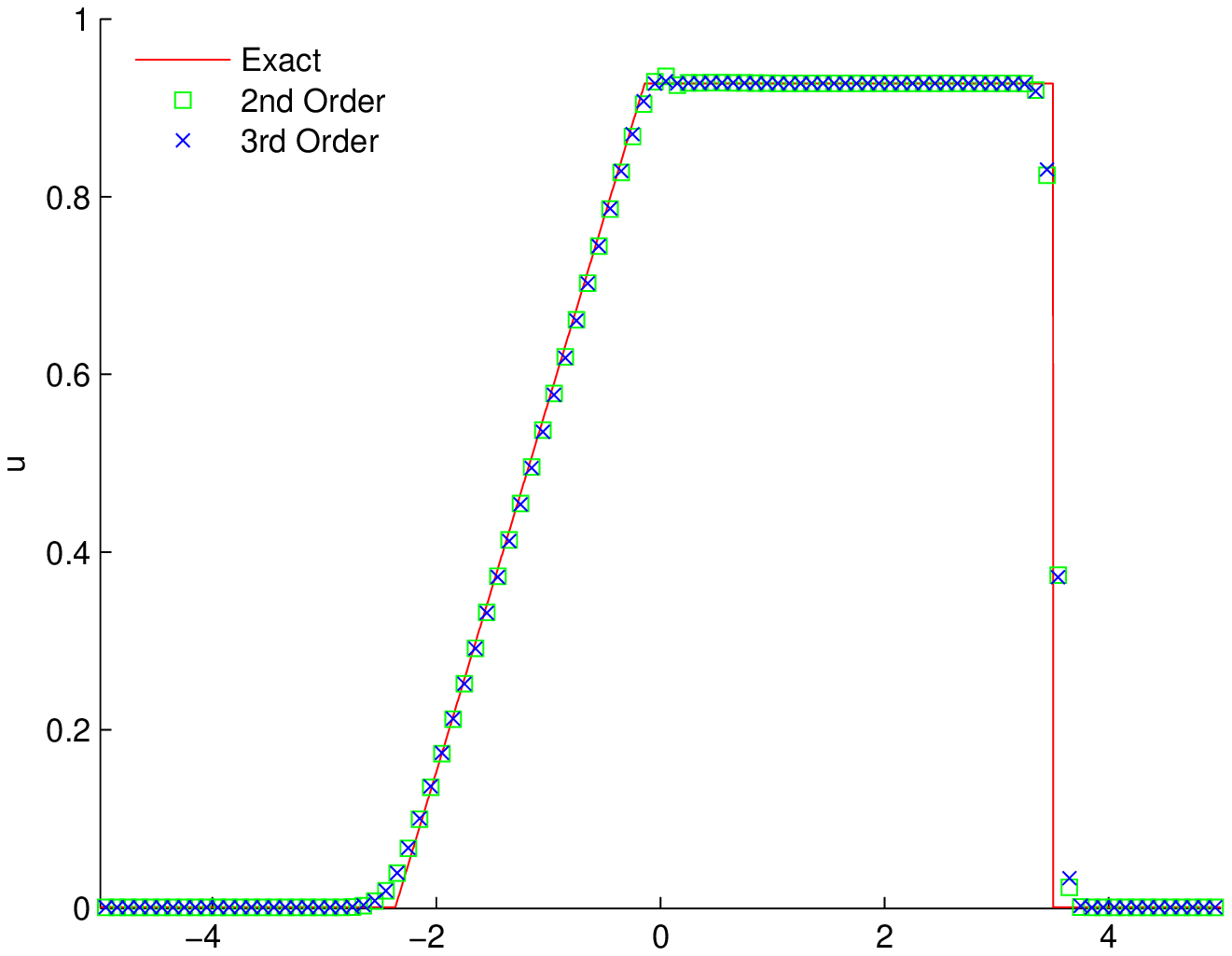}\\

    \includegraphics[height=2.0in, width=2.6in, trim=0 0 0 0, clip]{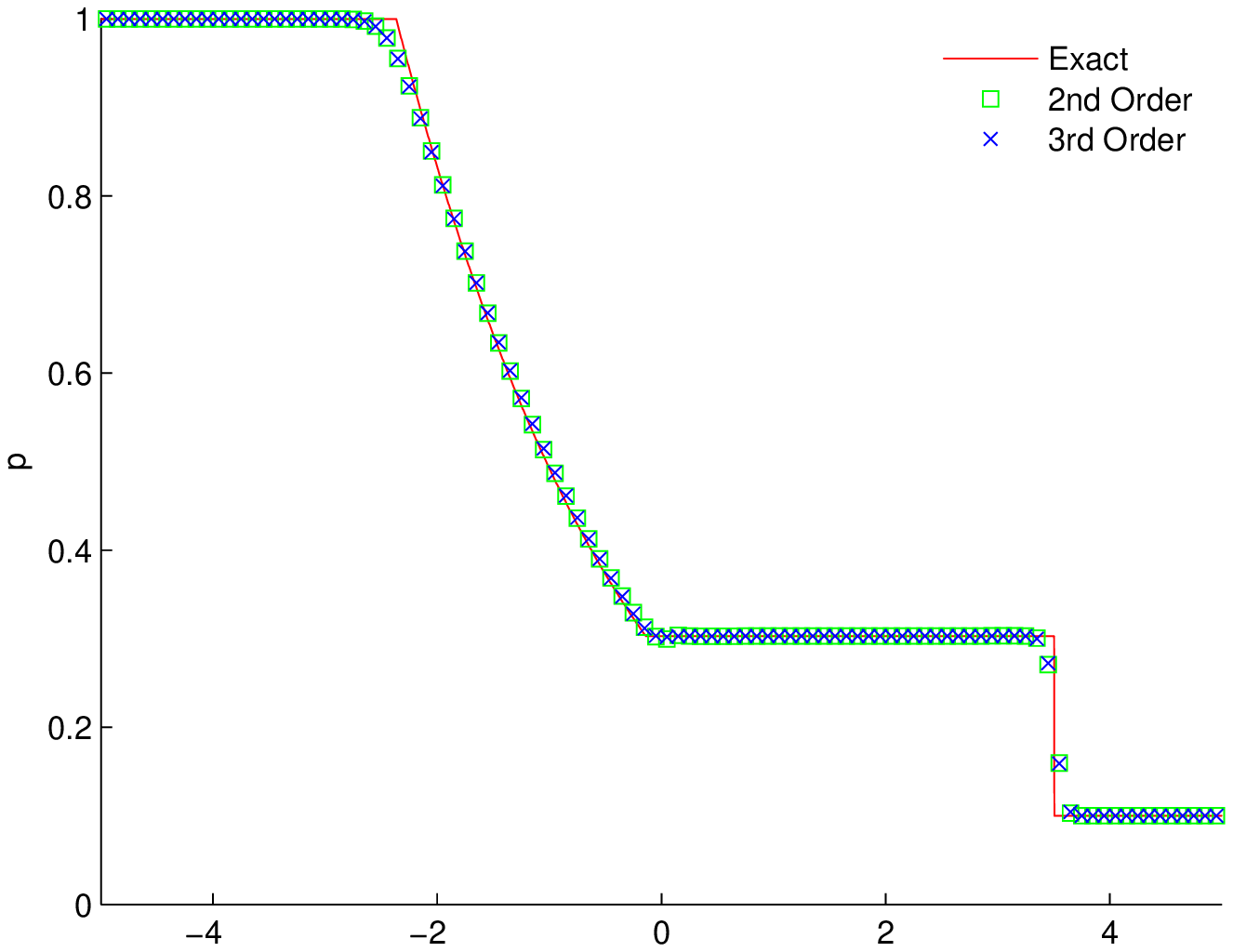}
    \includegraphics[height=2.0in, width=2.6in, trim=0 0 0 0, clip]{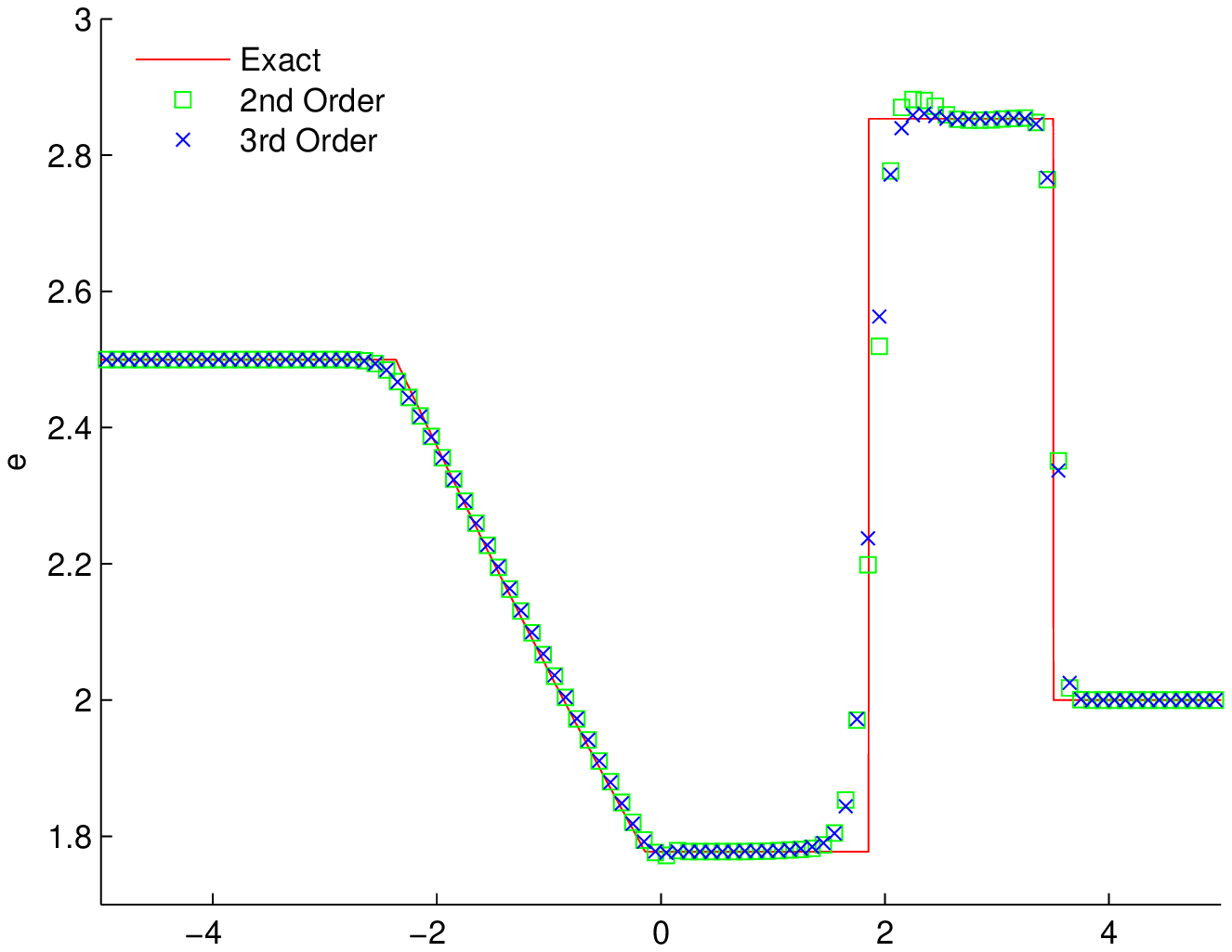}
  \end{center}
  \caption{ Numerical solutions of Sod problem: 100 grid point are used. }
  \label{fig:sod}
\end{figure}

The first test is the standard Riemann problem proposed by sod \cite{sod}.
The gas is
initially at rest with $\rho = 1$, $p = 1$ for $-5\leq x<0$ and $\rho = 0.125$, $p = 0.1$ for $0 \leq x < 5$.
At time  $t=2$, the numerical solutions with 100 points are shown in Fig. \ref{fig:sod}.
We can see that both of the computed solution agree well with the exact one and the third-order scheme shows better performance.

\subsection{123 problem}

\begin{figure}
  \begin{center}
    \includegraphics[height=2.0in, width=2.6in, trim=0 0 0 0, clip]{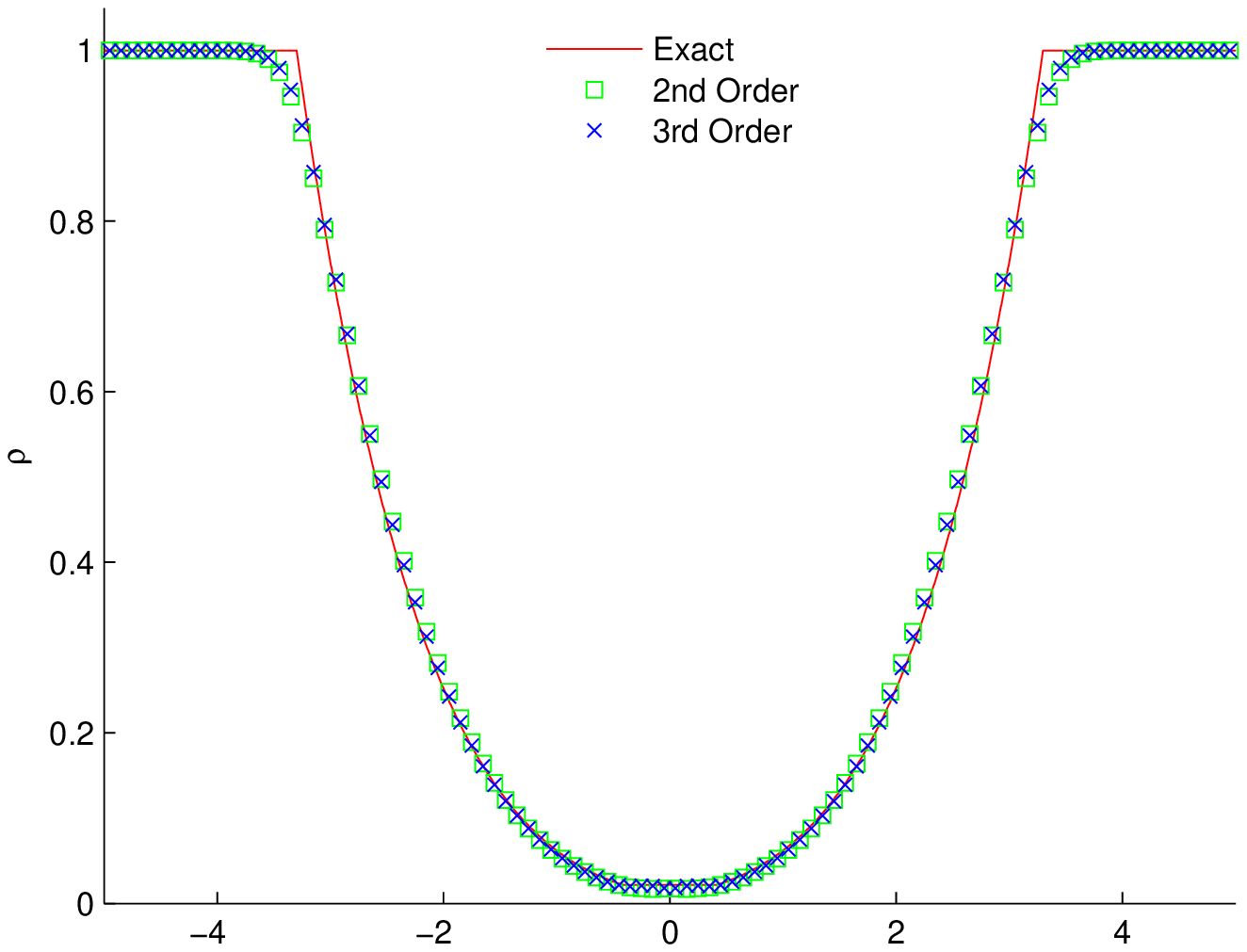}
    \includegraphics[height=2.0in, width=2.6in, trim=0 0 0 0, clip]{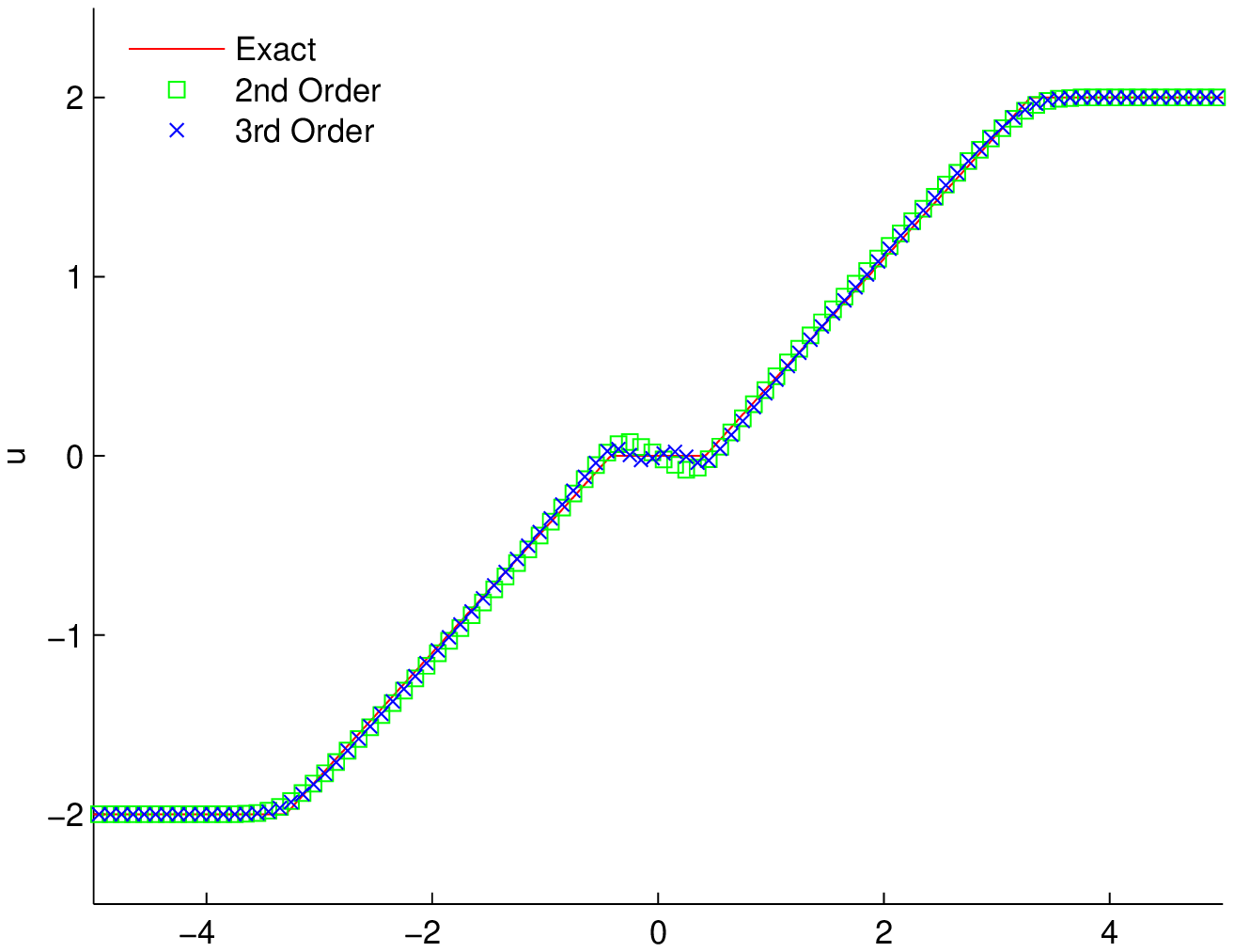}\\
    \includegraphics[height=2.0in, width=2.6in, trim=0 0 0 0, clip]{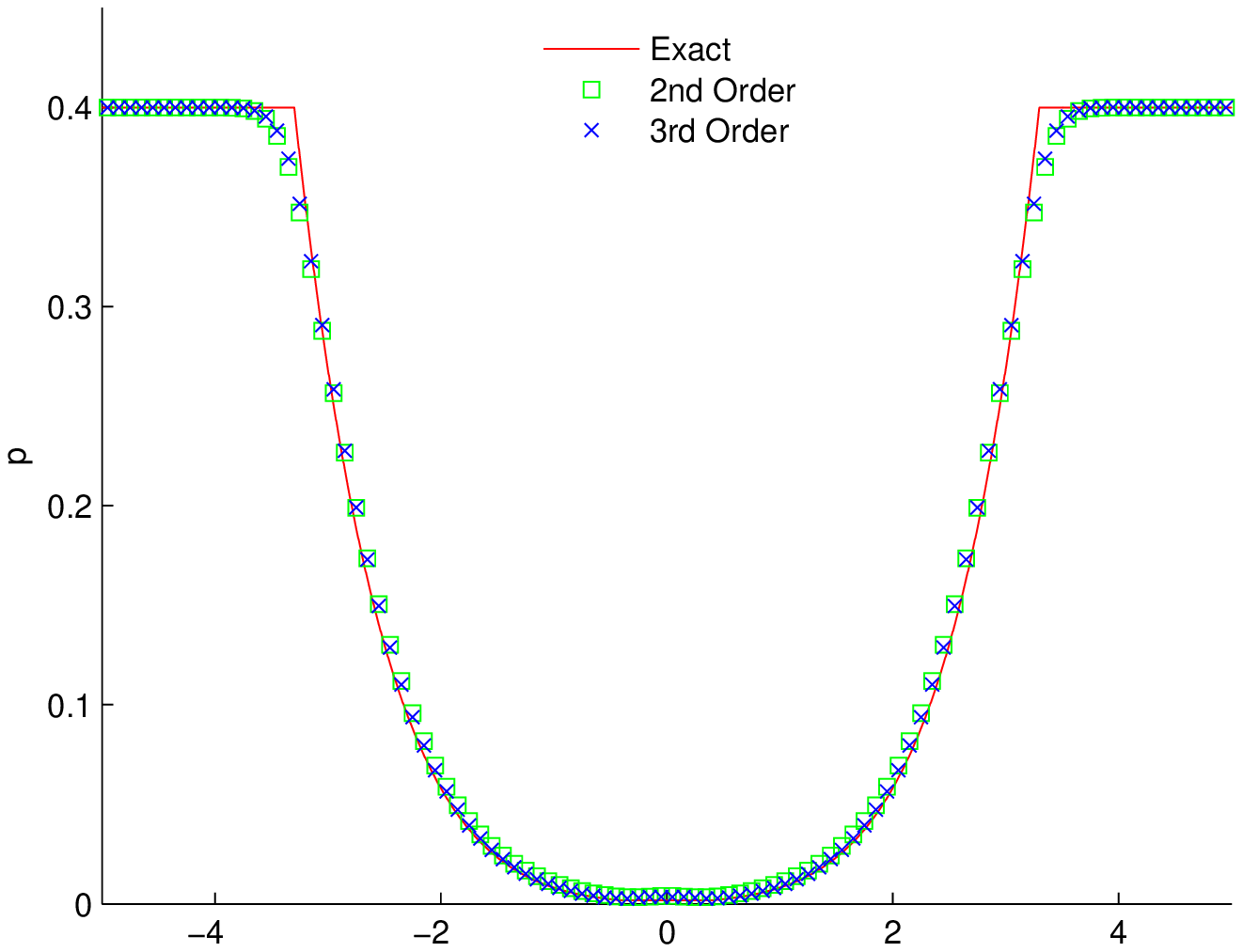}
    \includegraphics[height=2.0in, width=2.6in, trim=0 0 0 0, clip]{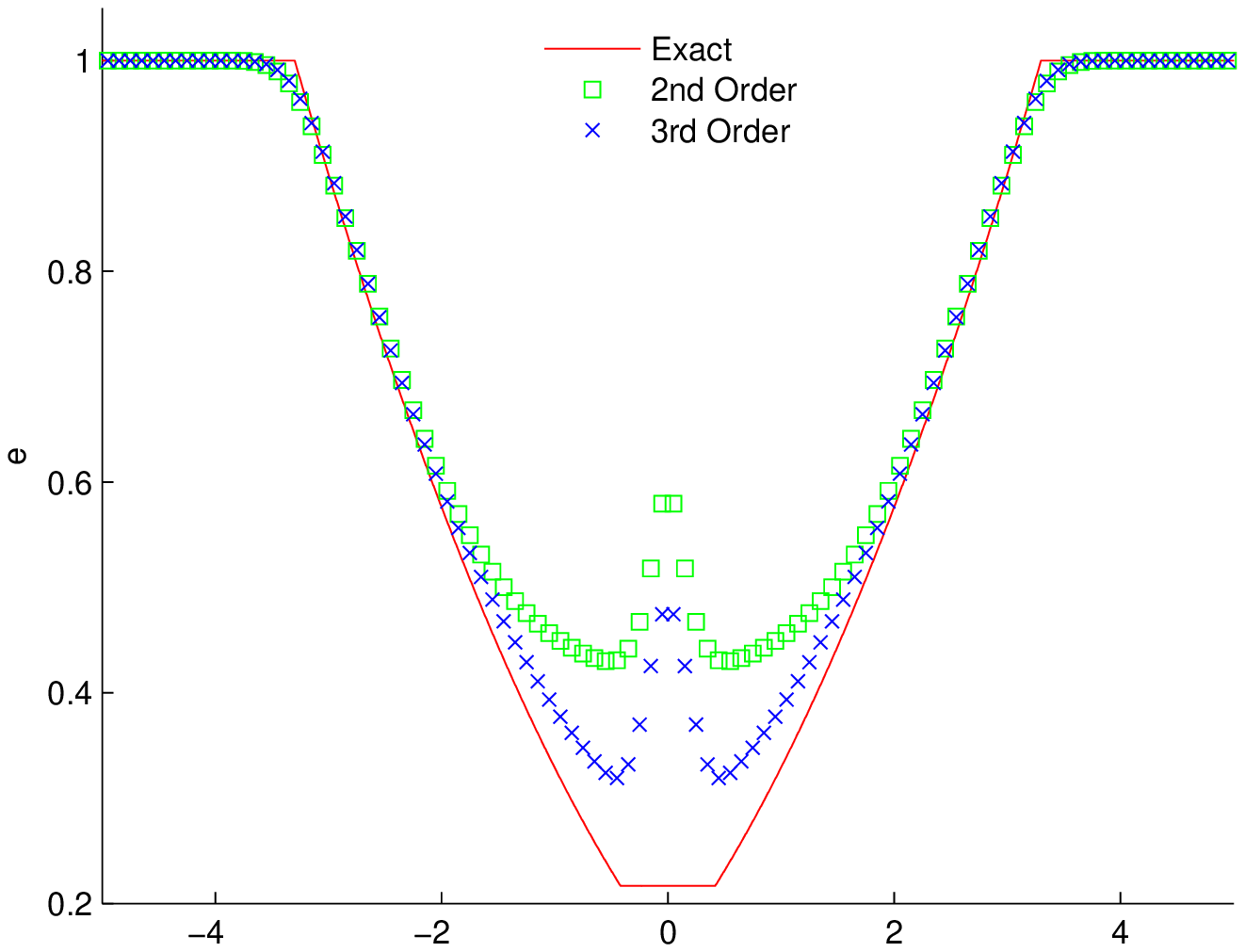}
  \end{center}
  \caption{ Numerical solutions of 123 problem. 100 grid point are used. }
  \label{fig:low_density}
\end{figure}

This example was first proposed by \cite{ein}.
The initial data is given with
$(\rho,u,p) = (1, 2, 0.4)$ for $-5 \leq x < 0$ and $(\rho,u,p) = (1,2,0.4)$ for $0 \leq x < 5$.
 The numerical solutions at time $t=1.2$ are shown in Fig. \ref{fig:low_density}. This test case demonstrates the ability of the GRP schemes to preserve the positivity of the density, pressure and internal energy. Again, the internal energy profile conforms the better performance of third-order scheme.

\subsection{Woodward-Colella blast wave problem}

\begin{figure}
  \begin{center}
    \includegraphics[height=2.0in, width=2.6in, trim=0 0 0 0, clip]{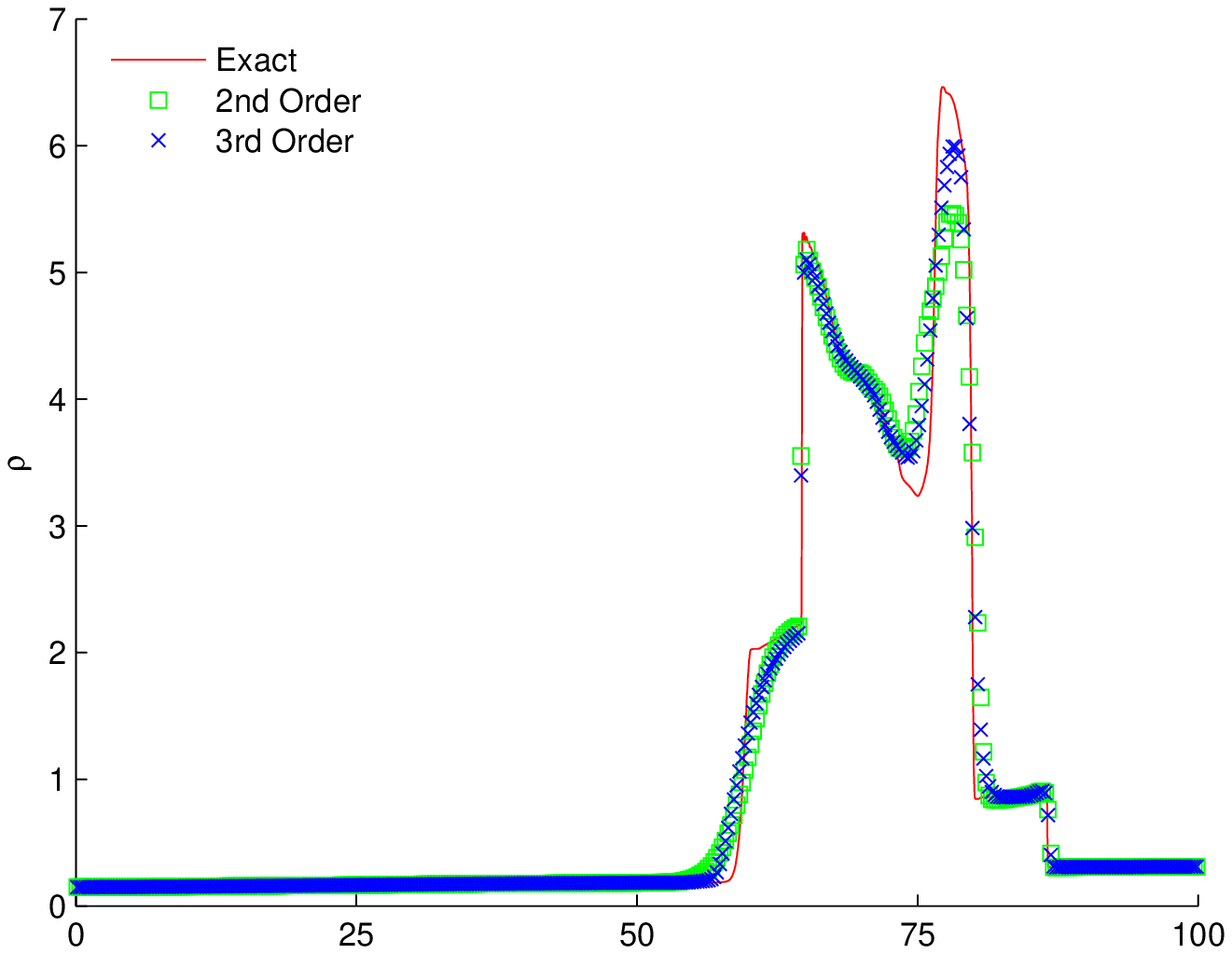}
    \includegraphics[height=2.0in, width=2.6in, trim=0 0 0 0, clip]{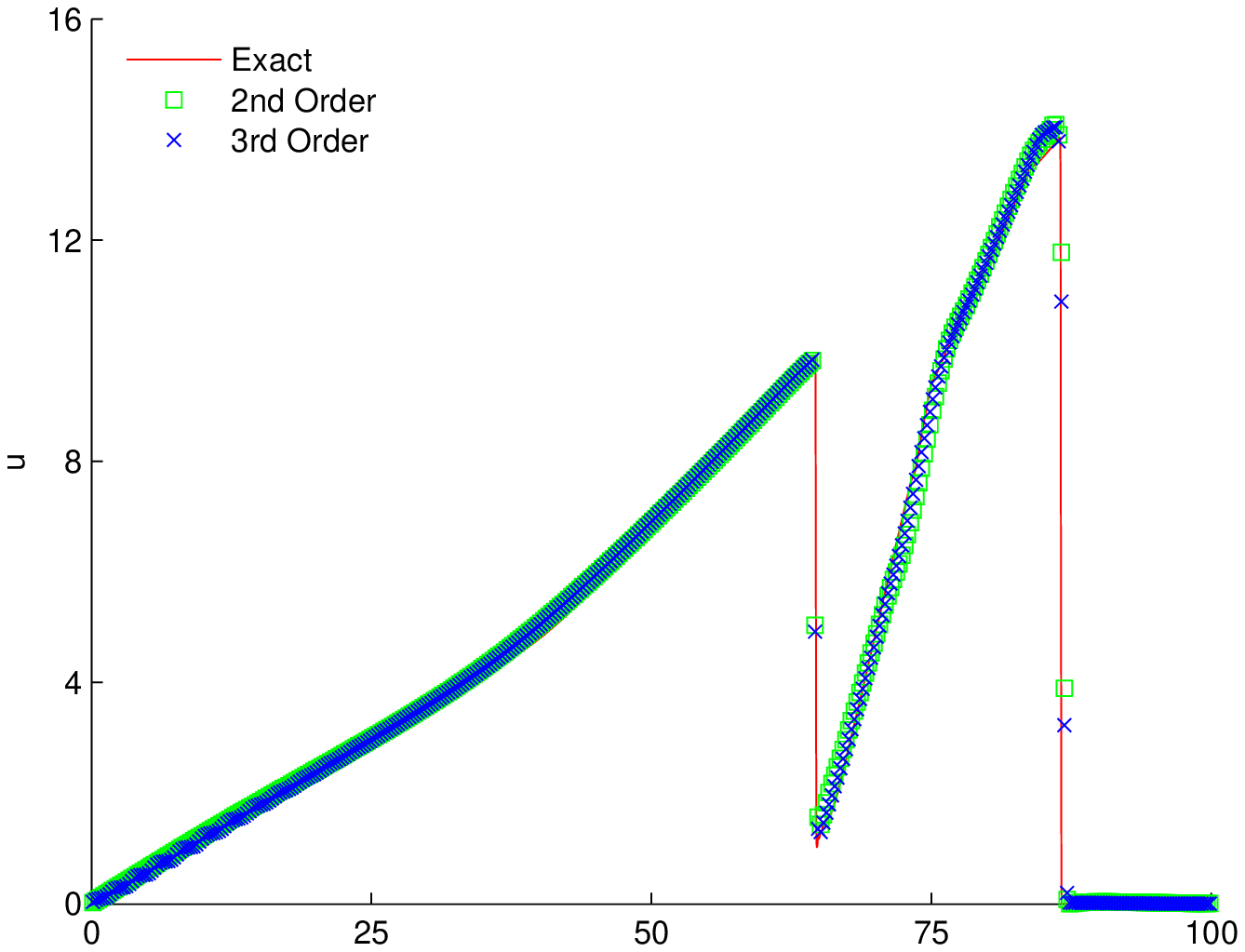}\\

    \includegraphics[height=2.0in, width=2.6in, trim=0 0 0 0, clip]{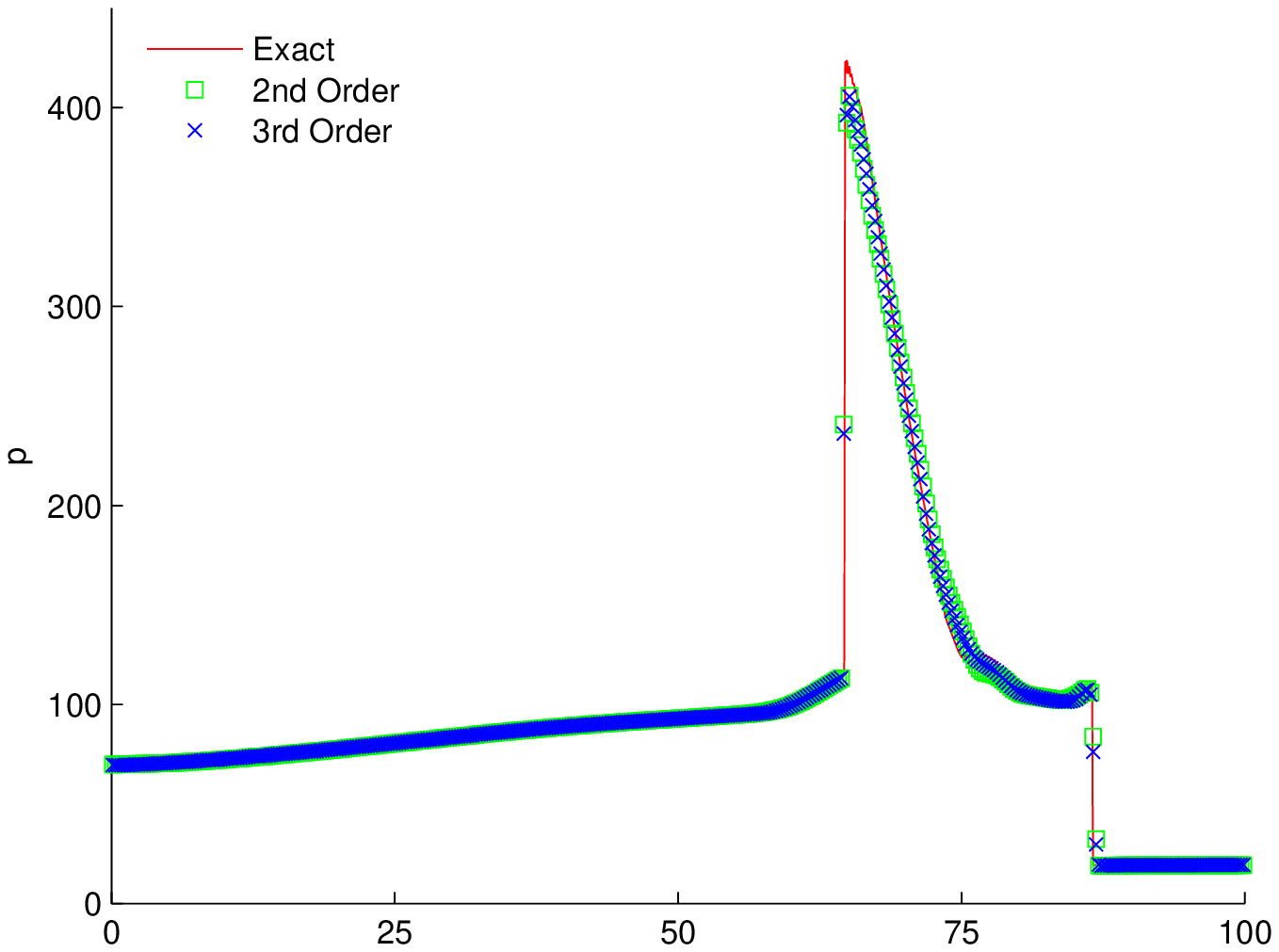}
    \includegraphics[height=2.0in, width=2.6in, trim=0 0 0 0, clip]{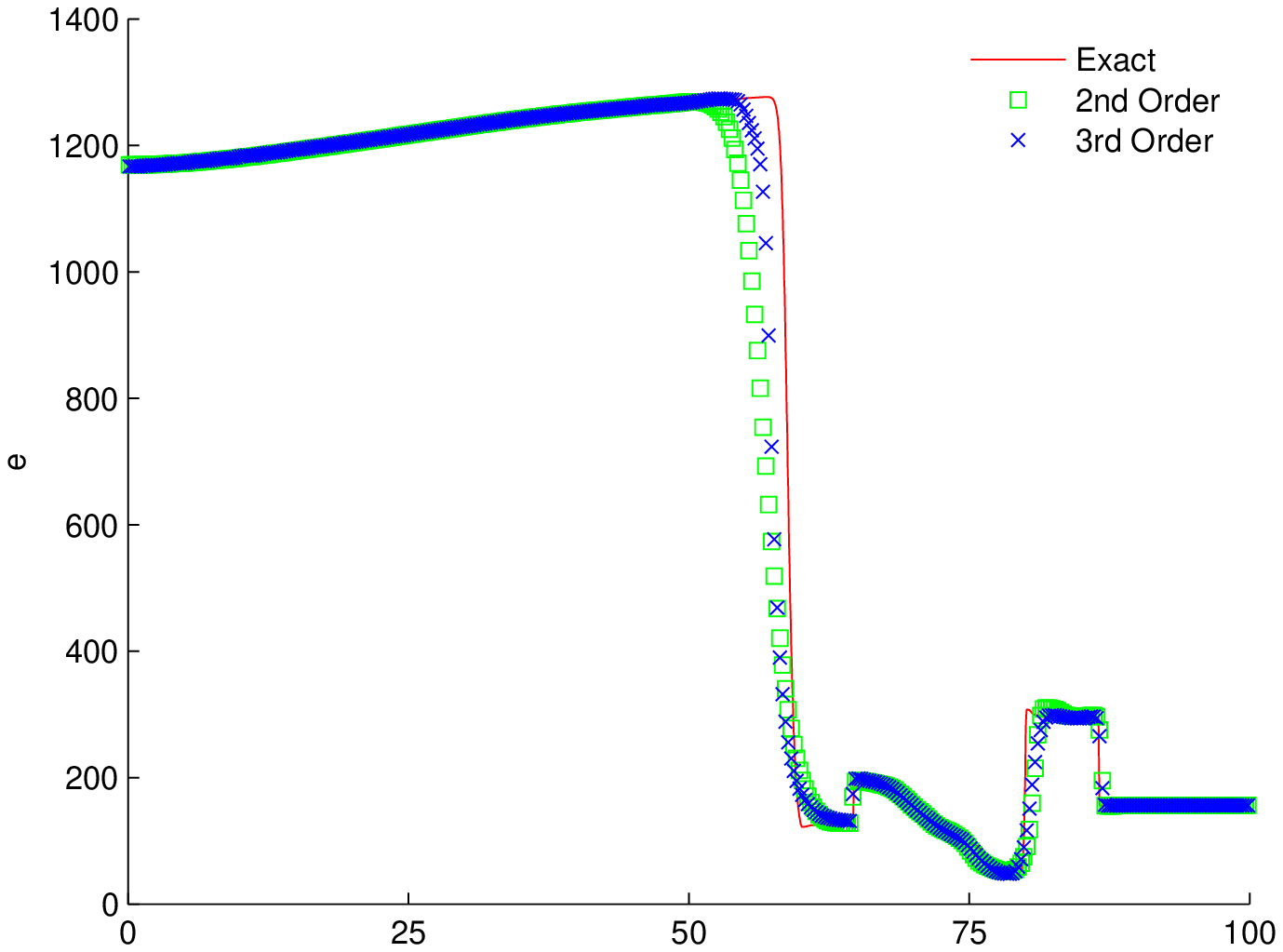}
 \end{center}
  \caption{ Numerical solutions of Woodward-Colella blast problem. 400 grid point are used.}
  \label{fig:blast_wave}
\end{figure}

This is a problem proposed by \cite{wood}.
The diatomic gas is initially at rest, and the density is unit
everywhere. The pressure is $p = 1000$ for $0 \leq x < 10$ and $p = 100$ for $90 \leq x < 100$, while it is only $p = 0.01$ in $10 \leq x < 90$. Reflecting boundary conditions are applied at both ends and
the output time is $t=3.8$.
Numerical solutions with 400 points are shown
in Fig. \ref{fig:blast_wave} to exhibit the performance of both schemes.
This test case clearly demonstrates the capability of both schemes
in the capturing of strong shock waves. The third-order scheme capture much sharper solution than the second-order scheme in the density and internal energy
distribution.

\subsection{Shock-density wave interaction}

\begin{figure}
  \begin{center}
    \includegraphics[height=2.5in, width=3.8in, trim=0 0 0 0, clip]{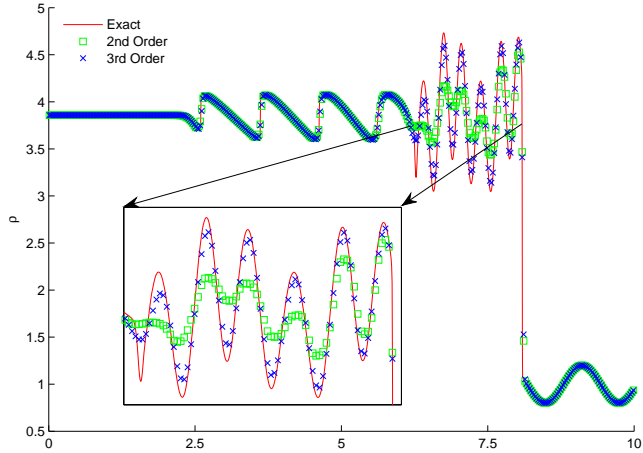}
  \end{center}
  \caption{ Numerical solutions of shock-density wave problem. 400 grid point are used.}
  \label{fig:shock-density}
\end{figure}

The Mach $3$ shock-entropy wave interaction \cite{shu2} is specified by the initial condition:
$(\rho, u, p)=(3.57134, 2.629369, 10.33333)$ for $0\leq x<1$ and $(\rho, u, p)=(1+0.2\sin(kx),0,1)$ for $1\leq x \leq 10$ with $k=5$.
The solution of this problem consists of a number of  shocklets and fine scales structures which are located behind a right-going main shock.
The computed density profile with 400 points, at $t = 2.0$, is shown in Fig. \ref{fig:shock-density}. Again the third-order scheme works better and captures much finer scale structures at high frequency waves behind the shock.

\subsection{Steady flow in a converging-diverging nozzle}

\begin{figure}
  \begin{center}
    \includegraphics[height=1.2in, width=4.0in, trim=0 0 0 0, clip]{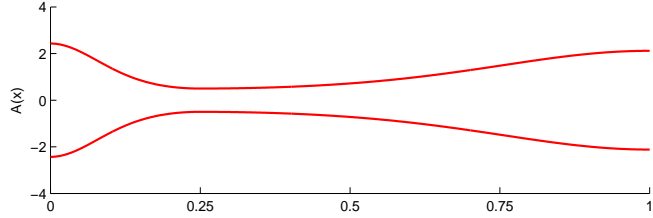}
  \end{center}
  \caption{ Nozzle contour. }
  \label{fig:A-x}
\end{figure}

We now use the examples in \cite[Sect. 6.5]{ben4} to test the ability of the GRP schemes to attain the steady state of a flow. Consider a flow in a converging-diverging nozzle, which occupies the internal $0 \leq x\leq 1$ and has a smooth cross-sectional area function $ A(x)$ given by the following expression:
\begin{equation}\label{eq:A-x}
A(x)=\begin{cases} A_{\textrm{in}} \exp\big( -\log(A_{\textrm{in}}) \sin^2(2 \pi x) \big), & 0\leq x<0.25;\\
A_{\textrm{ex}} \exp\l( -\log(A_{\textrm{ex}}) \sin^2\l( \f{2 \pi(1- x)}{3} \r) \r), & 0.25\leq x \leq 1,
\end{cases}
\end{equation}
where $A_{\textrm{in}}=4.8643$ and $A_{\textrm{ex}}=4.2346$. See Fig. \ref{fig:A-x}.
For a steady duct flow of a perfect gas, the Mach number $M(x)=u(x)/c(x)$ is determined by $A(x)$ through the algebraic relation
\begin{equation}\label{eq:M-x}
[A(x)]^2=\f{1}{[M(x)]^2} \l[ \f{2}{\ga+1} \l(1+\f{\ga-1}{2}[M(x)]^2 \r)\r]^{\f{\ga+1}{\ga-1}}.
\end{equation}
Then the steady flow profiles in the nozzle are given by
\begin{equation}\label{eq:duct-solu}
\begin{split}
&p(x)=p_0\l( 1+\f{\ga-1}{2} [M(x)]^2\r)^{-\f{\ga}{\ga-1}},\\
&\rho(x)=\rho_0\l( 1+\f{\ga-1}{2} [M(x)]^2\r)^{-\f{1}{\ga-1}},\\
&u(x)=M(x) \sqrt( \ga p(x)/ \rho(x)),
\end{split}
\end{equation}
for the flow being smooth, where $\rho_0$ and $p_0$ need to be specified.

\begin{figure}
  \begin{center}
    \includegraphics[height=1.1in, width=4.0in, trim=0 0 0 0, clip]{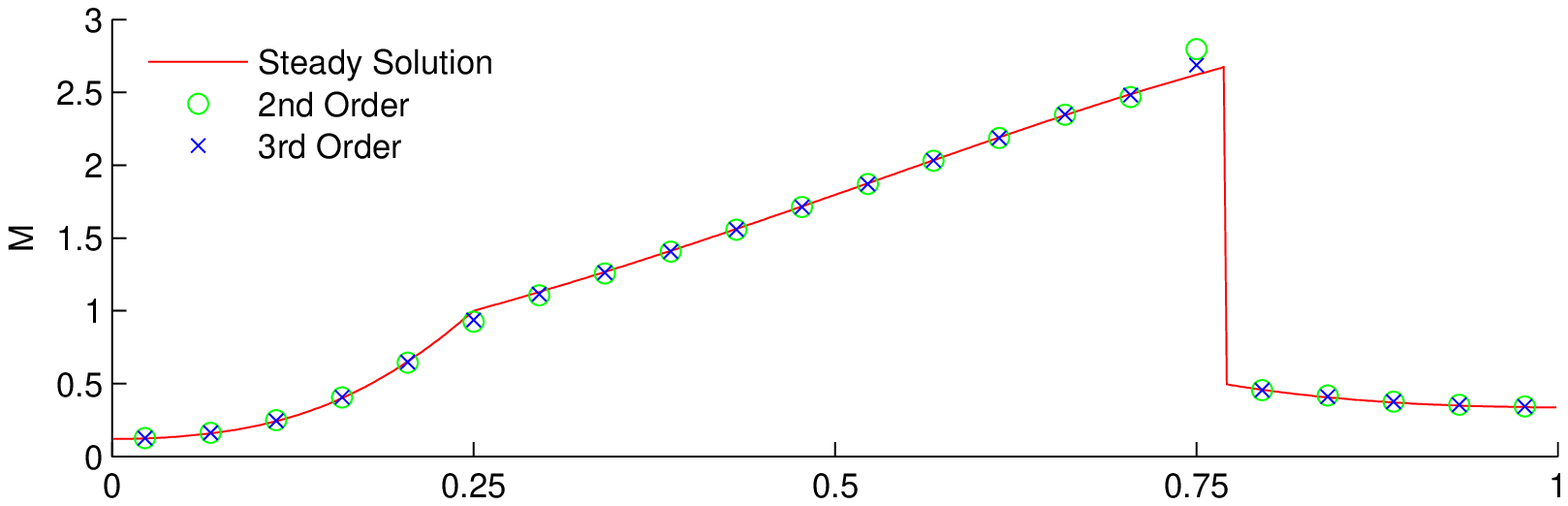}\\
    \includegraphics[height=1.1in, width=4.0in, trim=0 0 0 0, clip]{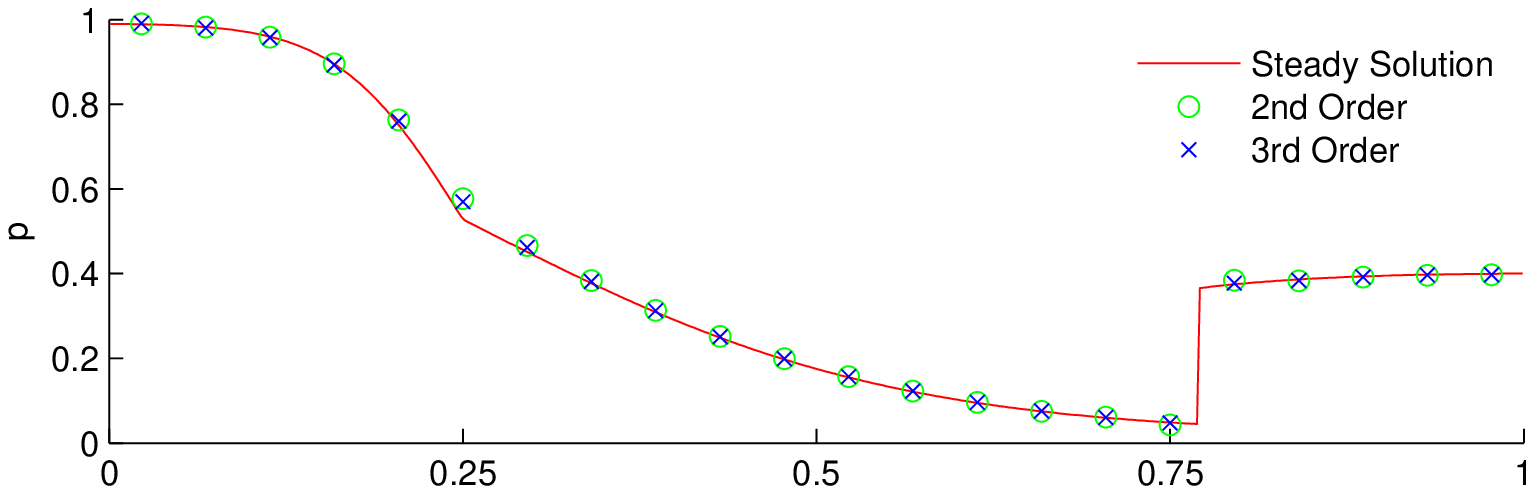}
  \end{center}
  \caption{  Large time flow in Laval nozzle: Case B. 22 grid points are used, at time $t=15.5$. }
  \label{fig:B15.5}
\end{figure}

The initial data we use are
\begin{equation}\label{eq:duct-ini}
U(x,0)=\begin{cases}
U_L=(\rho_0,0,p_b),& 0<x<0.25,\\
U_R=(\rho_0(p_b/p_0)^{1/\ga},0,p_b ),& 0.25<x<1,
\end{cases}
\end{equation}
where $p_b$ is a constant determined by the steady state solution at $x=1$.
We consider two cases. In both cases we take $\rho_0=p_0=0$ and $A(x)$ as in (\ref{eq:A-x}).
\begin{enumerate}
\item[(A)]  A smooth flow where $p(1)=0.0272237$ is obtained from (\ref{eq:duct-solu}) by taking $x=1$ in (\ref{eq:M-x}), leading to $M(1)=3$.
\item[(B)] Setting $p(1)=0.4$ leads to a discontinuous steady state solution, as shown by solid lines in Fig. \ref{fig:B15.5}.
\end{enumerate}

\begin{figure}
  \begin{center}
    \includegraphics[height=1.0in, width=4.0in, trim=0 0 0 0, clip]{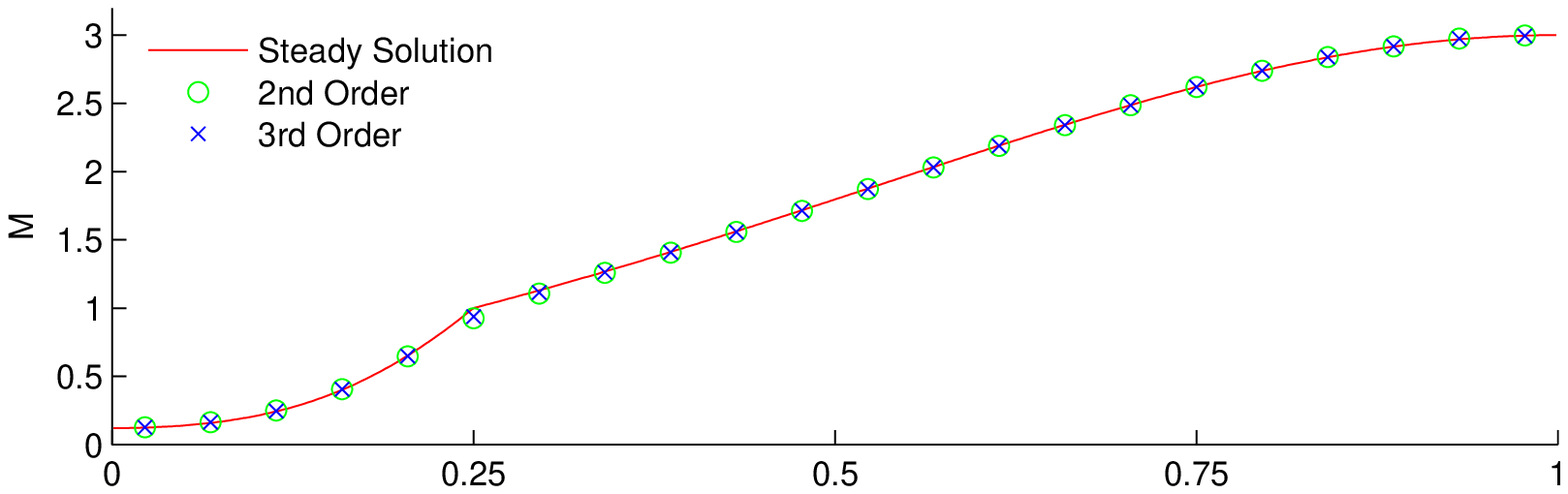}\\
    \includegraphics[height=1.0in, width=4.0in, trim=0 0 0 0, clip]{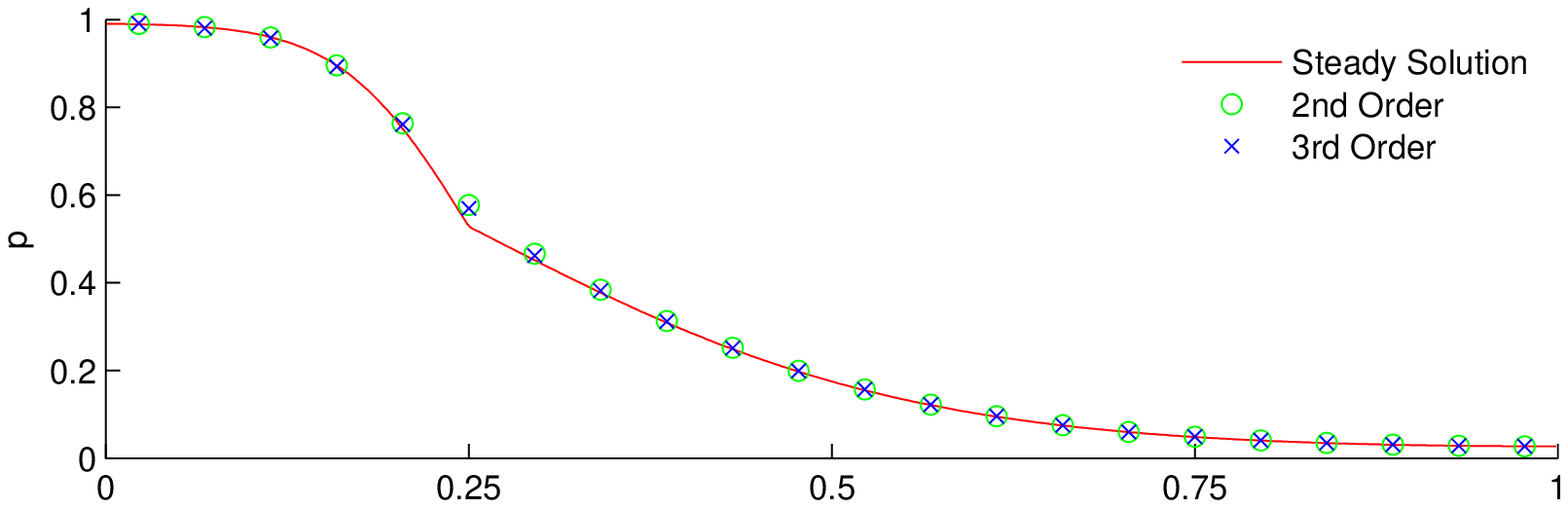}
  \end{center}
  \caption{ Large time flow in Laval nozzle: Case A. 22 grid points are used, at time $t=15.5$. }
  \label{fig:A15.5}
\end{figure}

\begin{figure}
  \begin{center}
    \includegraphics[height=1.1in, width=4.0in, trim=0 0 0 0, clip]{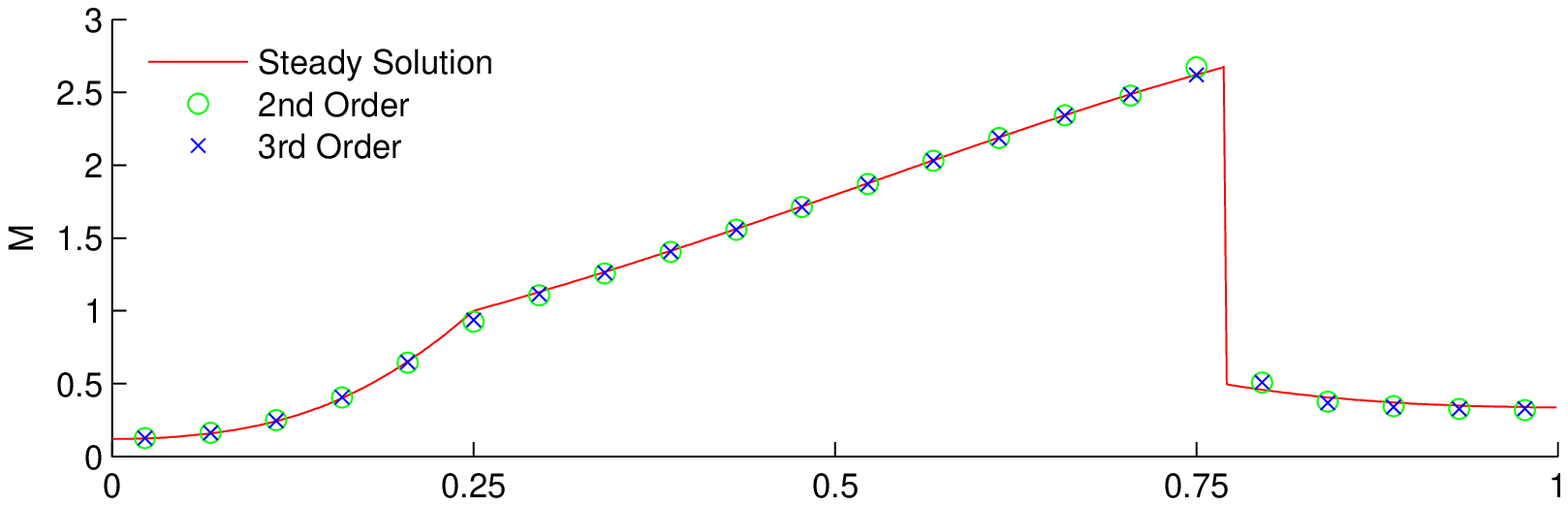}\\
    \includegraphics[height=1.1in, width=4.0in, trim=0 0 0 0, clip]{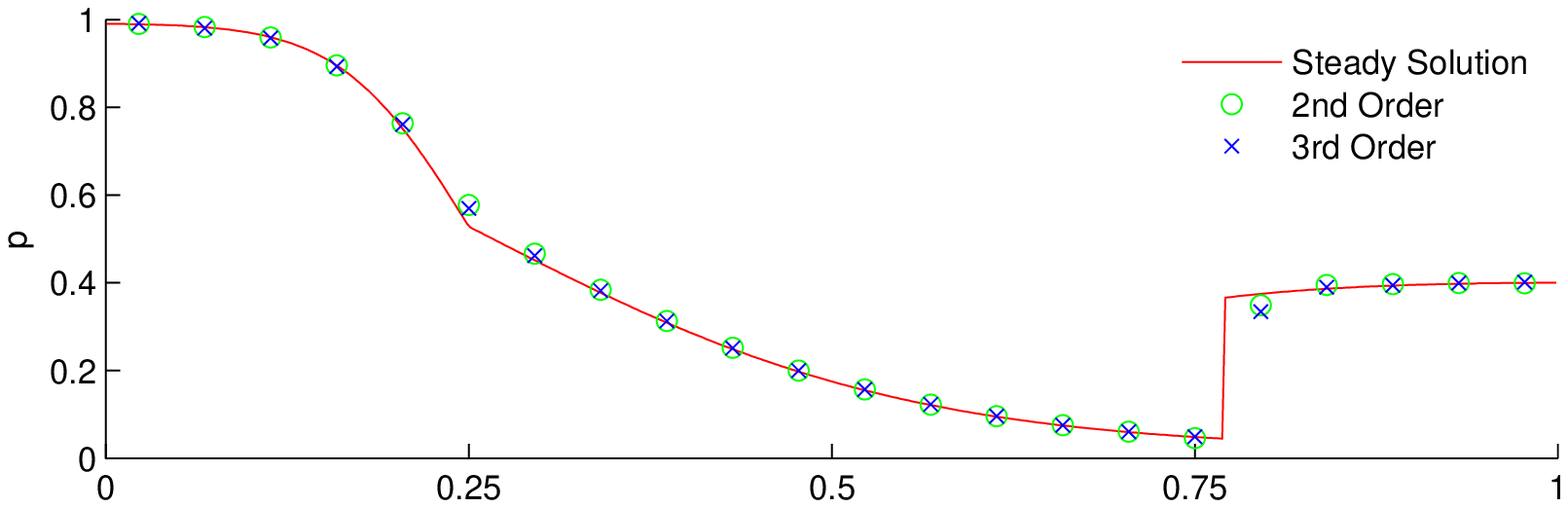}
  \end{center}
  \caption{  Large time flow in Laval nozzle: Case B. 22 grid points are used, at time $t=2.5$. }
  \label{fig:B2.5}
\end{figure}

We use the strategy in \cite[Sect. 6.5]{ben4} to deal with the boundary conditions at $x=0$ and $1$. In both cases, the number of grid points used are 22.
As shown in Fig. \ref{fig:A15.5} and \ref{fig:B15.5}, both of the GRP solutions at $t=15.5$ are good agreement with the exact solution. The third-order GRP solutions are closer to the analytical solutions than the second-order ones.
Moreover, as shown in Fig. \ref{fig:B2.5},  the GRP solutions almost attain the steady state at time $t=2.5$. This shows that the GRP solutions converges to steady solution quickly.

\appendix

\section{Formulae in Section \ref{sec:solver-form}}
\label{app:1}
For the general case of $\gamma>1$, the L-equations of $\v{w}_-=(S, \psi)$ in Section \ref{sec:solver-form} yield
\begin{equation}\label{eq:euler-l-g}
\begin{split}
& D_{\la_-} S(0,\be) = A_1 (\psi_L -\be )^{\f{\gamma+1}{\gamma-1}},\\
& D_{\la_-} \psi(0,\be)=A_2 (\psi_L -\be )^{\f{\gamma+1}{2(\gamma-1)}}+A_3(\be) (\psi_L -\be )^{\f{2\gamma}{\gamma-1}}+Z_0(\be) (\psi_L -\be )^{\f{\gamma+1}{2(\gamma-1)}},
\end{split}
\end{equation}
where $A_2=\tilde{A}_{2}-Z_0(\be_L)$ and
\begin{equation*}
Z_0(\be)=\begin{cases}
B_1(\psi_L -\be )^{\f{\gamma-3}{2(\gamma-1)}} + B_2 (\psi_L -\be )^{\f{3\gamma-5}{2(\gamma-1)}}, & \textrm{if}\ \gamma\neq 3, 5/3;\\[2mm]
 \f{\psi_L}{2} \f{A'(x)}{A(x)} \ln(\psi_L -\be)+B_2 (\psi_L -\be )^{\f{3\gamma-5}{2(\gamma-1)}}, & \textrm{if}\ \gamma=3;\\[2mm]
B_1 \psi_L(\psi_L -\be )^{\f{\gamma-3}{2(\gamma-1)}}  -\f{3}{8} \f{A'(x)}{A(x)} \ln(\psi_L -\be ), & \textrm{if}\ \gamma= 5/3,
\end{cases}
\end{equation*}
\begin{align*}
&A_1=(\psi_L-\be_L)^{-\f{\gamma+1}{\gamma-1}} D_{\la_-} S(0,\be_L),\\
&\tilde{A}_2=-\f{1}{\gamma(3\gamma-1) S_L} (\psi_L-\be_L)^{\f{\gamma-3}{2(\gamma-1)}} D_{\la_-} S(\be_L) + (\psi_L-\be_L)^{-\f{\gamma+1}{2(\gamma-1)}} D_{\la_-} \psi(\be_L).
\end{align*}

The coefficients $A_i$, $B_j$, $C_k$ are defined as in Table \ref{tab:AB} and \ref{tab:C}.

\begin{table}
\caption{\label{tab:AB} The coefficients $A_3,\cdots, A_{19}$ and $B_1,\cdots, B_{12}$ \hfill  }
\centering
\begin{tabular}{lclclcl}
\hline
$A_3$ & \hspace{0.8cm} & $\f{1}{\gamma(3\gamma-1) S_L} A_1$ & \hspace{0.8cm} & $B_1$ & \hspace{0.8cm} & $\l(\f{\ga-1}{\ga-3}\r) \psi_L \f{A'(x)}{A(x)}$ \\

$A_4$ &   &  $-\f{1}{\gamma(\ga+1) S_L} A_1$ &  &
$B_2$ &   & $-\f{2(\ga-1)}{(\ga+1)(3\ga-5)}\f{A'(x)}{A(x)}$\\

$A_5$ &   & $\f{\ga-1}{4} A_2$ &
& $B_3$ &  & $\f{\ga-1}{\gamma+1} \f{A'(x)}{A(x)}$  \\

$A_6$ &   & $\f{\ga-1}{4}( A_3- A_4)$ &
& $B_4$ &  &  $-\f{2(\ga-1)}{(\gamma+1)^2} \f{A'(x)}{A(x)}$ \\

$A_7$ &   & $-\f{(3-\ga)(\ga+1)}{8(\ga-1)} A_2$ &
 & $B_5$ &  & $\f{\ga-1}{4} (B_1-B_3)$  \\

$A_8$ &   & $-\f{2\ga}{\ga-1}\l(\f{3-\ga}{4} A_3 +\f{1+\ga}{4} A_4\r) $  &
& $B_6$ &  & $\f{\ga-1}{4}( B_2- B_4)$  \\

$A_{9}$ &  & $\f{1}{\ga(\ga-1)S_L^2} A_1$ &
& $B_7$ &  & $-\f{3-\ga}{4} B_1-\f{\ga+1}{4} B_3$   \\

$A_{10}$ &  & $-\f{1}{2\ga(\ga-1) S_L} A_7$ &
&  $B_8$ &  & $-\f{3-\ga}{2} B_2-\f{\ga+1}{2} B_4$  \\

$A_{11}$ &  & $-\f{1}{2\ga(\ga-1) S_L} A_8$ &
 & $B_9$ &  & $\f{1}{2}\l(\f{\ga+1}{\ga-1}\r)^2B_5-\f{\ga+1}{2(\ga-1)}B_7 $ \\

$A_{12}$ &  & $ \l(\f{\ga+1}{\ga-1} \r)^2 A_5-\f{\ga+1}{2(\ga-1)}A_7$ &
 &  $B_{10}$ & & $\f{1}{2}\l(\f{\ga+1}{\ga-1}\r)^2B_6-\f{\ga+1}{2(\ga-1)}B_8 $\\

$A_{13}$ &  & $ \l(\f{\ga+1}{\ga-1} \r)^2 A_6-\f{\ga+1}{2(\ga-1)}A_8$ &
& $B_{11}$ &  & $-\f{1}{2\ga(\ga-1) S_L} B_7$\\

$A_{14}$ &  & $ A_9+A_{11}$ &
 & $B_{12}$ &  & $-\f{1}{2\ga(\ga-1) S_L} B_8$ \\

$A_{15}$ &  & $2A_{1} A_{12}$ \\
$A_{16}$ &  & $2A_{1} A_{13}$ \\
$A_{17}$ &  & $A_1 A_{10}+A_2 A_{13}+A_3 A_{12} $ \\
$A_{18}$ &  & $A_1 A_{14}+A_3 A_{13} $ \\
$A_{19}$ &  & $A_2 A_{12}$ \\
\hline
\end{tabular}
\end{table}
\begin{table}
\caption{\label{tab:C} The coefficients $C_1,\cdots, C_{9}$ \hfill  }
\centering
\begin{tabular}{lcl}
\hline
$C_1$ & \hspace{2cm} & $ 2A_1B_9 $ \\
$C_2$ & \hspace{2cm} & $ 2A_1B_{10} $ \\
$C_3$ & \hspace{2cm} & $ A_2B_{9}+A_{12}B_{1}-\f{\psi_L}{2} A_7 \f{A'(x)}{A(x)}  $ \\
$C_4$ & \hspace{2cm} & $ A_3B_{9}+A_{1}B_{11}+A_{13}B_{1} -\f{\psi_L}{2} A_8 \f{A'(x)}{A(x)}  $ \\
$C_5$ & \hspace{2cm} & $ A_2B_{10}+A_{12}B_{2}+\l(\f{A_7}{\ga+1} -\f{A_2}{2}\r) \f{A'(x)}{A(x)}  $ \\
$C_6$ & \hspace{2cm} & $ A_3B_{10}+A_{1}B_{12}+A_{13}B_{2} +\l(\f{A_8}{\ga+1} -\f{A_3+A_4}{2}\r) \f{A'(x)}{A(x)}  $ \\
$C_7$ & \hspace{2cm} & $ B_1B_{9} -\f{\psi_L}{2} B_7 \f{A'(x)}{A(x)} - \psi_L^2 \l( \f{A'(x)}{A(x)} \r)^{'}$ \\
$C_8$ & \hspace{2cm} &
$ B_2B_{9}+B_1 B_{10}+\l( -\f{\psi_L}{2} B_8 +\f{B_7}{\ga+1} -\f{B_1+B_3}{2} \r) \f{A'(x)}{A(x)} + \f{\ga+3}{\ga+1}\psi_L\l( \f{A'(x)}{A(x)} \r)^{'}  $\\
$C_9$ & \hspace{2cm} & $ B_2B_{10}+\l( -\f{ B_8}{\ga+1} -\f{B_2+B_4}{2} \r) \f{A'(x)}{A(x)}-\f{2}{\ga+1} \l( \f{A'(x)}{A(x)} \r)^{'}$  \\
\hline
\end{tabular}
\end{table}

The function $Z_1(\be)$ and $Z_2(\be)$ in Proposition \ref{prop:euler-lq} are as follows
\begin{align*}
&Z_1(\be)=\f{2(\ga-1)}{3\ga-1} A_{15} (\psi_L-\be)^{\f{-3\ga+1}{2(\ga-1)}}-A_{16} \ln (\psi_L-\be)\\
&\hspace{2cm}+\f{\gamma-1}{\gamma+1} C_1 (\psi_L-\be)^{-\l(\f{\ga+1}{\ga-1}\r)}+\f{\gamma-1}{2} C_2 (\psi_L-\be)^{\f{-2}{\ga-1}} ,\\[3mm]
& Z_2(\be)=\f{1}{2\ga^2 S_L} (\psi_L-\be)^{-\f{2}{\ga-1}} D_{\la_-}^2 S(\be)+\f{\ga-1}{\ga+1}\l[\f{A_{15}}{\ga^2 S_L} -2 A_{17}\r](\psi_L-\be)^{\f{\ga+1}{2(\ga-1)}}\\
&\hspace{2cm}+ \f{\ga-1}{2 \ga} \l[ \f{A_{16}}{2 \ga^2S_L} - A_{18}\r](\psi_L-\be)^{\f{2\ga}{\ga-1}}+A_{19} (\psi_L-\be)^{-1}\\
&\hspace{2cm}+\l[\f{C_1 }{2\ga^2 S_L}-C_4 \r] (\psi_L-\be)+\l[ \f{C_2 }{4\ga^2 S_L} -\f{C_6}{2}\r] (\psi_L-\be)^2\\
&\hspace{2cm}+\f{2(\ga-1)}{\ga+1} C_3(\psi_L-\be)^{-\f{\ga+1}{2(\ga-1)}}
-\f{2(\ga-1)}{ \ga-3}C_5(\psi_L-\be)^{\f{\ga-3}{2(\ga-1)}}\\
&\hspace{2cm}+\f{\ga-1}{2} C_7(\psi_L-\be)^{-\f{2}{\ga-1}}-\f{\ga-1}{\ga-3} C_8(\psi_L-\be)^{\f{\ga-3}{\ga-1}}\\
&\hspace{2cm}-\f{\ga-1}{2(\ga-2)}C_9(\psi_L-\be)^{\f{2(\ga-2)}{\ga-1}}.
\end{align*}



\bibliographystyle{elsarticle-num}
\bibliography{<your-bib-database>}



\end{document}